\documentclass[a4paper,leqno]{article}
\usepackage{amsfonts}
\usepackage{amsmath}
\usepackage{amssymb}
\numberwithin{equation}{section}
\begin{document}
\newtheorem{theorem}{Theorem}[section]
\newtheorem{lemma}[theorem]{Lemma}
\newtheorem{corollary}[theorem]{Corollary}
\newtheorem{conjecture}[theorem]{Conjecture}
\newtheorem{remark}[theorem]{Remark}
\newtheorem{definition}[theorem]{Definition}
\newtheorem{problem}[theorem]{Problem}
\newtheorem{example}[theorem]{Example}
\newtheorem{proposition}[theorem]{Proposition}
\title{{\bf Canonical singular hermitian metrics \\ on relative log canonical bundles}}
\date{September 27, 2010}
\author{Hajime TSUJI}
\maketitle
\begin{abstract}
\noindent We introduce a new class of canonical analytic Zariski decompositions (AZD's in short)  called the supercanonical AZD's  on the canonical bundles of smooth projective varieties with pseudoeffective canonical classes.  
We study the variation of the supercanonical AZD $\hat{h}_{can}$ under 
projective deformations and give a new proof of the invariance of plurigenera. 
Moreover extending the results to the case of KLT pairs, we prove the invariance of logarithmic plurigenera for a family of 
KLT pairs. This paper supersedes \cite{tu8,tu}.  \\ 
MSC: 14J15,14J40, 32J18
\end{abstract}
\tableofcontents

\section{Introduction}

Let $X$ be a smooth projective variety and let $K_{X}$ be the canonical bundle 
of $X$.  
In algebraic geometry, the canonical ring $R(X,K_{X}) : = \oplus_{m=0}^{\infty}\Gamma\left(X,{\cal O}_{X}(mK_{X})\right)$ is one of the main objects to study. 
And it has been studied the variation of pluricanonical systems in terms of 
variation of Hodge structures(\cite{f,ka1,v1,v}).

The  purpose of this article is to study the variation of (log) canonical rings on a projective family by introducing a canonical singular hermitian metric on the relative (log) canonical bundles.  The important feature here is the semipositivity of  the relative (log) canonical bundles and the invariance of (log) plurigenera is obtained at the same time. Moreover we can deal with the adjoint line bundle of a pseudoeffective $\mathbb{Q}$-line bundle in a systematic way.   
\vspace{0mm} \\  
   
Let $X$ be a smooth projective variety such that $K_{X}$ is pseudoeffective.
In this article, we construct a singular hermitian metric $\hat{h}_{can}$ on $K_{X}$ such that 
\begin{enumerate}
\item[(1)] $\hat{h}_{can}$ is uniquely determined by $X$,    
\item[(2)] The curvature current $\sqrt{-1}\,\Theta_{\hat{h}_{can}}$ is semipositive,
\item[(3)] $H^{0}(X,{\cal O}_{X}(mK_{X})\otimes {\cal I}(\hat{h}_{can}^{m}))
\simeq H^{0}(X,{\cal O}_{X}(mK_{X}))$ holds for every $m\geqq 0$,  
\end{enumerate}
where 
${\cal I}(\hat{h}_{can}^{m})$ denotes the multiplier ideal sheaf of 
$\hat{h}_{can}^{m}$ as is defined in \cite{n}. 
We may summerize the 2nd and the 3rd conditions by introducing the following notion. 
  
\begin{definition}{\em {\bf (AZD)}}\label{azd}(\cite{tu,tu2})
Let $M$ be a compact complex manifold and let $L$ be a holomorphic line bundle
on $M$.  A singular hermitian metric $h$ on $L$ is said to be 
an analytic Zariski decomposition (AZD in short), if the followings hold.
\begin{enumerate}
\item[(1)] The curvature current $\sqrt{-1}\,\Theta_{h}$ is semipositive. 
\item[(2)] For every $m\geq 0$, the natural inclusion
\[
H^{0}(M,{\cal O}_{M}(mL)\otimes{\cal I}(h^{m}))\rightarrow
H^{0}(M,{\cal O}_{M}(mL))
\]
is an isomorphim.  \fbox{}  
\end{enumerate} 

\begin{remark}\label{rm}
A line bundle $L$ on a projective manifold  $X$ admits an AZD, if and only if 
$L$ is pseudoeffective (\cite[Theorem 1.5]{d-p-s}). \fbox{}
\end{remark}
\end{definition} 

\noindent In this sense, we construct an AZD $\hat{h}_{can}$ on $K_{X}$  depending only on $X$, 
when $K_{X}$ is pseudoeffective (by Remark \ref{rm} this is the minimal requirement for the existence of an AZD). 
In fact $\hat{h}_{can}$ is not only an AZD of $K_{X}$, but also a singular hermitian metric with minimal singularities on $K_{X}$ (cf. Definition \ref{minAZD}).
The important feature of this canonical metric $\hat{h}_{can}$ is that 
it naturally defines a singular hermitian metric on the relative canonical bundle on a smooth projective family of smooth projective varieties with pseudoeffective canonical bundles just by assigning the canonical metric 
on each smooth fiber  and taking lower-semicontinuous envelope and extension across singular fibers (cf. Theorem \ref{family}).  And the most important fact is 
that the resulting canonical metric $\hat{h}_{can}$ has semipositive curvature 
on the total space of the family. This immediately gives a new proof of  the invariance of plurigenera for smooth projective families (cf. Corollary \ref{pg}). And this result implies the existence of a canonical  hermitian metrics on  the direct image of a  relative pluricanonical system (cf. Theorem \ref{direct}) 
with ``Griffith semipositive'' curvature.  This  semipositivity result is similar to \cite{ka1,v1,v}.    

On the other hand,  it is natural to consider not only a single algebraic variety but also a pair of a variety and a divisor on it.  
One of the important class of such pairs  is the class of KLT pairs (cf. Definition \ref{KLT}).   In general it is a basic philosophy that the most of the results 
for the absolute case (the case of smooth projective varieties) can be 
generalized to the case of KLT pairs.  
Such generalization is   important because the log category is more natural 
to work.  For example the induction in dimension sometimes works more naturally  in the log category (see \cite{b-c-h-m} for example).    

In this paper we define a similar canonical singular hermitian metric on 
the log canonical bundle of  a KLT pair
 with pseudoeffective log canonical divisor.  And it satisfies a similar properties for a projective deformation of KLT pairs (cf. Theorem \ref{logfamily}). 
 By using this metric, we can deduce the invariance of logarithmic plurigenera (cf. Theorems \ref{logplurigenera} and \ref{logplurigeneraL}) and also local freeness and semipositivity of the direct images of pluri log canonical systems (cf. Theorems \ref{logplurigenera} and \ref{direct}).  

Moreover the construction here is applicable to the case of general noncompact complex manifolds such as bounded domains in $\mathbb{C}^{n}$(cf. Section \ref{OPEN}).  This seems to be an interesting topic in future. 

\subsection{Canonical AZD $h_{can}$}\label{hcan}

If we assume that $X$ has nonnegative Kodaira dimension, we have already konwn how to construct a canonical AZD for $K_{X}$.
Let us review the construction in \cite{tu9}.  
\begin{theorem} (\cite{tu9})\label{canazd}
Let $X$ be a smooth projective variety with nonnegative Kodaira dimension. 
We set for every point $x\in X$
\begin{equation}\label{CAZD}
K_{m}(x) :=  \sup \left\{\,\mid \sigma\mid^{\frac{2}{m}}\!\!(x)\, ;\,
\sigma \in \Gamma (X,{\cal O}_{X}(mK_{X})), \left|\int_{X}(\sigma\wedge\bar{\sigma})^{\frac{1}{m}}\right| = 1\right\}
\end{equation}
and 
\begin{equation}
K_{\infty}(x) := \limsup_{m\rightarrow\infty} K_{m}(x).
\end{equation} 
Then 
\begin{equation}\label{CAZD2}
h_{can}:= \mbox{\em the lower envelope of}\,\,K_{\infty}^{-1}
\end{equation}
is an AZD on $K_{X}$. \fbox{}
\end{theorem}
\begin{remark} By the ring structure of $R(X,K_{X})$, we see that
$\{K_{m!}\}$ is monotone increasing, hence  
\[
 \limsup_{m\rightarrow\infty} K_{m}(x) = \sup_{m\geqq 1}K_{m}(x)
\]
holds. \fbox{}
\end{remark}
\begin{remark}
Since $h_{can}$ depends only on $X$,   
the number 
\[
\int_{X}h_{can}^{-1}
\]
is an invariant of $X$.  But I do not know the properties of this number.  \fbox{}
\end{remark}
Apparently this construction is very canonical, i.e., $h_{can}$ depends only on  the complex structure of $X$. 
We call $h_{can}$ the {\bf canonical AZD} of $K_{X}$. 
But this construction works only if we know that the Kodaira dimension of $X$
is nonnegative apriori.  This is the main defect of $h_{can}$. 
For example, $h_{can}$ is useless to solve the abundance conjecture or to 
deduce the deformation invariance of plurigenera.  

Moreover although $h_{can}$ is an AZD of $K_{X}$, it is not clear that $h_{can}$ has minimal singularities  in the sense of 
 Definition \ref{minAZD} below.  But it is easy to see that $h_{can}$ has minimal singularities, if $K_{X}$ is abundant. 

\subsection{Supercanonical AZD $\hat{h}_{can}$}\label{construction}

To avoid the defect of $h_{can}$, we introduce the new AZD $\hat{h}_{can}$.
  Let us use the following terminology. 

\begin{definition}\label{pe}{\bf (Pseudoeffectivity)}  
Let $(L,h_{L})$ be a singular hermitian $\mathbb{Q}$-line bundle on a complex manifold $X$.
$(L,h_{L})$ is said to be  pseudoeffective, if the curvature current
 $\sqrt{-1}\,\Theta_{h_{L}}$ of $h_{L}$
is semipositive.  And a $\mathbb{Q}$-line bundle $L$ on a complex manifold $X$ 
is said to be pseudoeffecive, if there exists a singular hermitian metric 
$h_{L}$ on $L$ with semipositive curvature.  \fbox{}
\end{definition}

Let $X$ be a smooth projective $n$-fold such that the canonical bundle $K_{X}$ 
is pseudoeffective.  Let $A$ be a sufficiently  ample line bundle 
such that for every pseudoeffective singular hermitian line bundle $(L,h_{L})$
on $X$,
${\cal O}_{X}(A + L)\otimes {\cal I}(h_{L})$ and 
 ${\cal O}_{X}(K_{X}+ A + L)\otimes {\cal I}(h_{L})$ are globally generated. 
The existence of such an ample line bundle $A$ follows from Nadel's vanishing theorem (\cite[p.561]{n}). See Propositon \ref{suf} in Section \ref{suff} for detail.   

For every $x\in X$ we set 
\begin{equation}
\hat{K}_{m}^{A}(x) := \sup \left\{ \mid\sigma\mid^{\frac{2}{m}}(x)\mid 
\sigma \in \Gamma (X,{\cal O}_{X}(A + mK_{X})), \parallel\sigma\parallel_{\frac{1}{m}} = 1
\right\}, 
\end{equation}
where 
\begin{equation}\label{1/m}
\parallel\sigma\parallel_{\frac{1}{m}}:= \left|\int_{X}h_{A}^{\frac{1}{m}}\cdot(\sigma\wedge\bar{\sigma})^{\frac{1}{m}}\right|^{\frac{m}{2}}. 
\end{equation}
Here $\mid\sigma\mid^{\frac{2}{m}}$ is not a function on $X$, but the supremum
is takan as a section of the real line bundle $\mid\!A\!\mid^{\frac{2}{m}}\otimes
\mid\!K_{X}\!\mid^{2}$ in the obvious manner\footnote{We have abused the notations $\mid\!\!A\!\!\mid$, $\mid\!\!K_{X}\!\!\mid$ here.  These notations are similar to 
the notations of corresponding linear systems.  But we shall use the notation 
if without fear of confusion. }.  
Then $h_{A}^{\frac{1}{m}}\cdot\hat{K}_{m}^{A}$ is a continuous semipositive $(n,n)$-form on $X$. 
Under the above notations, we have the following theorem.

\begin{theorem}\label{main}
We set 
\begin{equation}
\hat{K}^{A}_{\infty}:= \limsup_{m\rightarrow\infty}h_{A}^{\frac{1}{m}}\cdot\hat{K}_{m}^{A}
\end{equation}
and 
\begin{equation}\label{hcanA}
\hat{h}_{can,A} := \mbox{\em the lower envelope of}\,\,\,(\hat{K}^{A}_{\infty})^{-1}.   
\end{equation}
Then $\hat{h}_{can,A}$ is an AZD of $K_{X}$.
And we define 
\begin{equation}
\hat{h}_{can} : =\mbox{\em the lower envelope of}\,\,\,\inf_{A}\hat{h}_{can,A},
\end{equation}
where $\inf$ denotes the pointwise infimum and $A$ runs all the  
ample line bundles on $X$. 
Then $\hat{h}_{can}$ is  a well defined AZD on $K_{X}$ 
with minimal singularities (cf. Definition \ref{minAZD}) depending only on $X$. \fbox{}
\end{theorem}
\begin{remark}
I believe that $\hat{h}_{can,A}$  is already independent of the sufficiently ample line bundle $A$. \fbox{} 
\end{remark}
\begin{remark} In \cite{kaehler}, I have defined a similar AZD of $K_{X}$ for 
a compact K\"{a}hler manifold $X$ with pseudoeffective canonical bundle.  The construction is even simpler than $\hat{h}_{can}$. But I have not yet proven the 
 semipositivity property corresponding to Theorem \ref{family*} below
  in the case of K\"{a}hler deformations. \fbox{}
\end{remark}
\begin{definition}{\em {\bf (Supercanonical AZD)}}\label{super} 
We call $\hat{h}_{can}$ in Theorem \ref{main} the supercanonical AZD 
of $K_{X}$.   
And we call the semipositive $(n,n)$-form $\hat{h}_{can}^{-1}$ 
the supercanonical volume form on $X$.  \fbox{}
\end{definition}
\begin{remark}
Here ``super'' means that corresponding volume form $\hat{h}_{can}^{-1}$
satisfies the inequality :
\begin{equation}
\hat{h}_{can}^{-1} \geqq h_{can}^{-1}, 
\end{equation}
if $X$ has nonnegative Kodaria dimension (cf. Theorem \ref{comparison}). \fbox{}\end{remark}
In the statement of Theorem \ref{main}, one may think that $\hat{h}_{can,A}$
may depend of the choice of the metric $h_{A}$. 
But later we prove that  $\hat{h}_{can, A}$ is independent of the choice of $h_{A}$
(cf. Lemma \ref{uniqueness}).  

\subsection{Variation of the supercanonical AZD $\hat{h}_{can}$}\label{intro}
Let $f : X \longrightarrow S$ be a fiber space such that  
$X,S$ are complex manifolds and  $f$ is a proper surjective projective morphism with connected fibers.  
Suppose that for every regular fiber $X_{s}:= f^{-1}(s)$, $K_{X_{s}}$ is 
pseudoeffective \footnote{This condition is equivalent to the one that 
for some regular fiber $X_{s}$, $K_{X_{s}}$ is pseudoeffective. This is well known. 
For the proof, see  Lemma \ref{ext} below and Remark \ref{psf}.}.   
In this case we may define a singular hermitian metric $\hat{h}_{can}$  
on $K_{X/S}$ similarly as above.   
Then $\hat{h}_{can}$ have  nice properties on $f : X \longrightarrow S$ 
as follows.  

\begin{theorem}\label{family*}
Let $f : X \longrightarrow S$ be a proper surjective projective morphism with connected fibers between complex manifolds such that 
for every regular fiber $X_{s}$, $K_{X_{s}}$ is pseudoeffective. 
We set $S^{\circ}$ be the maximal nonempty Zariski open subset of $S$ such that $f$ is smooth over $S^{\circ}$ and $X^{\circ} = f^{-1}(S^{\circ})$.
Then there exists a unique singular hermitian metric $\hat{h}_{can}$ on 
$K_{X/S}$ such that 
\begin{enumerate}
\item[(1)] $\hat{h}_{can}$ has  semipositive curvature on $X$, 
\item[(2)] $\hat{h}_{can}\!\mid\!\!X_{s}$ is an AZD of $K_{X_{s}}$ with minimal singularities for 
every $s \in S^{\circ}$,
\item[(3)] For every $s\in S^{\circ}$, 
 $\hat{h}_{can}\!\mid\!X_{s} \leqq \hat{h}_{can,s}$ 
holds, where $\hat{h}_{can,s}$ denotes the supercanonical AZD of $K_{X_{s}}$.
And $\hat{h}_{can}|X_{s} = \hat{h}_{can,s}$ holds outside of a set of measure $0$ on $X_{s}$ for almost every  $s\in S^{\circ}$.  
\end{enumerate} 
We call $\hat{h}_{can}$ in Theorem \ref{family} {\em the relative supercanonical  AZD} on $K_{X/S}$.  \fbox{}
\end{theorem}
To prove Theorem \ref{family*}, first we shall prove the following 
slightly weaker version. 
\begin{theorem}\label{family}
Let $f : X \longrightarrow S$, $S^{\circ}$ and $X^{\circ}:= f^{-1}(S^{\circ})$ 
as in Theorem \ref{family*}. 
 
Then there exists a unique singular hermitian metric $\hat{h}_{can}$ on 
$K_{X/S}$ such that 
\begin{enumerate}
\item[(1)] $\hat{h}_{can}$ has  semipositive curvature on $X$, 
\item[(2)] $\hat{h}_{can}\!\mid\!\!X_{s}$ is an AZD on $K_{X_{s}}$ for 
every $s \in S^{\circ}$,
\item[(3)] There exists the union $F$ of at most countable union of proper subvarieties
of $S^{\circ}$ such that for every $s\in S^{\circ}\,\,\backslash\,\,F$, 
 $\hat{h}_{can}\!\mid\!X_{s} \leqq \hat{h}_{can,s}$ 
holds, where $\hat{h}_{can,s}$ denotes the supercanonical AZD on $K_{X_{s}}$.
And $\hat{h}_{can}|X_{s} = \hat{h}_{can,s}$ holds outside of a set of measure $0$ on $X_{s}$ for almost every  $s\in S^{\circ}$. 
\fbox{} 
\end{enumerate} 
\end{theorem}
The only difference between Theorems \ref{family} and  \ref{family*} is 
the existence of the set $F$ in Theorem \ref{family}.  We prove Theorem \ref{family*} by using Theorem \ref{family} and the invariance of the twisted plurigenera: Corollary \ref{pg2} below (cf. Corollary \ref{empty}). 

In  Theorem \ref{family}, the assertions (1) and  (2) are very important in applications. 
By Theorem \ref{family} (or Theorem \ref{family*}) and the $L^{2}$-extension theorem (\cite[p.200, Theorem]{o-t}), we obtain the following corollary immediately (To make sure we give  
a proof in Section \ref{pgg}).
 
\begin{corollary}(\cite{s1,s2})\label{pg}
Let $f : X \longrightarrow S$ be a smooth projective family over 
a complex manifold $S$. 
Then for every positive integer $m$, the $m$-genus $P_{m}(X_{s}) := \dim H^{0}(X_{s},{\cal O}_{X_{s}}(mK_{X_{s}}))$ is a locally constant function on $S$ \fbox{} 
\end{corollary}

\subsection{Invariance of logarithmic plurigenera}

In Section 4, we shall generalize Theorems \ref{main} and \ref{family*} 
to the case of a projective families of KLT pairs (cf. Definition \ref{KLT}).  See Theorems \ref{main2} and 
\ref{logfamily} below.     
As a consequence we have the invariance of logarithmic plurigenera: 

\begin{theorem}\label{logplurigenera}
Let $f : X \longrightarrow S$ be a proper surjective projective morphism 
between  complex manifolds with connected fibers.  
Let $D$ be an effective $\mathbb{Q}$-divisor on $X$  such that 
\begin{enumerate}
\item[(a)] $D$ is $\mathbb{Q}$-linearly equivalent to a $\mathbb{Q}$-line bundle ($=$ a fractional power of a genuine line bundle) $B$,   
\item[(b)] The set: $S^{\circ} := \{ s\in S|\,\,\mbox{$f$ is smooth over $s$ and $(X_{s},D_{s})$ is
 KLT}\,\,\}$ ($D_{s} := D\!\mid\!X_{s}$) is nonempty. 
\end{enumerate}
Then for every positive integer $m$ such that $mB$ is Cartier, the logarithmic $m$-genus: 
\[
P_{m}(X_{s},B_{s}):= \dim H^{0}\!\left(X_{s},\mathcal{O}_{X_{s}}(m(K_{X_{s}}+B_{s}))\right)
\]
is locally constant on $S^{\circ}$, where $B_{s}:= B|X_{s}$. 
In particular the logarithmic Kodaira dimension of $(X_{s},D_{s})$ is  locally constant over $S^{\circ}$ and for such $m$,  $f_{*}\mathcal{O}_{X}(m(K_{X/S}+B))$ is 
locally free over $S^{\circ}$.  \fbox{} \vspace{3mm}
\end{theorem}
\begin{remark}
It would be interesting to consider a flat family $f : X \to S$ such  that the total space $X$ has only canonical singularities and $K_{X/S} + D$ is $\mathbb{Q}$-Cartier and pseudoeffective. 
I believe that the present proof of Theorem \ref{logplurigenera} works also in this case.  \fbox{}
\end{remark}
We note that in the special case that $B$ is a genuine line bundle, Theorem \ref{logplurigenera} has already been known  (\cite{c,va}).   
In Theorem \ref{logplurigenera}, the canonical choice of $B$ is the minimal positive multiple of $D$ so that the multiple has integral coefficients.  But in general, some smaller positive multiple of $D$ is $\mathbb{Q}$-linearly equivalent to a Cartier divisor. 
The following corollary is obvious. 

\begin{corollary}
Let $f : X \longrightarrow S$,$D$ 
and  $S^{\circ}$ be as in Theorem \ref{logplurigenera}. 
Then we have that for every positive integer $m$, 
\[
\dim H^{0}(X_{s},\mathcal{O}_{X}(\lfloor m(K_{X/S} + D)\rfloor|X_{s}))
\]
is locally constant on $S^{\circ}$ and $f_{*}\mathcal{O}_{X}(\lfloor m(K_{X/S} + D)\rfloor)$ 
is locally free over $S^{\circ}$. \fbox{}
\end{corollary}
\subsection{KLT line bundles and invariance of plurigenera for adjoint line bundles}
In the proof of Theorem \ref{logplurigenera}, for $s\in S^{\circ}$, we consider the singular hermitian metric: 
\[
h_{D,s}:= \frac{1}{|\sigma_{D}|^{2}}\left|X_{s}\right. 
\] 
on $B_{s}$ (see (\ref{sigmaDD}) for the notation), where $\sigma_{D}$ is a multivalued holomorphic section of $B$ (see the convention below)  
with divisor $D$.  The singular hermitian $\mathbb{Q}$-line bundle 
$(B_{s},h_{D,s})$ is an example of the following notion. 

\begin{definition}\label{singKLT}
Let $(L,h_{L})$ be a singular hermitian $\mathbb{Q}$-line bundle 
on a smooth projective variety $X$.  $(L,h_{L})$ said to be KLT (Kawamata log terminal), if the curvature current $\sqrt{-1}\,\Theta_{h_{L}}$ is semipositive and $\mathcal{I}(h_{L}) = \mathcal{O}_{X}$.  For an open subset $U$ of 
$X$, a pseudoeffective line bundle $(L,h_{L})$ on $X$ is said to be KLT over $U$, if $\mathcal{I}(h_{L})|U = \mathcal{O}_{U}$ 
holds.   

A $\mathbb{Q}$-line bundle $L$ on a smooth projective variety is said to be KLT, if 
there exists a singular hermitian metric $h_{L}$ such that 
$(L,h_{L})$ is KLT. \fbox{}
\end{definition}

\noindent Roughly speaking a KLT $\mathbb{Q}$-line bundle is a $\mathbb{Q}$-line bundle which admits 
a singular hermitian metric with semipositive curvature and relatively small singularities.  In this sense,  KLT $\mathbb{Q}$-line bundles  are  somewhere between 
semiample $\mathbb{Q}$-line bundles and pseudoeffective $\mathbb{Q}$-line bundles.  
The notion of  KLT $\mathbb{Q}$-line bundles is a natural generalization of the notion of KLT pairs.   

A very important example of KLT $\mathbb{Q}$-line bundle is the Hodge $\mathbb{Q}$-line  bundle associated with an Iitaka 
fibtration.  
Let $f : X \to Y$ be an  Iitaka fibration such that $X,Y$ are smooth 
and $f$ is a morphism. 
Then by \cite[p.169,Proposition 2.2]{f-m}
$f_{*}\mathcal{O}_{X}(m!K_{X/Y})^{**}$  is invertible on $Y$ for every  sufficiently large $m$, where $**$ denotes the double dual.
$f_{*}\mathcal{O}_{X}(m!K_{X/Y})^{**}$ is of rank $1$ for every sufficiently large
$m$.  
We define the $\mathbb{Q}$-line bundle 
\begin{equation}
L := \frac{1}{m!}\,f_{*}\mathcal{O}_{X}(m!K_{X/Y})^{**} 
\end{equation}
on $Y$ and call it the Hodge $\mathbb{Q}$-line bundle associated with 
$f: X \to Y$.  And for every $y\in Y$ such that $f$ is smooth over $y$, we set 
\begin{equation}
h_{L}^{m!}(\sigma,\sigma)(y) = \left|\int_{X_{y}}(\sigma\wedge\overline{\sigma})^{\frac{1}{m!}}\,\right|^{m!} \hspace{10mm} (\sigma \in L_{y}). 
\end{equation}
and call it the Hodge metric on $L$ at $y$.  
Then $h_{L}$ extends to  a singular hermitian metric on $L$ and $(L,h_{L})$ is KLT 
by the theory of variation of Hodge structures (\cite{sch}). 

Using this new notion, we have a further generalization of Theorem \ref{logplurigenera} as follows. 
 
\begin{theorem}\label{logplurigeneraL}
Let $f : X \longrightarrow S$ be a proper surjective projective morphism 
between  complex manifolds with connected fibers 
and let $(L,h_{L})$ be a pseudoeffective singular hermitian $\mathbb{Q}$-line bundle (cf. Definition \ref{pe} below) on $X$  such that for a general fiber $X_{s}$, $(L,h_{L})|X_{s}$ is KLT, 
We set 
\[
S^{\circ}:= \left\{s\in S|\,\mbox{$f$ is smooth over $s$ and 
$(L,h_{L})|X_{s}$ is well defined and KLT}\,\right\}.
\]
Then for every positive integer $m$ such that $mL$ is Cartier, the twisted $m$-genus: 
\[
P_{m}(X_{s},L_{s}):= \dim H^{0}\!\left(X_{s},\mathcal{O}_{X_{s}}(m(K_{X_{s}}+L_{s}))\right)
\]
is locally constant on $S^{\circ}$, where $L_{s}:= L|X_{s}$.  \fbox{} \vspace{3mm}
\end{theorem}
\noindent In fact Theorem \ref{logplurigenera} follows from Theorem \ref{logplurigeneraL} by taking $(L,h_{L})$ to be $(B,1/|\sigma_{D}|^{2})$.  
The main feature of Theorem \ref{logplurigeneraL} 
is that the singularities  of $h_{L}$ do not appear in the statement 
as long as the singularities are KLT. 
  In this sense, KLT singularities are negligible in this case.  The reason why we can neglect KLT singularities is  that the  construction of the  AZD on $K_{X} + L$ (cf. Theorem \ref{KLTAZD}) does not involve any high powers of $h_{L}$.
 In the special case that $L$ is a genuine line bundle, Theorem \ref{logplurigeneraL} has already been known  (\cite{c,va}).   The main difficulty to deal with  a singular hermitian $\mathbb{Q}$-line bundle is that if we take a multiple
 to make it a genuine line bundle, then we may get a nontrivial multiplier 
 ideal sheaf. 
 
\subsection{A conjecture for K\"{a}hler fibrations}\label{KF}

The invariance of plurigenera is an important consequence of Theorem \ref{family*} or Theorem \ref{family}.  But anyway it has been already known by other methods.  Actually 
one of  the main significance of Theorem \ref{family*} is that it gives a perspective in the case of  K\"{a}hler fibrations as follows. 

Let $f : X \to S$ be a surjective proper K\"{a}hler morphism with connected fibers between connected complex manifolds.   Let $S^{\circ}$ denote the complement of the discriminant locus of $f$.
Let $(L,h_{L})$ be a hermitian line bundle on $X$ with semipositive curvature.
Suppose that $K_{X_{s}}+L|X_{s}$ is pseudoeffective for every $s\in S^{\circ}$. 
We shall consider an analogy of $\hat{h}_{can}$ as follows. For  $s \in S^{\circ}$ we set 
\begin{equation}
\hspace{-10mm} dV_{max}((L,h_{L})|X_{s}):= \mbox{the upper semicontinuous envelope of}
\end{equation}
\[
\sup \left\{ h^{-1}|\,\,\mbox{$h$: a singular hermitian metric on $K_{X}$ 
such that $\sqrt{-1}\,(\Theta_{h} + \Theta_{h_{L}}) \geqq 0$}, \,\,\,
\int_{X}h^{-1} = 1\right\},
\]
where $\sup$ means the poitwise supremum.   We call $dV_{max}((L,h_{L})|X_{s})$ 
the maximal volume form  of $X_{s}$ with respect to $(L,h_{L})|X_{s}$.  
Then it is easy to see that $h_{min,s}:= dV_{max}((L,h_{L})|X_{s})^{-1}\cdot h_{L}$ is an AZD of $K_{X_{s}}+L|X_{s}$ with minimal singularities (see Definition \ref{minAZD} and Section \ref{appendix}).    
This definition can be generalized to the case of noncompact complex manifolds 
such as bounded domains in $\mathbb{C}^{n}$. 
Now I would like to propose the following conjecture. 
\begin{conjecture}\label{CONJ}
In the above notations,  we define the relative volume form 
$dV_{max,X/S}(L,h_{L})$ on $f^{-1}(S^{\circ})$ by $dV_{max,X/S}(L,h_{L})|X_{s} :=  dV_{max}((L,h_{L})|X_{s})$ for $s\in S^{\circ}$ and we define the singular hermitian metric 
$h_{X/S}$ on $(K_{X/S}+ L)|f^{-1}(S^{\circ})$ by 
\[
h_{min,X/S}(L,h_{L}) := \mbox{the lower semicontinuous envelope of}\,\,\, dV_{max,X/S}(L,h_{L})^{-1}\cdot h_{L}.   
\]
Then $h_{min,X/S}(L,h_{L})$ extends to a singular hermitian metric on $K_{X/S} + L$ over $X$ and has semipositive curvature.  \fbox{}
\end{conjecture}
We call $h_{min,X/S}$ the minimal singular hemitian metric on $K_{X/S}+L$ with respect to $h_{L}$. 
This conjecture is very similar to Theorem \ref{family*} and the recent result of Berndtsson(\cite{be}). 
If this conjecture is affirmative, we can prove the defomation ivariance of 
plurigenera for K\"{a}hler deformations.   One can consider also the case that $(L,h_{L})$ is 
pseudoeffective with KLT singularities.  I have also other conjectures 
which relates  the abundance of canonical bundle and the minimal singular hermitian metric on canonical bundle.   I would like to discuss about 
Conjecture \ref{CONJ} and other conjectures in the subsequent papers.  \vspace{3mm}\\

The organization of this article is as follows.   In Section 2, we prove 
Theorem \ref{main}.  In Section 3, we prove Theorem \ref{family*} by using 
a result in \cite[Corollary 4.2]{b-p}.  In Section 4, we generalize the 
reuslts in Sections 2 and 3 to the case of KLT pairs.  Here the new ingredient is the use of dynamical systems of singular hermitian metrics.  
\vspace{3mm}\\  
\noindent{\bf Conventions} 
\begin{itemize}
\item In this paper all the varieties are defined over $\mathbb{C}$.  
\item We frequently use the classical result  that the supremum of a family of plurisubharmonic functions locally uniformly bounded from above is again plurisubharmonic, if we take 
the upper-semicontinuous envelope of the supremum (\cite[p.26, Theorem 5]{l}). 
\item For simplicity, we denote the upper(resp. lower)semicontinuous envelope
simply by the upper(resp. lower) envelope.  
\item In this paper all the singular 
hermitian metrics are supposed to be lower-semicontinuous.
\end{itemize} 
\noindent{\bf Notations}
\begin{itemize}
\item  For a real number $a$, $\lceil a\rceil$ denotes the minimal integer greater than or equal to $a$ and $\lfloor a\rfloor$ denotes the maximal integer smaller than 
or equal to $a$. 
We set $\{ a\} := a -\lfloor a\rfloor$ and call it the fractional part of $a$.   
\item Let $X$ be a projective variety and let $D$ be a Weil divsor on $X$.
Let $D = \sum d_{i}D_{i}$ be the irreducible decomposition.  
We set 
\begin{equation}\label{round}
\lceil D\rceil := \sum \lceil d_{i}\rceil D_{i}
,\,\lfloor D \rfloor := \sum \lfloor d_{i}\rfloor D_{i},\, \{ D\} := \sum \{ d_{i}\}D_{i}.
\end{equation}
\item For a positive integer $n$, $\Delta^{n}$ denotes the unit open polydisk 
in $\mathbb{C}^{n}$ with radius $1$, i.e.,
\[
\Delta^{n} := \{ (t_{1},\cdots ,t_{n})\in\mathbb{C}^{n}; |t_{i}| < 1(1\leqq i\leqq n)\}. 
\] 
We denote $\Delta^{1}$ simply by $\Delta$. 
\item Let $L$ be a $\mathbb{Q}$-line bundle on a compact complex manifold $X$, i.e., $L$ is a formal fractional power of a genuine line bundle on $X$.  
A  singular hermitian metric $h$ on $L$ is given by
\[
h = e^{-\varphi}\cdot h_{0},
\]
where $h_{0}$ is a $C^{\infty}$-hermitian metric on $L$ and 
$\varphi\in L^{1}_{loc}(X)$.
We call $\varphi$ the  weight function of $h$ with respect to $h_{0}$.  We note that even though $L$ is not a genuine line bundle, $h$ makes sense, since a hermitian metric is a real object. 

The curvature current $\Theta_{h}$ of the singular hermitian $\mathbb{Q}$-line
bundle $(L,h)$ is defined by
\[
\Theta_{h} := \Theta_{h_{0}} + \partial\bar{\partial}\varphi,
\]
where $\partial\bar{\partial}\varphi$ is taken in the sense of current.
We define the multiplier ideal sheaf  ${\cal I}(h)$ of $(L,h)$  
by 
\[
\mathcal{I}(h)(U) := \{ f \in \mathcal{O}_{X}(U); \,\, |f|^{2}\,e^{-\varphi}
\in L^{1}_{loc}(U)\},
\]
where $U$ runs open subsets of $X$. 
\item  For  a Cartier divisor $D$, we denote the 
corresponding line bundle by the same notation. 
Let $D$ be an effective $\mathbb{Q}$-divisor on a smooth projective variety $X$.  Let $a$ be a positive integer such that $aD \in \mbox{Div}(X)$. 
We identify $D$ with a formal $a$-th root of the line bundle $aD$.  
We say that $\sigma$ is a multivalued global holomorphic section of $D$ 
with divisor $D$, if  $\sigma_{D}$ is the formal $a$-th root of a nontrivial global holomorphic section of $aD$ with divisor $aD$. And $1/|\sigma_{D}|^{2}$ denotes the singular hermitian metric on $D$ defined by 
\begin{equation}\label{sigmaDD}
\frac{1}{|\sigma_{D}|^{2}} := \frac{h_{0}}{h_{0}(\sigma_{D},\sigma_{D})}, 
\end{equation}
where $h_{0}$ is an arbitrary $C^{\infty}$-hermitian metric on $D$.
\item For a singular hermitian line bundle $(F,h_{F})$ on a compact complex 
manifold $X$ of dimension $n$.   $K(K_{X}+F,h_{F})$ denotes the diagonal part of the Bergman kernel of $H^{0}(X,{\cal O}_{X}(K_{X} + F)\otimes {\cal I}(h_{F}))$ with respect to the $L^{2}$-inner product: 
\begin{equation}\label{inner}
(\sigma ,\sigma^{\prime}) := (\sqrt{-1})^{n^{2}}\int_{X}h_{F}\cdot\sigma\wedge \overline{{\sigma}^{\prime}}, 
\end{equation}  
i.e., 
\begin{equation}\label{BergmanK}
K(K_{X}+F,h_{F}) = \sum_{i=0}^{N}|\sigma_{i}|^{2}, 
\end{equation}
where $\{\sigma_{0},\cdots ,\sigma_{N}\}$ is a complete orthonormal basis 
of $H^{0}(X,{\cal O}_{X}(K_{X} + F)\otimes {\cal I}(h_{F}))$. 
It is clear that $K(K_{X}+F,h_{F})$ is independent of the choice of 
the complete orthonormal basis. \vspace{3mm}
\end{itemize}
\noindent {\bf Acknowledement}:  
I would like to thank to the referee who suggested me the use of Fujita's 
elementary argument  instead of Schmidt's theory of variation of Hodge structure(\cite{sch}) in Section \ref{confam}.        

\section{Proof of Theorem \ref{main}}
In this section we shall prove Theorem \ref{main}. 
We shall use the same notations as in Section \ref{construction}. 
The upper estimate of $\hat{K}_{m}^{A}$ is almost the same as in \cite{tu9}, but the lower estimate of $\hat{K}_{m}^{A}$ requires the $L^{2}$-extension theorem
(\cite{o-t,o}). 

\subsection{Upper estimate of $\hat{K}_{m}^{A}$}\label{up}
Let $X$ be as in Theorem \ref{main} and let $n$ denote $\dim X$.
Let $x\in X$ be an arbitrary point.  
Let $(U,z_{1},\cdots ,z_{n})$ be a coordinate neighborhood of 
$X$ which is biholomorphic to the unit open polydisk $\Delta^{n}$ 
such that $z_{1}(x)= \cdots = z_{n}(x) = 0$.  

Let $\sigma \in \Gamma (X,{\cal O}_{X}(A + mK_{X}))$. 
Taking $U$ sufficiently small, we may assume that $(z_{1},\cdots ,z_{n})$
is a holomorphic local coordinate on a neighborhood of the closure of $U$ and 
there exists 
a local holomorphic frame $\mbox{\bf e}_{A}$ of $A$ on a neighborhood
of the closure of $U$. 
Then there exists a bounded holomorphic function $f_{U}$ on $U$ such that 
\begin{equation}
\sigma = f_{U}\cdot \mbox{\bf e}_{A}\cdot (dz_{1}\wedge\cdots \wedge dz_{n})^{m}\end{equation}
holds.  
Suppose that 
\begin{equation}
\left|\int_{X}h_{A}^{\frac{1}{m}}\cdot(\sigma\wedge\bar{\sigma})^{\frac{1}{m}}\right|= 1
\end{equation}
holds.   
Then we see that 
\begin{eqnarray}\label{submean}
\int_{U}\mid f_{U}(z)\mid^{\frac{2}{m}}d\mu (z) 
& \leqq &
\left(\inf_{U} h_{A}(\mbox{\bf e}_{A},\mbox{\bf e}_{A})\right)^{-\frac{1}{m}}\cdot \int_{U}h_{A}\left(\mbox{\bf e}_{A}\,\mbox{\bf e}_{A}\right)^{\frac{1}{m}}\mid f_{U}\mid^{2}d\mu (z) \\
& \leqq & \left(\inf_{U} h_{A}(\mbox{\bf e}_{A},\mbox{\bf e}_{A})\right)^{-\frac{1}{m}} \nonumber 
\end{eqnarray}
hold, where $d\mu (z)$ denotes the standard Lebesgue measure on the coordinate. Hence by the submeanvalue property of plurisubharmonic functions, 
\begin{equation}
h_{A}^{\frac{1}{m}}\cdot \mid\sigma \mid^{\frac{2}{m}}(x)
\leqq \left(\frac{h_{A}(\mbox{\bf e}_{A},\mbox{\bf e}_{A})(x)}{\inf_{U} h_{A}(\mbox{\bf e}_{A},\mbox{\bf e}_{A})}\right)^{\frac{1}{m}}\!\!\cdot\pi^{-n}\cdot |dz_{1}\wedge\cdots\wedge dz_{n}|^{2}(x)
\end{equation}
holds. 
Let us fix a $C^{\infty}$-volume form $dV$ on $X$. 
Since $X$ is compact and every line bundle on a contractible Stein manifold is trivial, we have the following lemma.
\begin{lemma}\label{upper}
 There exists a positive constant $C$ 
independent of the line bundle $A$ and the $C^{\infty}$-metric $h_{A}$ such that 
\[
\limsup_{m\rightarrow\infty}h_{A}^{\frac{1}{m}}\cdot\hat{K}_{m}^{A} \leqq C\cdot dV
\]
holds on $X$.  
\fbox{} 
\end{lemma}
\subsection{Lower estimate of $\hat{K}_{m}^{A}$}\label{low}
The lower estimate of $\hat{K}_{m}^{A}$ is the essential part of the proof of 
Theorem \ref{main}. 

Let $h_{X}$ be any $C^{\infty}$-hermitian metric on $K_{X}$. 
Let $h_{0}$ be an AZD on $K_{X}$ defined by  
\begin{equation}\label{hzero}
h_{0}:= \mbox{the lower envelope of}\,\,\inf \left\{ h \mid \mbox{$h$ is a singular hermitian metric}\right.
\end{equation}
\[
\left.\mbox{on $K_{X}$ with  $\sqrt{-1}\,\Theta_{h}\geqq 0$,\,$h \geqq h_{X}$}\right\},
\]
where $\inf$ denotes  the pointwise infimum. 
Then by the classical theorem (\cite[p.26, Theorem 5]{l}), $h_{0}$ is an AZD with minimal singularities 
in the sense of Definition \ref{minAZD} below.

Let us compare $h_{0}$ and $\hat{h}_{can}$.
   By the $L^{2}$-extension theorem (\cite{o}), we have the following lemma. 
\begin{lemma}\label{l2ext}
There exists a positive constant $C$ independent of $m$ such that 
\begin{equation}
K(A + mK_{X}, h_{A}\cdot h_{0}^{m-1}) \geqq C\cdot (h_{A}\cdot h_{0}^{m})^{-1}
\end{equation}
holds, where $K(A + mK_{X}, h_{A}\cdot h_{0}^{m-1})$ is the diagonal part of 
Bergman kernel of $A+ mK_{X}$ with respect to the $L^{2}$-inner product:
\begin{equation}
(\sigma ,\sigma^{\prime}):= (\sqrt{-1})^{n^{2}}\int_{X}\sigma\wedge\overline{\sigma^{\prime}}
\cdot h_{A}\cdot h_{0}^{m-1}, 
\end{equation}
where we have considered $\sigma,\sigma^{\prime}$ as  $A+(m-1)K_{X}$ valued
canonical forms (see (\ref{BergmanK})). \fbox{}
\end{lemma}
{\em Proof of Lemma \ref{l2ext}}.
By the extremal property of the Bergman kernel (see for example \cite[p.46, Proposition 1.4.16]{kr}), we have that
\begin{equation}\label{extremal}
K(A + mK_{X}, h_{A}\cdot h_{0}^{m-1})(x) = \sup \left\{ \mid\!\sigma (x)\!\mid^{2}
; \sigma\in\Gamma (X,{\cal O}_{X}(A + mK_{X})\otimes {\cal I}(h_{0}^{m-1})), 
\parallel\!\sigma\!\parallel = 1\right\}
\end{equation}
holds for every $x\in X$, where $\parallel\!\sigma\!\parallel$ denotes the norm $(\sigma,\sigma)^{\frac{1}{2}}$.  
Let $x$ be a point such that $h_{0}$ is not $+\infty$ at $x$.
Let $dV$ be an arbitrary $C^{\infty}$-volume form on $X$ as in Section \ref{construction}.  
Then by the $L^{2}$-extension theorem (\cite{o,o-t}) and 
the sufficient ampleness of $A$ (see Sections \ref{construction} and \ref{suff}),
we  may extend  any  $\tau_{x} \in (A + mK_{X})_{x}$ with $h_{A}\cdot h_{0}^{m-1}\cdot dV^{-1}
(\tau_{x},\tau_{x}) = 1$ to a global section 
$\tau \in  \Gamma (X,{\cal O}_{X}(A + mK_{X})\otimes {\cal I}(h_{0}^{m-1}))$
such that 
\begin{equation}
\parallel\!\tau\!\parallel \leqq C_{0},
\end{equation}
where $C_{0}$ is a positive constant independent of $x$ and $m$.
Let $C_{1}$ be a positive constant such that 
\begin{equation}
h_{0} \geqq C_{1}\cdot dV^{-1}
\end{equation}
holds on $X$.  
By (\ref{extremal}), we obtain the lemma by 
taking $C = C_{0}^{-1}\cdot C_{1}$. \fbox{}\vspace{3mm} \\

Let $\sigma \in \Gamma (X,{\cal O}_{X}(A + mK_{X})\otimes {\cal I}(h_{0}^{m-1}))$ such that 
\begin{equation}
\left|\,\int_{X}\sigma\wedge\bar{\sigma}\cdot h_{A}\cdot h_{0}^{m-1}\right|
= 1
\end{equation}
and 
\begin{equation}
\mid\sigma\mid^{2}(x) = K(A + mK_{X}, h_{A}\cdot h_{0}^{m-1})(x)
\end{equation}
hold, i.e., $\sigma$ is a peak section at $x$. 
Then by the H\"{o}lder inequality we have that  
\begin{eqnarray}\label{holder}
\,\,\,\,\,\,\left|\int_{X}h_{A}^{\frac{1}{m}}\cdot (\sigma\wedge\bar{\sigma})^{\frac{1}{m}}\right|
& \leqq &\left(\int_{X}h_{A}\cdot h_{0}^{m}\cdot \mid\sigma\mid^{2}\cdot h_{0}^{-1}\right)^{\frac{1}{m}}
\cdot \left(\int_{X}h_{0}^{-1}\right)^{\frac{m-1}{m}} \\
& \leqq & \left(\int_{X}h_{0}^{-1}\right)^{\frac{m-1}{m}} \nonumber
\end{eqnarray}
hold. 
Hence we have the inequality:  
\begin{equation}\label{b}
\hat{K}_{m}^{A}(x) \geqq  K(A + mK_{X}, h_{A}\cdot h_{0}^{m-1})(x)^{\frac{1}{m}}
\cdot \left(\int_{X}h_{0}^{-1}\right)^{-\frac{m-1}{m}}
\end{equation}
holds.  Now we shall consider the limit: 
\begin{equation}
\limsup_{m\rightarrow\infty}h_{A}^{\frac{1}{m}}\cdot K(A + mK_{X}, h_{A}\cdot h_{0}^{m-1})^{\frac{1}{m}}.
\end{equation}
Let us recall the following result.

\begin{lemma}(\cite[p.376, Proposition 3.1]{dem})\label{dem}
Let $M$ be a smooth projective variety and let $H$ be a sufficiently ample line bundle on $M$ and let $h_{H}$ be a $C^{\infty}$-hermitian metric on $H$ with strictly positive curvature.  Then for every pseudoeffective singular hermitian line bundle 
$(L,h_{L})$ on  $M$, 
\[
\limsup_{m\rightarrow\infty}h_{H}^{\frac{1}{m}}\cdot K(K_{M}+ H+ mL, h_{H}\cdot h_{L}^{m})^{\frac{1}{m}}= h_{L}^{-1}
\]
holds. \fbox{}
\end{lemma}
\begin{remark}\label{DDDD}
In (\cite[p.376, Proposition 3.1]{dem}, J.P. Demailly only considered the local  version of Lemma \ref{dem}.  But the same proof works in our case 
by the sufficiently ampleness of $H$.  More precisely if we take $H$ to be  sufficiently ample, by the $L^{2}$-extension theorem \cite{o-t,o},  there exists an interpolation operator:
\[
I_{x} : A^{2}(x,(K_{M} + H+ L)_{x},h_{H}\cdot h_{F},\delta_{x})
\rightarrow A^{2}(M,K_{M}+H +L,h_{A}\cdot h_{L})
\]
such that the operator norm of $I_{x}$ is bounded by 
a positive constant independent of $x$ and $(L,h_{L})$, 
where $A^{2}(M,K_{M}+ H + F,h_{H}\cdot h_{L})$ denotes the Hilbert space defined by
\[
 A^{2}(M,K_{M}+ H+ L,h_{H}\cdot h_{L}) := \left\{\,\sigma\in \Gamma (M,{\cal O}_{M}(K_{M}+H + L)), \left|\int_{M}h_{H}\cdot h_{L}\cdot\sigma\wedge\bar{\sigma}\,\right| < + \infty\right\} 
\]
with the $L^{2}$-inner product:  
\[
(\sigma ,\sigma^{\prime}):= (\sqrt{-1})^{n^{2}}\int_{M}h_{H}\cdot h_{L}\cdot\sigma\wedge\overline{\sigma^{\prime}} \hspace{15mm} (n:= \dim M)
\]
and  $A^{2}(x,(K_{M}+H + L)_{x},h_{H}\cdot h_{L},\delta_{x})$ is defined 
similarly, where $\delta_{x}$ is the Dirac measure supported at $x$. 
This is the precise meaning of sufficiently ampleness of $H$ in Lemma \ref{dem}. 
We note that if $h_{L}(x) = + \infty$, then 
$A^{2}(x,(K_{M}+H + L)_{x},h_{H}\cdot h_{L},\delta_{x}) = 0$.
   In this setting,  for a Stein local coordinate neighborhood $U$, 
\[
\hspace{-15mm}\limsup_{m\rightarrow\infty}\,h_{H}^{\frac{1}{m}}\cdot K(K_{M}+H+mL,h_{H}\cdot h_{L}^{m})^{\frac{1}{m}}|\,U = 
\]
\[
\limsup_{m\rightarrow\infty}\,h_{H}^{\frac{1}{m}}\cdot K(K_{M}+H+mL|U, h_{H}\cdot h_{L}^{m}|U)^{\frac{1}{m}} 
 = h_{L}^{-1}|\,U
\] 
hold. This kind of localization principle 
of Bergman kernels is quite standard.  
Moreover for an arbitrary pseudoeffective singular hermitian line bundle $(F,h_{F})$ on $M$, 
\[
\limsup_{m\rightarrow\infty}\,(h_{H}\cdot h_{F})^{\frac{1}{m}}\cdot K(K_{M}+ H+F+ mL, h_{H}\cdot h_{F}\cdot h_{L}^{m})^{\frac{1}{m}}= h_{L}^{-1}
\]
holds almost everywhere on $M$. 

In fact the $L^{2}$-extension theorem (\cite{o-t,o}) implies the 
inequality: 
\begin{equation}
\limsup_{m\rightarrow\infty}\,(h_{H}\cdot h_{F})^{\frac{1}{m}}\cdot K(K_{M}+ H+F+ mL, h_{H}\cdot h_{F}\cdot h_{L}^{m})^{\frac{1}{m}}
\geqq  h_{L}^{-1}
\end{equation}
and the converse inequality is elementary.  See \cite{dem} for details 
and applications.

The reason why we can take $H$ independent of $(L,h_{L})$ is the fact that the  $L^{2}$-extension theorem (\cite{o-t,o}) is uniform with respect to plurisubharmonic weights.  Moreover the extension norm is independent of the weights.     
\fbox{} 
\end{remark} 
We may and do assume that $A$ is sufficiently ample in the sense of 
Lemma \ref{dem}.  Anyway to define $\hat{h}_{can}$ we will replace $A$ by $\ell A$  and take the upper limit as $\ell$ tends to infinity.  
Then by Lemma \ref{dem} letting $m$ tend to infinity in (\ref{b}), we have the following lemma.

\begin{lemma}\label{lower}
\[
\limsup_{m\rightarrow\infty}h_{A}^{\frac{1}{m}}\cdot \hat{K}_{m}^{A} \geqq \left(\int_{X}h_{0}^{-1}\right)^{-1}\cdot h_{0}^{-1}
\]
holds. \fbox{} 
\end{lemma}

By  Lemmas \ref{upper} and \ref{lower}, we see that 
\begin{equation}
\hat{K}_{\infty}^{A} := \limsup_{m\rightarrow\infty}h_{A}^{\frac{1}{m}}\cdot\hat{K}_{m}^{A}
\end{equation}
exists as a bounded semipositive $(n,n)$-form on $X$ ($n= \dim X$).
We set  
\begin{equation}
\hat{h}_{can,A} := \mbox{the lower envelope of}\,\,\,(K_{\infty}^{A})^{-1}. 
\end{equation}

\subsection{Independence of $\hat{h}_{can,A}$ from $h_{A}$}

In the above construction, $\hat{h}_{can,A}$ depends on 
the choice of the $C^{\infty}$-hermitian metric $h_{A}$ apriori.
But we have the following lemma.
\begin{lemma}\label{uniqueness}
$\hat{K}_{\infty}^{A}= \limsup_{m\rightarrow\infty}h_{A}^{\frac{1}{m}}\cdot\hat{K}_{m}^{A}$ is independent of the choice of the $C^{\infty}$-hermitian metric
$h_{A}$. 
Hence $h_{can,A}$ is independent of the choice of the $C^{\infty}$-hermitian metric $h_{A}$.  \fbox{} 
\end{lemma}
 
\noindent{\em Proof of Lemma \ref{uniqueness}}. 
Let $h_{A}^{\prime}$ be another $C^{\infty}$-hermitian metric on $A$. 
We define for $x\in X$ 
\begin{equation}
(\hat{K}_{m}^{A})^{\prime}(x)
:= \sup \left\{ \mid\sigma\mid^{\frac{2}{m}}\!(x) ;\, 
\sigma \in \Gamma (X,{\cal O}_{X}(A + mK_{X})), \left|\int_{X}(h_{A}^{\prime})^{\frac{1}{m}}\cdot (\sigma\wedge\bar{\sigma})^{\frac{1}{m}}\right| = 1\right\}.
\end{equation}
We note that the ratio $h_{A}/h^{\prime}_{A}$ is a positive $C^{\infty}$-function 
on $X$ and 
\begin{equation}
\lim_{m\rightarrow\infty}\left(\frac{h_{A}}{h^{\prime}_{A}}\right)^{\frac{1}{m}} = 1
\end{equation}
uniformly on $X$.   
Since the definitions of $\hat{K}_{m}^{A}$ and $(\hat{K}^{A}_{m})^{\prime}$ 
use the extremal properties, we see easily that for every positive number 
$\varepsilon$, there exists a positive integer $N$ such that 
for every $m \geqq N$
\begin{equation}
(1-\varepsilon )(\hat{K}_{m}^{A})^{\prime}\leqq \hat{K}_{m}^{A} \leqq (1 + \varepsilon )(\hat{K}_{m}^{A})^{\prime}
\end{equation}
holds on $X$.  This completes the proof of Lemma \ref{uniqueness}.
\fbox{}  

\subsection{Completion of the proof of Theorem \ref{main}}

By Lemma \ref{lower} we have the following lemma. 

\begin{lemma}\label{MAZD}
$\hat{h}_{can,A}$ is an AZD with minimal singularities on $K_{X}$. \fbox{}
\end{lemma}
{\em Proof}. Let $h_{0}$ be an AZD  of $K_{X}$ with minimal singularities (cf.
 Definition \ref{minAZD} ) constructed as in (\ref{hzero}).
Then by Lemma \ref{lower} we see that 
\begin{equation}\label{inth}
\hat{h}_{can,A} \leqq \left(\int_{X}h_{0}^{-1}\right)\cdot h_{0}
\end{equation}
holds. 
Hence we see  
\begin{equation}
{\cal I}(\hat{h}_{can,A}^{m}) \supseteq {\cal I}(h_{0}^{m})
\end{equation}
holds for every $m\geqq 1$. This implies that 
\[
H^{0}(X,{\cal O}_{X}(mK_{X})\otimes{\cal I}(h_{0}^{m}))
\subseteq  H^{0}(X,{\cal O}_{X}(mK_{X})\otimes{\cal I}(\hat{h}_{can,A}^{m}))
\subseteq H^{0}(X,{\cal O}_{X}(mK_{X}))
\]
hold, hence 
\[
H^{0}(X,{\cal O}_{X}(mK_{X})\otimes{\cal I}(\hat{h}_{can,A}^{m}))
\simeq  H^{0}(X,{\cal O}_{X}(mK_{X}))
\]
holds for every $m\geqq 1$. 
And by the construction and the classical theorem of Lelong (\cite[p.26, Theorem 5]{l}) 
stated in Section \ref{intro}, $\hat{h}_{can,A}$ has semipositive curvature current.  
Hence $\hat{h}_{can,A}$ is an AZD on $K_{X}$ and depends only on $X$ and $A$
by Lemma \ref{uniqueness}.  We note that $\hat{h}_{can,A}$ is less singular than 
$h_{0}$ by (\ref{inth}).    Since $h_{0}$ has minimal singularities, $\hat{h}_{can,A}$ has minimal singularities, too. 
\fbox{} \vspace{3mm} \\
Let us consider 
\begin{equation}
\hat{K}_{\infty} : = \sup_{A} \hat{K}_{\infty,A}, 
\end{equation}
where $\sup$ means the pointwise supremum and $A$ runs all the 
sufficiently ample line bundle on $X$. 
Then by Lemma \ref{upper}, we see that 
$\hat{K}_{\infty}$ is a well defined semipositive $(n,n)$-form on $X$.
We set 
\begin{equation}
\hat{h}_{can} := \mbox{the lower envelope of}\,\,\, \hat{K}_{\infty}^{-1}. 
\end{equation} 
Then by the construction, $\hat{h}_{can} \leqq \hat{h}_{can,A}$ 
for every ample line bundle $A$. 
Since $\hat{h}_{can,A}$ is an AZD on $K_{X}$, $\hat{h}_{can}$ 
is also an AZD on $K_{X}$ indeed (again by \cite[p.26, Theorem 5]{l}). 
Since $\hat{h}_{can,A}$ depends only on $X$ and $A$, $\hat{h}_{can}$ is uniquely determined by $X$.  By Lemma \ref{lower}, it is clear that $\hat{h}_{can}$ is 
an AZD on $K_{X}$ with minimal singularities in the sense of Definition \ref{minAZD} below.  
This completes the proof of Theorem \ref{main}. \fbox{} 
\begin{remark}\label{minimality} 
As  in Section \ref{low}, we see that $\hat{h}_{can}$ is 
an AZD on $K_{X}$ with minimal singularities (cf. Definition \ref{minAZD}). \fbox{}
\end{remark}
\subsection{Comparison of $h_{can}$ and $\hat{h}_{can}$}
Suppose that $X$ has nonnegative Kodaira dimension. 
Then by Theorem \ref{canazd}, we can define the canonical AZD 
$h_{can}$ on $K_{X}$.   We shall compare $h_{can}$ and $\hat{h}_{can}$. 

\begin{theorem}\label{comparison}
$\hat{h}_{can,A} \leqq h_{can}$ 
holds on $X$.  In particular
 $\hat{h}_{can} \leqq h_{can}$
 holds on $X$  \fbox{} 
\end{theorem}
{\em Proof of Theorem \ref{comparison}}.
If $X$ has  negative Kodaira dimension, then the right hand side 
is infinity.  Hence the inequality is trivial. 
 Suppose that $X$ has nonnegative Kodaira dimension. 
Let $\sigma\in \Gamma (X,{\cal O}_{X}(mK_{X}))$ be an element such that 
\begin{equation}
\left|\int_{X}(\sigma\wedge\bar{\sigma})^{\frac{1}{m}}\right| = 1. 
\end{equation}
Let $x\in X$ be an arbitrary point on $X$. Since 
 ${\cal O}_{X}(A)$ is globally generated by the definition of $A$, there exists an element 
 $\tau \in \Gamma (X,{\cal O}_{X}(A))$ such that 
 $\tau (x) \neq 0$ and  $h_{A}(\tau ,\tau )\leqq 1$ on $X$. 
Then we see that 
\begin{equation}
\left|\int_{X}h_{A}(\tau ,\tau )^{\frac{1}{m}}\cdot (\sigma\wedge\bar{\sigma})^{\frac{1}{m}}\right|\leqq 1
\end{equation}   
holds. 
This implies that 
\begin{equation}
\hat{K}_{m}^{A}(x) \geqq \,|\tau (x)|^{\frac{2}{m}}\cdot K_{m}(x)
\end{equation}
holds at $x$.   Noting $\tau (x) \neq 0$, letting $m$ tend to infinity, 
we see that 
\begin{equation}
\hat{K}^{A}_{\infty}(x) \geqq K_{\infty}(x)
\end{equation}
holds.   Since $x$ is arbitrary, this completes the proof of Theorem \ref{comparison}.  \fbox{}  
\begin{remark}
The equality $h_{can} = \hat{h}_{can}$ implies the abundance of $K_{X}$, if 
the numerical dimension of $(K_{X},\hat{h}_{can})$ is equal to the 
numerical dimension of $K_{X}$.   This problem will be treated in \cite{num}.  
\fbox{}
\end{remark}

By the same proof we obtain the following comparison theorem 
(without assuming $X$ has nonnegative Kodaira dimension).
\begin{theorem}\label{monotonicity}
Let $A, B$ a sufficiently ample line bundle on $X$. 
Suppose that $B - A$ is globally generated, then 
\[
\hat{h}_{can,B} \leqq \hat{h}_{can,A}
\]
holds.  \fbox{} 
\end{theorem}
\begin{remark}\label{any}
Theorem \ref{monotonicity} implies that 
\begin{equation}
\hat{h}_{can} = \lim_{\ell\rightarrow\infty}\hat{h}_{can,\ell A}
\end{equation}
holds for any  ample line bundle $A$ on $X$. \fbox{}
\end{remark}
\begin{remark}\label{independence}
By Kodaira's lemma and Theorem \ref{monotonicity}, we see that $\hat{h}_{can,A}$ is independent of $A$ when $K_{X}$ is big. 
But it is not clear wheter $\hat{h}_{can,A}$ is independent of $A$,  when $K_{X}$ is pseudoeffective but not big.
But by Lemma \ref{lower}, one can easily deduce that for  any two members of  $\{\hat{h}_{can,A}\}$,  the ratio of these metrics is uniformly bounded on $X$, where $A$ runs all the ample line bundles on $X$.  \fbox{}
\end{remark}

\subsection{Canonical volume forms on open manifolds}\label{OPEN}

The construction in Section \ref{hcan} can be generalized to an arbitrary 
complex manifold.   This is just a formal generalization.  But it arises the 
many interesting problems and also is important to consider the degeneration, 
This subsection is not essential in the later argument.  Hence one may skip it. 
Let $M$ be a complex manifold.  
For every positive integer $m$, we set 
\[
Z_{m}:= \left\{\, \sigma\in \Gamma(M,\mathcal{O}_{M}(mK_{M})); \,
\left|\int_{M}(\sigma\wedge\overline{\sigma})^{\frac{1}{m}}\right| < + \infty\,\right\}
\]
and 
\[
K_{M,m} := \sup \left\{\,|\sigma|^{\frac{2}{m}}\,;\,\,\sigma\in \Gamma(M,\mathcal{O}_{M}(mK_{M})), \,
\left|\int_{M}(\sigma\wedge\overline{\sigma})^{\frac{1}{m}}\right| \leqq 1\,\right\},
\]
where $\sup$ denotes the pointwise supremum. 
\begin{proposition}
\[
K_{M,\infty}:= \limsup_{m\to\infty}K_{M,m}
\]
exists and if $Z_{m}\neq 0$ for some $m > 0$, then 
$K_{M,\infty}$ is not identically $0$ and 
\[
h_{can,M}:= \mbox{the lower envelope of}\,\,\,\, \frac{1}{K_{M,\infty}} 
\]
is a well defined singular hermitian metric on $K_{M}$ with semipositive curvature current.  \fbox{}
\end{proposition}
By definition, $h_{can,M}$ is invariant under the automorphism group $\mbox{Aut}(M)$.  
Hence we obtain the following: 
\begin{proposition}
Let $\Omega$ be a homogeneous bounded domain in $\mathbb{C}^{n}$. Then 
$h_{can,\Omega}^{-1}$ is a constant multiple of the Bergman volume form on $\Omega$. 
\fbox{} 
\end{proposition}
For a general bouded domain in $\mathbb{C}^{n}$ it seems to be very difficult to calculate the invariant volume form $h_{can}^{-1}$.   
Let us consider the punctured disk 
\[
\Delta^{*}:= \{ t\in \mathbb{C}|\,\, 0< |t| < 1\}. 
\]
Then one sees that unlike the Bergman kernel, $h_{can,\Delta^{*}}$ reflects the  puncture. 
The following conjecture seems to be very plausible. 
But at this moment I do not know how to solve.    
\begin{conjecture}\label{delta*}
\[
h_{can,\Delta^{*}}^{-1} = O\left(\frac{\sqrt{-1}dt\wedge d\bar{t}}{|t|^{2}\left(\log|t|\right)^{2}}\right)
\]
holds.  
\fbox{} 
\end{conjecture}
Conjecture \ref{delta*} is very important in many senses (see Proposition \ref{FB} below and  Remark \ref{extdeg} for example). In particular it seems to be 
the key to extend Theorem \ref{main2} to LC pairs.      

Next we shall consider the following situation. 
Let $X$ be a smooth projective variety and let $D$ be a divisor with simple normal crossings on $X$.  We set  $M:= X \backslash D$.  Let $A$ be a sufficiently ample line bundle on $X$ in the sense of Proposition \ref{suf}. 
Let $h_{A}$ be a  $C^{\infty}$-hermitian metric on $A$ with strictly positive curvature. 
We define 
\begin{equation}
\hat{K}^{A}_{m}:= \sup \{\,|\sigma|^{\frac{2}{m}};\, \sigma \in \Gamma (M,\mathcal{O}_{M}(A+ mK_{M})), \,
\parallel\sigma\parallel_{\frac{1}{m}} = 1\},
\end{equation}
where  
\begin{equation}
\parallel\sigma\parallel_{\frac{1}{m}}:= \left|\int_{X}h_{A}^{\frac{1}{m}}\cdot(\sigma\wedge\bar{\sigma})^{\frac{1}{m}}\right|^{\frac{m}{2}}. 
\end{equation}
And as (\ref{hcanA}) we define 
\begin{equation}
\hat{h}_{can,A} := \mbox{the lower envelope of}\,\,\,\liminf_{m\to\infty}(\hat{K}^{A}_{m})^{-1}.   
\end{equation}
As Lemma \ref{uniqueness}, we see that $\hat{h}_{can,A}$ is independent of the 
choice of $h_{A}$. 
We set
\[
\hat{h}_{can,M} :=  \mbox{the lower envelope of}\,\,\inf_{A}\hat{h}_{can,A}, 
\]
where $A$ runs all the ample line bundle on $X$. 
We note that 
\[
\{\,\sigma \in \Gamma (M,\mathcal{O}_{M}(A+ mK_{M})),\, 
\parallel\sigma\parallel_{\frac{1}{m}} < \infty\} 
\simeq \Gamma (X,\mathcal{O}_{X}(mK_{X}+(m-1)D))
\]
holds by a simple calculation. 
\begin{definition}\label{FV}
Let $X$ be a smooth projective variety and let $D$ be a divisor with simple normal crossings on $X$. Let $\sigma_{D}$ be a global holomorphic section of 
$\mathcal{O}_{X}(D)$ with divisor $D$.   
$M:= X \backslash D$ is said to be of finite volume,  
if there exists an AZD $h$ of $K_{X} + D$ such that 
\[
\int_{M}h^{-1}\cdot h_{D}
\]
is finite, where $h_{D}:= |\sigma_{D}|^{-2}$.  \fbox{}
\end{definition}
\begin{remark} In the above definition, $h$ is not an AZD of minimal singularities (cf. Definition \ref{minAZD}), when $K_{X} + D$ is ample. \fbox{}  
\end{remark}
\begin{example}
Let $\omega_{E}$ be a complete K\"{a}hler-Einstein form on $M$ such that 
$-\mbox{Ric}_{\omega_{E}} = \omega_{E}$ (\cite{ko}).
We set $n:= \dim X$.  Then $h = (\omega_{E}^{n})\cdot h_{D}$ is an AZD 
on $K_{X} + D$ such that 
\[
\int_{M}h^{-1}\cdot h_{D} = \int_{M}\omega_{E}^{n} < + \infty.
\]
Hence $M$ is of finite volume. \fbox{}
\end{example}
\begin{theorem}\label{open}
Let $X$ be a smooth projective variety and let $D$ be a divisor with simple normal crossings on $X$.    
We set $M := X \backslash D$.  Suppose that $M$ is of finite volume. Then $h_{can,M}\cdot h_{D}$ is an AZD on $K_{X} + D$. 
\fbox{}
\end{theorem}
{\em Proof.}
Let $h_{0}$ be an AZD on $K_{X}+D$ such that $(M,h_{0}^{-1}\cdot h_{D})$ is of finite volume.   
Then by  H\"{o}lder's inequality,  for every $\sigma\in \Gamma (X,\mathcal{O}_{X}(mK_{X}+ (m-1)D))$, 
\begin{equation}\label{holderopen}
\left|\int_{M}h_{A}^{\frac{1}{m}}\cdot h_{D}\cdot (\sigma\wedge\bar{\sigma})^{\frac{1}{m}}\right|
 \leqq \left(\int_{M}h_{A}\cdot h_{0}^{m}\cdot \mid\sigma\mid^{2}\cdot h_{0}^{-1}\cdot h_{D}\right)^{\frac{1}{m}}
\cdot \left(\int_{M}h_{0}^{-1}\cdot h_{D}\right)^{\frac{m-1}{m}} 
\end{equation}
holds.   By the assumption, we know that the second factor on the right-hand side is finite. We shall show the fisrt factor in the right-hand side is finite.  We note that for a local generator $\tau$ of $K_{X} + D$ on an open set 
$U$, 
$\log h_{0}(\tau,\tau)$ is locally bounded from below on $U$, since 
$\log h_{0}(\tau,\tau)$ is plurisuperharmonic. Let $\sigma_{D}$ be as in 
Definition \ref{FV}.   
Then $\parallel\sigma_{D}\parallel^{2}\cdot (h_{0}^{-1}\cdot h_{D})$ is equal 
to $ (|\tau|^{2}\cdot h_{0}(\tau,\tau)^{-1})\cdot(\parallel\sigma_{D}\parallel^{2}\cdot h_{D})$ and is a bounded volume form on $X$, where  $\parallel\sigma_{D}\parallel$ 
denotes the hermitian norm of $\sigma_{D}$ with respect to a fixed 
$C^{\infty}$-hermitian metric on $\mathcal{O}_{X}(D)$.  
 Then since $\sigma$ belongs to $\Gamma (X,\mathcal{O}_{X}(A+mK_{X}+(m-1)D))$, 
\[
\int_{M}h_{A}\cdot h_{0}^{m}\cdot \mid\sigma\mid^{2}\cdot h_{0}^{-1}\cdot h_{D} < \infty
\] 
holds. 
Hence by (\ref{holderopen}), we have the inequality:  
\begin{equation}\label{bopen}
\hat{K}_{m}^{A}(x) \geqq  K(A + mK_{M}, h_{A}\cdot h_{0}^{m-1}\cdot h_{D}^{m-1})(x)^{\frac{1}{m}}
\cdot \left(\int_{M}h_{0}^{-1}\cdot h_{D}\right)^{-\frac{m-1}{m}}
\end{equation}\label{PAZD}
holds.  Since $(M,h_{0}^{-1}\cdot h_{D})$ has finite volume, by the 
same argument as in Section 2.2, letting $m$ tend to infinity, by Lemma \ref{dem} and Remark \ref{DDDD}. 
 , we see that 
\begin{equation}\label{PPP}
\hat{h}_{can,A} \leqq h_{0}\cdot h_{D}^{-1}\cdot\left(\int_{M}h_{0}^{-1}\cdot h_{D}\right)
\end{equation}
holds.   

On the other hand, we obtain the upper estimate of $\hat{K}_{m}^{A}$ 
as follows. 
Let $dV$ be a $C^{\infty}$-volume form on $X$. 
By the submeanvalue inequality for plurisubharmonic functions as in Section \ref{up}, we see that there exists a positive number $C^{\prime}$ independent of $m$ such that 
\[
h_{A}^{\frac{1}{m}}\cdot\hat{K}_{m}^{A}  \leqq C^{\prime}\cdot\frac{dV}{\parallel\sigma_{D}\parallel^{2}}
\]
holds on $X$. 
Hence $\hat{h}_{can,A}$ exists as a well defined singular hermitian metric 
on $K_{X} + D$ and by the construction $\hat{h}_{can,A}$ has semipositive curvature current.     
By (\ref{PPP}) $\hat{h}_{can,A}\cdot h_{D}$ is an AZD on $K_{X}+D$. 
This implies that $\hat{h}_{can,M}\cdot h_{D}$ is an AZD on $K_{X}+D$.  
This completes the proof of Theorem \ref{open}.   
\fbox{} \vspace{3mm} \\
The following problem seems to be interesting. 
\begin{problem}\label{EIN}
Let $X$ be a smooth projective variety and let $D$ be a divisor with only normal crossings on $X$ such that $K_{X} + D$ is ample.   We set $M:= X\backslash D$. Is $\hat{h}_{can,M}^{-1}$ a constant multiple of the K\"{a}hler-Einstein volume form on $X$ constructed in \cite{ko} ?  \fbox{}
\end{problem}
If the above problem is affimative $(M,\hat{h}_{can,M}^{-1})$ is of finite 
volume.  The following is the first step to solve the problem. 
\begin{proposition}\label{FB}
Let $(X,D)$, $\hat{h}_{can}$ be as in Problem \ref{EIN}. 
If Conjecture \ref{delta*} holds, then 
\[
\int_{M}\hat{h}^{-1}_{can,M} < \infty
\]
holds.  \fbox{} 
\end{proposition}
{\em Proof}. 
Let $(U,t_{1},\cdots,t_{n})$ be a local coordinate such that 
\begin{enumerate}
\item $U$ is biholomorphic to $\Delta^{n}$ by $(t_{1},\cdots ,t_{n})$,  
\item $U \cap D = \{(t_{1},\cdots,t_{n})\in \Delta^{n}| t_{1}\cdots t_{k} = 0\}$ holds for some $k$. 
\end{enumerate}
We note that the equality:
\[
\hat{h}_{can,U\backslash D} = h_{can,U\backslash D}
\]
holds as Lemma \ref{uniqueness}, since $A|U$ is trivial with smooth metric.   
For every subset $V$ of $M$, we see that the monotonicity: 
\[
\hat{h}_{can,M}^{-1}\leqq \hat{h}_{can,V}^{-1}
\]
holds by the above construction. 
 Hence we see that  
\[
\hat{h}_{can,M}^{-1}\leqq \hat{h}^{-1}_{can,U\backslash D} 
\] 
holds.  Then by Conjecture \ref{delta*}, we see that 
$\hat{h}_{can,M}$ is of locally of finite volume at every point on 
$U \cap D$.  
This completes the proof.  \fbox{}

\section{Variation of $\hat{h}_{can}$ under projective deformations} 
In this section we shall prove Theorem \ref{family}.
The main ingredient of the proof is the plurisubharmonic variation property of Bergman kernels (\cite{be,b-p,tu8}).

\subsection{Construction of $\hat{h}_{can}$ on a family}\label{confam}

Let $f : X \longrightarrow S$ be a proper surjective projective morphism 
with connected fibers between complex manifolds  
as in Theorem \ref{family}. 

The construction of $\hat{h}_{can}$ can be performed 
simultaeneously on the family as follows.  The same construction works 
for flat projective family with only canonical singularities. 
But for simplicity we shall work on smooth category.

Let $S^{\circ}$  be the maximal nonempty Zariski open subset of $S$ 
such that $f$ is smooth over $S^{\circ}$ and let us set  
$X^{\circ} := f^{-1}(S^{\circ})$. 

Hereafter we shall assume that $\dim S = 1$.  The general case of Theorem \ref{family} easily follows from just by cutting down $S$ to curves (cf. Section \ref{dimgeq1} below). 
 Let $A$ be a sufficiently ample line bundle on $X$  such that for every pseudoeffective  singular hermitian line bundle $(L,h_{L})$, 
${\cal O}_{X}(A + L)\otimes {\cal I}(h_{L})$ and ${\cal O}_{X}(K_{X} + A+ L)
\otimes {\cal I}(h_{L})$ are  globally generated 
and ${\cal O}_{X_{s}}(A + L\!\mid\!\!X_{s})\otimes {\cal I}(h_{L}\!\mid\!\!X_{s})$ and ${\cal O}_{X_{s}}(K_{X_{s}}+ A + L\!\mid\!\!X_{s})\otimes {\cal I}(h_{L}\!\mid\!\!X_{s})$ 
are globally generated for every $s\in S^{\circ}$ as long as $h_{L}\!\!\mid\!X_{s}$ is well defined (cf. Proposition \ref{suf}). 
Let $h_{A}$ be a $C^{\infty}$-hermitian metric on $A$ with strictly positive 
curvature. 
We set 
\begin{equation}\label{em}
E_{m} := f_{*}{\cal O}_{X}(A + mK_{X/S}). 
\end{equation}
Since we have assumed that $\dim S = 1$, $E_{m}$ is a vector bundle on $S$ for every 
$m\geqq 1$.  We denote the fiber of the vector bundle over $s\in S$ by 
$E_{m,s}$. 
Then we shall define the sequence of $\frac{1}{m}A$-valued relative volume forms by  
\begin{equation}\label{KK}
\hat{K}^{A}_{m,s} := \sup \left\{ \mid\!\sigma\!\mid^{\frac{2}{m}} ;
\sigma\in E_{m,s}, 
\left|\int_{X_{s}}h_{A}^{\frac{1}{m}}\cdot(\sigma\wedge\bar{\sigma})^{\frac{1}{m}}\right| = 1\right\}
\end{equation}
for every $s\in S^{\circ}$, where  $\sup$ denotes the 
pointwise supremum.   
This fiberwise construction is different from that in Section \ref{construction}  at the point that we use $E_{m,s}$ instead of $\Gamma (X_{s},{\cal O}_{X_{s}}(A\!\!\mid\!\!X_{s}
+ mK_{X_{s}}))$. 
We note that the difference occurs only over at most countable union of 
proper analytic subsets in $S^{\circ}$ by the upper-semicontinuity theorem 
of cohomologies. \vspace{3mm} \\

\noindent We define the relative $\mid\!A\!\mid^{\frac{2}{m}}$ valued volume form $\hat{K}^{A}_{m}$ by  
\begin{equation}\label{kma}
\hat{K}_{m}^{A}\!\!\mid\!\!X_{s} := \hat{K}^{A}_{m,s} \hspace{5mm}(s\in S^{\circ})
\end{equation}
and a relative volume form $\hat{K}^{A}_{\infty}$ by 
\begin{equation}\label{kinfty}
\hat{K}_{\infty}^{A}\!\!\mid\!\!X_{s} := \limsup_{m\rightarrow\infty}h_{A}^{\frac{1}{m}}\cdot\hat{K}^{A}_{m,s} \hspace{5mm} (s\in S^{\circ}).
\end{equation}
We define a singular hermitian metrics on 
$\frac{1}{m}A + K_{X/S}$ by 
\begin{equation}\label{hma}
\hat{h}_{m,A}:= \mbox{the lower envelope of}\,\,\,(\hat{K}_{m}^{A})^{-1}.
\end{equation}
We set 
\begin{equation}
\hat{h}_{can,A} := \mbox{the lower envelope of}\,\,\,\liminf_{m\rightarrow\infty}
h_{A}^{-\frac{1}{m}}\cdot\hat{h}_{m,A}.
\end{equation}
Then we define  
\begin{equation}
\hat{h}_{can} := \mbox{the lower envelope of}\,\,\, \inf_{A}\hat{h}_{can,A},
\end{equation}
where $A$ runs all the ample line bundles on $X$. 
At this moment, $\hat{h}_{can}$ is defined  only on $K_{X/S}\!\!\mid\!\!X^{\circ}$.
The extension of $\hat{h}_{can}$ to the singular hermitian metric on the whole $K_{X/S}$ will be discussed later.
 
\subsection{Semipositivity of the curvature current of $\hat{h}_{m,A}$}\label{semipos}

Let $E_{m,s}$ denote the fiber of the vector bundle $E_{m}$ at $s$.
For $s\in S^{\circ}$, we define the pseudonorm $\parallel\!\sigma\!\parallel_{\frac{1}{m}}$ 
of $\sigma \in E_{m,s}$ by 
\begin{equation}\label{taut}
\parallel\!\sigma\!\parallel_{\frac{1}{m}}:= \left|\int_{X_{s}}h_{A}^{\frac{1}{m}}\cdot (\sigma\wedge\bar{\sigma})^{\frac{1}{m}}\right|^{\frac{m}{2}}.
\end{equation}


Now we quote the following crucial result.  This is a generalization 
of \cite[p.57,Theorem 1]{ka1}.    
\begin{theorem}(\cite[Corollary 4.2]{b-p})\label{NSB-P}
Let $p : X \to Y$  be a smooth projective fibraition and let 
$(L,h_{L})$ be a pseudoeffective singular hermitian line bundle on $X$.  Let $m$ be a positive integer. 
Suppose that $E:= p_{*}\mathcal{O}_{X}(mK_{X/Y}+L)$ is locally free. 
We set  
\[
K_{m}^{L}(x):= \sup \left\{\,|\sigma|^{\frac{2}{m}}(x)\,;\, 
\sigma\in E_{p(x)}, \left|\int_{X_{p(x)}}h_{L}\cdot(\sigma\wedge \bar{\sigma})^{\frac{1}{m}}\right| = 1\right\}.
\]
Then $h_{NS} := (K^{L}_{m})^{-m}$ is a singular hermitian metric on $mK_{X/Y} + L$ with semipositive curvature current.   \fbox{}  
\end{theorem} 
\begin{remark}\label{BBP}
In  \cite[Corollary 4.2]{b-p}, they have assumed that 
for every $s\in S^{\circ}$, every global holomorphic section of $(mK_{X/Y}+L)|X_{s}$ extends locally to a holomorphic section of $mK_{X/Y}+L$ 
on a neighborhood of $X_{s}$. Apparently they have misunderstood that 
this extension property is equivalent to the local freeness of 
the direct image $E = p_{*}\mathcal{O}_{X}(mK_{X/Y}+L)$. 
Actually without assuming such an extension property, the local freeness of $E$ is automatic in the case of $\dim Y = 1$, since the direct image $E$ is always torsion free. 
 In fact for $y\in Y$ the fiber  $E_{y}$ of the vector bundle $E$ at $y$ is a subspace of $H^{0}(X_{y},\mathcal{O}_{X_{y}}(mK_{X_{y}}+L|X_{y}))$  such that every element of $E_{y}$ is locally holomorphically extendable on a neighborhood of $X_{y}$. 
Hence the proof of  \cite[Corollary 4.2]{b-p} is valid in this case   
and Theorem \ref{NSB-P} is valid as it stands.   

Or we can argue as follows.  By the upper-semicontinuity of 
$h^{0}(X_{y},\mathcal{O}_{X_{y}}(mK_{X_{y}}+L|X_{y}))$, there exists a nonempty 
Zariski open subset $Y_{0}$ of $Y$ such that for every $y\in Y_{0}$, every element of \\ $H^{0}(X_{y},\mathcal{O}_{X_{y}}(mK_{X_{y}}+L|X_{y}))$ extends on a neighborhood of $X_{y}$ as a holomorphic section of $mK_{X/Y}+L$.   
We note the  $h_{L}$ dominates a $C^{\infty}$-metric of $L$ by the assumption 
and $L^{2/m}$-pseudonorm is lower-semicontinuous on $E$ because 
$h_{L}$ is lower-semicontinuous. Hence $K_{m}^{L}(x)$ is locally bounded from above as a section of the real line bundle $|K_{X/Y}|^{2}\otimes |\,L\,|^{\frac{2}{m}}$ over $X$. Let $U$ be an open subset of $X$ such that 
$mK_{X/Y}+L$ has holomorphic frame $\mbox{\bf e}$ on $U$.  Then $\log (K_{m}^{L}/|\mbox{\bf e}|^{2/m})$ is plurisubharmonic on $U\cap f^{-1}(Y_{0})$ by \cite[Corollary 4.2]{b-p}.    
Then since $\log (K_{m}^{L}/|\mbox{\bf e}|^{2/m})$ is locally bounded from above, we may apply the classical extension theorem for plurisubharmonic 
functions (\cite[p.704, Theorem 1.2 (b)]{h-p}). Hence we may extend $h_{NS}$ as a singular hermitian metric with 
semipositive curvature on the whole $mK_{X/Y} + L$.  This argument is better, 
since we do not really use the local freeness of $E$. Hence 
Theorem \ref{NSB-P} holds without assuming the local freeness of  $E$. 
\fbox{}
\end{remark}
Theorem \ref{NSB-P} immediately implies the semipositivity of the curvature 
current of $\hat{h}_{can,A}$. 
\begin{corollary}\label{mApos}
$\hat{h}_{m.A}$ has semipositive curvature current on $X^{\circ}$. \fbox{}
\end{corollary} 
Now let us consider the behavior of $\hat{h}_{m,A}$ along $X \backslash X^{\circ}$.
Let $p$ be a point in $S\backslash S^{\circ}$. 
Since the problem is local, we may and do assume $S$ is the unit open
disk $\Delta$ in $\mathbb{C}$ with center $0$ for the time being 
and $p$ is the origin $0$. 

The following argument is taken from \cite[p. 782,Lemma (1,11)]{f}. Let $\sigma \in \Gamma (X,{\cal O}_{X}(A + mK_{X/S}))$ be a section 
such that $\sigma|X_{0}\neq 0$.   
We consider the (multivalued) $\frac{1}{m}A$-valued relative canonical form:
\begin{equation}
\eta : = \sigma^{\frac{1}{m}}.
\end{equation}
We may and do assume that the support of the fiber $X_{0}$ is a divisor with simple normal crossings.  Let 
\begin{equation}
X_{0} = \sum_{i} \nu_{i}X_{0,i}
\end{equation}
be the irreducible decomposition.  We note that the (multivalued) $\frac{1}{m}A$-valued canonical form
\begin{equation}
f^{*}dt \wedge \eta
\end{equation} 
does not vanish identically on $X_{0}$ by the assumption: $\sigma|X_{0}\neq 0$.   Since the zero divisor of 
$f^{*}dt$ is $\sum_{i}(\nu_{i}-1)X_{i,0}$,  we see that for some $i$, 
$\eta|X_{i,0}$ is a nonzero $\frac{1}{m}A$-valued meromorphic canonical form.    
Hence 
\begin{equation}
\parallel\!\sigma\!\parallel_{\frac{1}{m}}:= \left|\int_{X_{s}}h_{A}^{\frac{1}{m}}(\sigma\wedge\bar{\sigma})^{\frac{1}{m}}\right|^{\frac{m}{2}}
\end{equation}
has a positive lower bound around  $p = 0 \in S$. 
This implies that $h_{m,A}$ has positive lower bound around $X_{0}$. 
Then we see that (fixing a local holomorphic frame of $A + mK_{X/S}$) 
$-\log h_{m,A}$ is locally bounded from above 
around $X_{0}$ and extends across $X_{0}$ as a (local) plurisubhramonic function by \cite[p.704, Theorem 1.2 (b)]{h-p}.    
This implies that $\hat{h}_{m,A}$ is bounded from below by a smooth metric 
along the boundary $X\,\backslash\, X^{\circ}$. 
Hence $\hat{h}_{m,A}$ extends to a singular hermitian metric of $\frac{1}{m}A+K_{X/S}$ with semipositive curvature on the whole $X$ by the same manner as above.   
Now we set 
\begin{equation}
\hat{h}_{can,A} := \mbox{the lower envelope of}\,\,\,\liminf_{m\rightarrow\infty}h_{A}^{-\frac{1}{m}}\cdot\hat{h}_{m,A}.  
\end{equation}
To extend $\hat{h}_{can,A}$ across $X\,\backslash\, X^{\circ}$, we use the 
following useful lemma.

\begin{lemma}(\cite[Corollary 7.3]{b-t})\label{b-t}
Let $\{ u_{j}\}$ be a sequence of plurisubharmonic functions locally bounded above 
on the bounded open set $\Omega$ in $\mathbb{C}^{n}$. Suppose further 
\[
\limsup_{j\rightarrow\infty}u_{j}
\]
is not identically $-\infty$ on any component of $\Omega$. 
Then there exists a plurisubharmonic function $u$ on $\Omega$
such that the set of points: 
\[
\{x \in \Omega \mid u(x) \neq (\limsup_{j\rightarrow\infty}u_{j})(x)\}
\]
is pluripolar.
\fbox{}
\end{lemma}

Since $\hat{h}_{m,A}$ extends to a singular hermitian metric on $\frac{1}{m}A + K_{X/S}$
with semipositive curvature current  on the whole $X$ and 
\begin{equation}\label{HAD}
\hat{h}_{can,A}:= \mbox{the lower envelope of}\,\,\,\liminf_{m\rightarrow\infty}h_{A}^{-\frac{1}{m}}\cdot\hat{h}_{m,A}
\end{equation}
exists as a singular hermitian metric on $K_{X/S}$ on $X^{\circ} = f^{-1}(S^{\circ})$, we see that $\hat{h}_{can,A}$ extends to a singular hermitian metric 
on the whole $X$ with semipositive curvature current  by Lemma \ref{b-t}.

Repeating the same argument we see that $\hat{h}_{can}$ is a well defined 
singular hermitian metric on 
$K_{X/S}\!\mid\!X^{\circ}$ with semipositive curvature current and it extends to a singular hermitian metric 
on $K_{X/S}$ with semipositive curvature current on the whole $X$.

\subsection{Case $\dim S > 1$}\label{dimgeq1}

In Sections \ref{confam},\ref{semipos}, 
we have assumed that $\dim S = 1$.  
In this subsection, we shall extend $\hat{h}_{can}$ as a singular hermitian metric on $K_{X/S}$ over $X$ with semipositive curvature in the case of $\dim S > 1$.
The proof is done just by slicing, i.e., 
we slice  the base $S$ by families of curves and apply classical extension theorems for plurisubharmonic functions or closed semipositve currents.  Let us assume that $\dim S > 1$ holds. 
In this case $E_{m} = f_{*}{\cal O}_{X}(A + mK_{X/S})$ may not be 
locally free on $S^{\circ}$.   If $E_{m}$ is not locally free at $s_{0}\in S^{\circ}$, then $\hat{K}_{\infty}^{A}$ may not be well defined or may be discontinuous  at $s_{0}$, 
because in this case the fiber $E_{m,s_{0}}$ is defined as a maximal linear subspace of $\Gamma (X_{s_{0}},\mathcal{O}_{X_{s_{0}}}(A|X_{s_{0}}+mK_{X_{s_{0}}}))$ such that every element of the subspace is extendable to a holomorphic section of $A + mK_{X/S}$ on a neighborhood of $X_{s_{0}}$. 
See (\ref{EMS}) below.    
We set for $m\geqq 1$. 
\begin{equation}
V_{m}:= \{s \in S^{\circ}\mid E_{m} \,\,\,\mbox{is not locally free at $s$}\}. 
\end{equation}
Then $V_{m}$ is of codimension $\geq 2$, since $E_{m}$ is torsion free. 
By the construction, apriori $\hat{h}_{can}$ is well defined only on 
$S^{\circ} \backslash \cup_{m=1}^{\infty}V_{m}$.  
We note that apriori $\hat{h}_{m,A}$ defined only on  
$f^{-1}(S^{\circ}\backslash V_{m})$.  
Then since $f^{-1}(V_{m})$ is of codimension $\geqq 2$ in $X^{\circ}$, by the 
Hartogs type extension \cite[p.71, Theorem 6]{h}, we may extend $\hat{h}_{m,A}$ across $f^{-1}(V_{m})$ as a singular hermitian metric of $\frac{1}{m}A+ K_{X/S}|X^{\circ}$ with 
semipositive curvature current.  Or more directly, one may use the argument in Remark\ref{BBP} to extend $\hat{h}_{m,A}$ across $f^{-1}(V_{m})$.  The extension theorem \cite[p.71, Theorem 6]{h} is stated for closed semipositive $(1,1)$ currents.   In our case, we need the extension of plurisubharmonic functions.  But these two extensions are obviously related 
by $\partial\bar{\partial}$-Poincar\'{e} lemma (and the Hartogs extension for pluriharmonic functions).  Hence by the costruction, $\hat{h}_{can}$ is  extended to $X^{\circ}$ as a singular hermitian metric on $K_{X/S}|X^{\circ}$. 

Next we shall extend $\hat{h}_{can}$ across $X\backslash X^{\circ}$. We note that the problem is local and birationally invariant (because the pushforwad of a closed semipositive current is again a closed semipositive).  
Hence by taking a suitable modification of $f: X \to S$, we may assume the followings:
\begin{enumerate}
\item[(1)] $S$ is the unit open polydisk:
$\Delta^{k}:= \{(s_{1},\cdots ,s_{k})\in \mathbb{C}^{k};\,|s_{i}|< 1, i=1,\cdots,k\} (k= \dim S > 1)$. 
\item[(2)] $D:= S \backslash S^{\circ}$  is a divisor with normal crossings
on $S$.  
\end{enumerate} 


 Let $C$ be a smooth irreducible curve in $S$ satisfying:  
\begin{enumerate}
\item[(C1)] $C \cap S^{\circ} \neq \emptyset$,
\item[(C2)] $f^{-1}(C)$ is smooth.    
\end{enumerate}
Then by the adjunction formula, we see that 
\begin{equation}\label{bc}
K_{X/S}|f^{-1}(C) = K_{f^{-1}(C)/C}
\end{equation}
holds. For such a curve $C$, noting (\ref{bc}), we may extend $\hat{h}_{can}|f^{-1}(C)\cap X^{\circ}$ 
 to a singular hermitian metric on $K_{X/S}|f^{-1}(C)$ with semipositive 
 curvature by the case of $\dim S = 1$.   

First we shall assume that $f : X \to S$ is flat and  $D$ is smooth.  In this case we may assume that 
\begin{equation}
D = \{ s_{1} = 0\} 
\end{equation}  
holds without loss of generality.
We set
\begin{equation}\label{curve}
C(d_{2},\cdots ,d_{k}):= \{ (s_{1},\cdots,s_{k})\in\Delta^{k}| s_{2}= d_{2},\cdots ,s_{k} = d_{k}\}.  
\end{equation}  
By Bertini's theorem, for a general $(d_{2},\cdots ,d_{k})\in \Delta^{k-1}$ (here ``general'' means outside of a proper analytic subset),   
$C(d_{2},\cdots ,d_{k})$ satisfies the above conditions (C1) and (C2) 
 and $\{f^{-1}(C(d_{2},\cdots ,d_{k}))| (d_{2},\cdots ,d_{k})\in \Delta^{k-1}\}$
 is a flat family over $\Delta^{k-1}$ . 
Let $(s_{1})$ denote the divisor of $s_{1}$ and let 
\begin{equation}
f^{*}(s_{1}) = \sum \nu_{i}X_{i}
\end{equation}
be the irreducible decomposition.  
Let $x \in X_{i,reg} \backslash (\cup_{j\neq i}X_{j})$ be a general (here ``general'' means outside of some proper algebraic subset) point such that there exists a member $C$ in $\{ C(d_{2},\cdots ,d_{k})| (d_{2},\cdots ,d_{k})\in \Delta^{k-1}\}$ such that 
\begin{enumerate}
\item[(1)] $C$ satisfies (C1),(C2),   
\item[(2)] $f^{-1}(C)$ intersects $X_{i,reg}$ transversally at $x$. 
\end{enumerate}
Then (a branch of) $f^{*}s_{1}^{1/\nu_{i}}$ is a local defining function of $X_{i}$ on a neighborhood $W$ of $x$.  And if we take $W$ sufficiently small, we may find  holomorphic functions $z_{1},\cdots ,z_{n} (n = \dim X - \dim S)$ on $W$ such that 
\begin{equation}
(f^{*}s_{1}^{1/\nu_{i}},f^{*}s_{2},\cdots ,f^{*}s_{k},z_{1},\cdots ,z_{n})
\end{equation}
is a local coordinate on $W$.  Since  $\hat{h}_{can}|f^{-1}(C(d_{2},\cdots,d_{k})\cap S^{\circ})\cap W$ extends to a singular hemitian metric on the whole slice $f^{-1}(C(d_{2},\cdots,d_{k}))\cap W$ for every $(d_{2},\cdots,d_{k})\in \Delta^{k-1}$, we see that by \cite[p.710, Theorem2.1 (c)]{h-p},  $\hat{h}_{can}$ extends to a singular hermitian metric with semipositive curvature current on $W$.  In this way we see that $\hat{h}_{can}$ extends to a singular 
hermitian metric across a nonempty Zariski open subset of $X_{i}$ for every $i$.  Then by \cite[p.71, Theorem 6]{h}, we may extend $\hat{h}_{can}$  across 
the whole $\sum_{i} X_{i}$.  Hence in this case we may extend $\hat{h}_{can}$ 
across the boundary $f^{-1}(D)$.  

If $D = S \backslash S^{\circ}$ is reducible and $f : X \to S$ is flat, we extend $\hat{h}_{can}$ across  $f^{-1}(D_{reg})$ as above 
 and then by \cite[p.71, Theorem 6]{h} we extend $\hat{h}_{can}$ across $f^{-1}(D_{sing})$ which is of codimension $\geqq 2$ in $X$ thanks to the flatness of $f$. 
   
If $f : X \to S$ is not flat, we shall take 
a flattening $\hat{f} : \hat{X} \to \hat{S}$ of $f : X \to S$ (cf. \cite{hiro}). In this case $\hat{X}$ and $\hat{S}$ may be singular, but we may and do take them to be normal.
Let $C$ be a curve on $\hat{S}_{reg}$ 
such that $\hat{f}^{-1}(C)\cap \hat{X}_{reg}$ is smooth and $C\cap \hat{S}_{reg}^{\circ}\neq \emptyset$. 
Although $\hat{f}^{-1}(C)$ may be  singular,  taking a resolution of $\hat{f}^{-1}(C)$, by the adjunction formula and the proof in the case of $\dim S = 1$, we may extend $\hat{h}_{can}|f^{-1}(C\cap\hat{S}_{reg}^{\circ})$ to a singular hermitian metric on $K_{\hat{X}_{reg}/\hat{S}_{reg}}|f^{-1}(C)\cap \hat{X}_{reg}$ with semipositive curvature.  
Hence by the above argument, we see that $\hat{h}_{can}$ is a well defined 
singular hermitian metric (with semipositive curvature current) on $K_{\hat{X}_{reg}/\hat{S}_{reg}}$ over $\hat{X}_{reg}\cap \hat{f}^{-1}(\hat{S}_{reg})$. Here we have  abused the same notation $\hat{h}_{can}$ for the metric on the different space. But the metric $\hat{h}_{can}$ is birationally invariant.
 
Let $Z$ be image of $\hat{X}_{sing}\cup \hat{f}^{-1}(\hat{S}_{sing})$ by the natural morphism $\hat{X} \to X$.  Then $Z$ is of codimension at least $2$ in $X$.
Then the above argument in the flat case, $\hat{h}_{can}$ extends to a singular hermitian metric on $K_{X/S}|X \backslash Z$ with semipositive  curvature current.   Then again by \cite[p.71, Theorem 6]{h}, we see that $\hat{h}_{can}$ extends to a 
 singular hermitian metric on $K_{X/S}$ with semipositive curvature current on the whole $X$.     
This completes the proof of the assertion (1) in Theorem \ref{family}.

\subsection{Completion of the proof of Theorem \ref{family}}\label{CL}

To complete the proof of Theorem \ref{family}, we need to show 
that $\hat{h}_{can}$ defines an AZD for $K_{X_{s}}$ for every $s\in S^{\circ}$. 
To show this fact, we modify the construction of $\hat{K}_{m}^{A}$ (cf. (\ref{KK})). 
Here we do not assume  $\dim S = 1$. 

Let us fix $s \in S^{\circ}$ and let $h_{0,s}$ be an AZD with minimal singularities of $K_{X_{s}}$ constructed as (\ref{hzero}), i.e.,
\begin{equation}\label{0s}
h_{0,s}: = \mbox{the upper envelope of}\hspace{50mm}
\end{equation}
\[
\inf \left\{ h| \,\mbox{$h$ is a singular hermitian metric on $K_{X_{s}}$ such that $\sqrt{-1}\,\Theta_{h}\geqq 0$ and $h \geqq h_{s}$}\right\},
\]   
where $h_{s}$ is a fixed $C^{\infty}$-hermitian metric on $K_{X_{s}}$. 
Let $U$ be a neighborhood of $s\in S^{\circ}$ in $S^{\circ}$ which is biholomorphic to the unit open polydisk $\Delta^{k}$ in $\mathbb{C}^{k} (k:= \dim S)$.  On $f^{-1}(U)$, we shall identify $K_{X/S}|U$ with $K_{X}|U$ by tensoring $f^{*}(dt_{1}\wedge\cdots\wedge dt_{k})$, where $(t_{1},\cdots ,t_{k})$ denotes the standard coordinate on $\Delta^{k}$. 
By the $L^{2}$-extension theorem (\cite{o-t,o})
 and the argument modeled after \cite{s1}, we have the following lemma 
 which asserts that $\Gamma (X_{s},{\cal O}_{X_{s}}(A|X_{s} + mK_{X_{s}}))$ 
 contains a ``large'' linear subspace whose elements are extendable on 
 a neighborhood of $X_{s}$.   
\begin{lemma}\label{ext}
Every element of 
$\Gamma (X_{s},{\cal O}_{X_{s}}(A|X_{s} + mK_{X_{s}})
\otimes {\cal I}(h_{0,s}^{m-1}))$ 
extends to an element of  
$\Gamma (f^{-1}(U),{\cal O}_{X}(A + mK_{X}))$
for every positive integer $m$.    
\fbox{} 
\end{lemma}
\begin{remark}\label{psf}
In the proof of Lemma \ref{ext}, we only use the pseudoeffectivity of $K_{X_{s}}$.  
Hence this lemma implies that all the fiber over $U$ has pseudoeffective 
canonical bundles. \fbox{}
\end{remark}

\noindent{\em Proof of Lemma \ref{ext}}. 
We prove the lemma by induction on $m$.
If $m=1$, then the $L^{2}$-extension theorem (\cite{o-t,o}) 
implies that every element of 
$\Gamma (X_{s},{\cal O}_{X_{s}}(A + K_{X_{s}}))$ 
extends to an element of $\Gamma (f^{-1}(U),{\cal O}_{X}(A + K_{X}))$.
 Let $\{\sigma_{1,s}^{(m-1)},\cdots ,\sigma_{N(m-1)}^{(m-1)}\}$
be a basis of  $\Gamma (X_{s},{\cal O}_{X_{s}}(A|X_{s} + (m-1)K_{X_{s}})
\otimes {\cal I}(\tilde{h}_{0,s}^{m-2}))$ for some $m \geqq 2$. 
Suppose that we have already constructed holomorphic extensions:
\begin{equation}
\{\tilde{\sigma}_{1,s}^{(m-1)},\cdots ,\tilde{\sigma}_{N(m-1),s}^{(m-1)}\}
\subset \Gamma (f^{-1}(U),{\cal O}_{X}(A + (m-1)K_{X}))
\end{equation}
of   $\{\sigma_{1,s}^{(m-1)},\cdots ,\sigma_{N(m-1),s}^{(m-1)}\}$
to $f^{-1}(U)$. 
We define the singular hermitian metric $\tilde{h}_{m-1}$ on $(A + (m-1)K_{X})|f^{-1}(U)$ by 
\begin{equation}
\tilde{h}_{m-1} :=
 \frac{1}{\sum_{j=1}^{N(m-1)}|\tilde{\sigma}^{(m-1)}_{j,s}|^{2}}.
\end{equation}
We note that by the choice of $A$, 
${\cal O}_{X_{s}}(A\!\mid\!\!X_{s} + mK_{X_{s}})\otimes {\cal I}(h_{0,s}^{m-1})$ 
is globally generated. 
Hence we see that 
\begin{equation}
{\cal I}(h_{0,s}^{m})\subseteq {\cal I}(h_{0,s}^{m-1})\subseteq {\cal I}(\tilde{h}_{m-1}|X_{s}) 
\end{equation}
hold on $X_{s}$.  Apparently $\tilde{h}_{m-1}$ has a semipositive curvature current. 
Hence by the $L^{2}$-extension theorem (\cite[p.200, Theorem]{o-t}), we  may extend every element of
\begin{equation} 
\Gamma (X_{s},{\cal O}_{X_{s}}(A + mK_{X_{s}})\otimes {\cal I}(h_{0,s}^{m-1}))
\end{equation}
to an element of 
\begin{equation}
\Gamma (f^{-1}(U),{\cal O}_{X}(A + mK_{X})\otimes {\cal I}(\tilde{h}_{m-1})).
\end{equation}
This completes the proof of Lemma \ref{ext} by induction. \fbox{}  \vspace{3mm} \\

\noindent
We set 
\begin{equation}\label{Xi}
\Xi_{m,s}^{A}
:=
\sup \left\{ \mid\sigma\mid^{\frac{2}{m}}
; \,\sigma\in \Gamma (X_{s},{\cal O}_{X_{s}}(A|X_{s} + mK_{X_{s}})\otimes
{\cal I}(h_{0,s}^{m-1})),\,\,\, \left|\int_{X_{s}}h_{A}^{\frac{1}{m}}\cdot 
(\sigma\wedge \bar{\sigma})^{\frac{1}{m}}\right| = 1\right\},
\end{equation}
where  $\sup$ denotes the pointwise supremum. 

Next we shall compare $\Xi_{m,s}^{A}$ with $\hat{K}^{A}_{m.s}$.  
But since $\dim S > 1$, we need to generalize the definition 
of $\hat{K}^{A}_{m,s}$.  Recall that in  the case of $\dim S = 1$, we have defined $\hat{K}^{A}_{m,s}$ as (\ref{KK}).   
In this case  $E_{m}= f_{*}\mathcal{O}_{X}(A + mK_{X/S})$ (cf. (\ref{em})) may not be locally free on 
$S^{\circ}$.   
 For $s\in S^{\circ}$  we define   $E_{m,s}$  by 
 \begin{equation}\label{EMS}
 E_{m,s} := \left\{\,\sigma \in \Gamma (X_{s},\mathcal{O}_{X_{s}}(A|X_{s}+mK_{X_{s}}))| \,\,\mbox{$\sigma$ is extendable to }\right.
\end{equation}  
\[
\mbox{a holomorphic section of $A + mK_{X/S}$ on a neighborhood of $X_{s}$}\}.
\]
This is the right substitute of the fiber 
of $E_{m}$  at $s$ in this case.  
For every $s\in S^{\circ}$, we define  $\hat{K}^{A}_{m,s}$ by  
\begin{equation}
\hat{K}^{A}_{m,s} = \sup \left\{ \mid\sigma\mid^{\frac{2}{m}}
; \,\sigma\in E_{m,s},\,\,\, \left|\int_{X_{s}}h_{A}^{\frac{1}{m}}\cdot 
(\sigma\wedge \bar{\sigma})^{\frac{1}{m}}\right| = 1\right\}.   
\end{equation}
This is the extension of the definition (\ref{KK}) in Section \ref{confam}, where we have assumed that $\dim S = 1$.    
And we set 
\begin{equation}\label{KINF}
\hat{K}^{A}_{\infty,s}:= \limsup_{m\to\infty}h_{A}^{\frac{1}{m}}\cdot\hat{K}^{A}_{m,s}.  
\end{equation}
On the other hand we have already defined $\hat{h}_{can,A}$ over $X$ (cf. (\ref{HAD})).   
And we set 
\begin{equation}
\hat{K}^{A}_{\infty} = \hat{h}_{can,A}^{-1}. 
\end{equation}
By the definition of $E_{m,s}$ (cf. (\ref{EMS})) and the lower-semicontinuity of 
$\hat{h}_{can,A}$, we have that 
\begin{equation}\label{TRR}
\hat{K}_{\infty,s}^{A} \leqq \hat{K}^{A}_{\infty}|X_{s}
\end{equation}
holds for every $s\in S^{\circ}$. 
By Lemma \ref{ext}, we obtain the following lemma immediately. 
\begin{lemma}\label{limit}
\begin{equation}
\limsup_{m\rightarrow\infty}h_{A}^{\frac{1}{m}}\cdot\Xi^{A}_{m,s} \leqq \hat{K}^{A}_{\infty}|X_{s} 
\end{equation}
holds.  \fbox{}
\end{lemma}
{\em Proof}. 
By the definition of $\Xi^{A}_{m,s}$ above and Lemma \ref{ext}
we have that 
\begin{equation}\label{3}
\Xi^{A}_{m,s} \leqq \hat{K}^{A}_{m,s} 
\end{equation}
holds on $X_{s}$. 
On the other hand, by (\ref{KINF}) and (\ref{TRR}), we see that  
\begin{equation}\label{4}
\limsup_{m\rightarrow\infty}h_{A}^{\frac{1}{m}}\cdot\hat{K}^{A}_{m,s} = \hat{K}_{\infty,s}^{A} \leqq \hat{K}^{A}_{\infty}|X_{s}
\end{equation}
hold. 
Hence combining (\ref{3}) and (\ref{4}), we complete the proof of Lemma \ref{limit}. \fbox{} \vspace{3mm}\\

\noindent We set  
\begin{equation}
H_{m,A,s} := (\Xi^{A}_{m,s})^{-1}. 
\end{equation}
Then we have the following lemma. 
\begin{lemma}\label{NN}
If we define 
\begin{equation}
\Xi^{A}_{\infty ,s}:= \limsup_{m\rightarrow\infty}h_{A}^{\frac{1}{m}}\cdot\Xi^{A}_{m,s}
\end{equation}
and 
\begin{equation}
H_{\infty,A,s} := \mbox{{\em the lower envelope of}}\,\,\Xi_{\infty .A,s}^{-1},
\end{equation}
$H_{\infty,A,s}$ is an AZD on $K_{X_{s}}$ with minimal singularities. \fbox{} 
\end{lemma}
{\em Proof}. 
Let $h_{0,s}$ be an AZD on $K_{X_{s}}$ with minimal singularities as (\ref{0s}). 
We note that  \\
${\cal O}_{X_{s}}(A\!\!\mid\!\!X_{s} + mK_{X_{s}})\otimes {\cal I}(h_{0,s}^{m-1})$
is globally generated by the definition of $A$. 
Then by the definition of $\Xi_{m,s}^{A}$,  
\begin{equation}
{\cal I}(h_{0,s}^{m}) \subseteq {\cal I}(H_{m,A,s}^{m})
\end{equation}
holds for every $m\geqq 1$.  Hence by repeating the argument in Section \ref{low},
similar to  Lemma \ref{lower}, we have that 
\begin{equation}
H_{\infty,A,s} \leqq \left(\int_{X_{s}}h_{0,s}^{-1}\right)\cdot h_{0,s}
\end{equation}
holds.    Hence $H_{\infty,A,s}$ is an AZD on $K_{X_{s}}$ with minimal singularities. 
\fbox{}  \vspace{7mm} \\
By the construction of $\hat{h}_{can}$ and Lemma 3.5   
\begin{equation}\label{infer}
\hat{h}_{can}|X_{s} \leqq H_{\infty,A,s}
\end{equation}
holds on $X_{s}$.  
Hence by Lemma \ref{NN} and (\ref{infer}), we see that $\hat{h}_{can}|X_{s}$ is an AZD on $K_{X_{s}}$ with minimal singularities. 
Since $s\in S^{\circ}$ is arbitrary, we see that 
$\hat{h}_{can}\!\!\mid\!\!X_{s}$ is an AZD on $K_{X_{s}}$ 
with minimal singularities for every $s\in S^{\circ}$.  This completes the proof of the assertion (2) in Theorem \ref{family}.
 
We have already seen that the singular hermitian metric $\hat{h}_{can}$ has semipositive curvature  current (cf. Section \ref{semipos}). 
For every $\ell,m \geqq 1$, we set  
\begin{equation}
E^{(\ell)}_{m} = f_{*}{\cal O}_{X_{s}}(\ell A + mK_{X_{s}}).
\end{equation} 
We note that there exists the union $F$ of at most countable 
proper subvarieties of $S^{\circ}$ such that for every $s \in S^{\circ}\,\backslash\, F$ and every $\ell , m\geqq 1$, $E^{(\ell)}_{m}$ is locally free at $s$ and 
\begin{equation}\label{elm}
E^{(\ell)}_{m,s} = \Gamma (X_{s},{\cal O}_{X_{s}}(\ell A|X_{s} + mK_{X_{s}}))
\end{equation}
holds, where $E^{(\ell)}_{m,s}$ denotes the fiber 
of the vector bundle $E^{(\ell)}_{m}$ at $s$. 
Then by the construction and Theorem \ref{monotonicity}(see Remark \ref{any})\footnote{\noindent Theorem \ref{monotonicity} is used because some ample line bundle on the fiber may not extends to an ample line bundle on $X$ in general.} for every $s\in S^{\circ}\,\backslash\,F$, 
\begin{equation}
\hat{h}_{can}\!\!\mid\!X_{s} \leqq  \hat{h}_{can,s}
\end{equation}
holds, where $\hat{h}_{can,s}$ is the supercanonical AZD on $K_{X_{s}}$. 
This completes the proof of the first half of the assertion (3) in Theorem \ref{family}. 
Here the strict inequality may occur on $S^{\circ}$ by the 
effect of the fact that we have taken the lower-semicontinuous envelope in the
construction of $\hat{h}_{can}$.  By the construction it is clear that 
the latter half of the assertion (3) holds.  This completes the proof of Theorem \ref{family}. \fbox{} 

\subsection{Proof of Corollary \ref{pg}}\label{pgg}

Although  I believe that Corollary \ref{pg} is a immediate consequence of  Theorem \ref{family}, to avoid unnecessary misunderstanding, I give a brief proof here. 

Let $f : X \to S$ be a smooth projective family such that $K_{X_{s}}$ is 
pseudoeffective for every $s\in S$. 
We may and do assume that $S$ is the unit open disk $\Delta$ in $\mathbb{C}$. 
We note that there exists a Stein Zariski open subset $U$ of $X$ such that  $K_{X/S}|U$ is trivial.  Then by the $L^{2}$-extension theorem (\cite[p.200, Theorem]{o-t}) and the assertion (1) of Theorem \ref{family}, every element of 
$H^{0}(X_{s},\mathcal{O}_{X_{s}}(mK_{X_{s}})
\otimes \mathcal{I}(\hat{h}_{can}^{m-1}|X_{s}))$ 
extends to an element of 
$H^{0}(X,\mathcal{O}_{X}(K_{X}+ (m-1)K_{X/S})
\otimes \mathcal{I}(\hat{h}_{can}^{m-1}))$ for every $s\in S$. 
By the assertion (2) of Theorem \ref{family}, we see that
\begin{equation}
H^{0}(X_{s},\mathcal{O}_{X_{s}}(mK_{X_{s}})
\otimes \mathcal{I}(\hat{h}_{can}^{m-1}|X_{s}))
\simeq 
H^{0}(X_{s},\mathcal{O}_{X_{s}}(mK_{X_{s}}))
\end{equation}
holds for every $s\in S$.  
Hence every element of $H^{0}(X_{s},\mathcal{O}_{X_{s}}(mK_{X_{s}}))$ extends to an element of $H^{0}(X,\mathcal{O}_{X}(K_{X}+ (m-1)K_{X/S})
\otimes \mathcal{I}(\hat{h}_{can}^{m-1}))$. 
Then since $s$ is arbitrary, by the upper-semicontinuity of cohomologies, we see that 
 the  $m$-genus $h^{0}(X_{s},{\cal O}_{X_{s}}(mK_{X_{s}}))$ is locally constant on $S$.  \fbox{}

\subsection{Tensoring  semipositive $\mathbb{Q}$-line bundles}\label{T1}

In this subsection, we shall consider a minor generalization of 
Theorems \ref{main} and \ref{family} and complete the proof of Theorem \ref{family*}.  

Let $X$ be a smooth projective $n$-fold such that the canonical bundle $K_{X}$ 
is pseudoeffective.  Let $A$ be a sufficiently  ample line bundle 
such that for every pseudoeffective singular hermitian line bundle $(L,h_{L})$
on $X$,
${\cal O}_{X}(A + L)\otimes {\cal I}(h_{L})$ and 
 ${\cal O}_{X}(K_{X}+ A + L)\otimes {\cal I}(h_{L})$ are globally generated. 

Bet $(B,h_{B})$ be a $\mathbb{Q}$-line bundle on $X$ with $C^{\infty}$-hermitian metric with semipositive curvature.  
For every $x\in X$ and a positive integer $m$ such that $mB$ is Cartier, we set 
\begin{equation}
\hat{K}_{m}^{A}(B,h_{B})(x) :=  \sup \left\{\,\mid\sigma\mid^{\frac{2}{m}}\!\!(x);\, \sigma \in \Gamma (X,{\cal O}_{X}(A + m(K_{X}+B))), \parallel\sigma\parallel_{\frac{1}{m}}= 1\right\},
\end{equation}
where
\begin{equation}  
\parallel\sigma\parallel_{\frac{1}{m}} := \left|\int_{X}h_{A}^{\frac{1}{m}}\cdot h_{B}\cdot(\sigma\wedge\bar{\sigma})^{\frac{1}{m}}\right|^{\frac{m}{2}}. 
\end{equation}
Then $h_{A}^{\frac{1}{m}}\cdot\hat{K}^{A}_{m}(B,h_{B})$ is a continuous semipositive $|B|^{2}$-valued $(n,n)$-form on $X$. 
Under the above notations, we have the following theorem.

\begin{theorem}\label{mainB}
We set 
\begin{equation}
\hat{K}_{\infty}^{A}(B,h_{B}):= \limsup_{m\rightarrow\infty}h_{A}^{\frac{1}{m}}\cdot\hat{K}_{m}^{A}(B,h_{B})
\end{equation}
and 
\begin{equation}
\hat{h}_{can,A}(B,h_{B}) := \mbox{\em the lower envelope of}\,\,\,\hat{K}^{A}_{\infty}(B,h_{B})^{-1}.   
\end{equation}
Then $\hat{h}_{can,A}(B,h_{B})$ is an AZD on $K_{X}+B$.
And we define 
\begin{equation}\label{HB}
\hat{h}_{can}(B,h_{B}) : =\mbox{\em the lower envelope of}\,\,\,\inf_{A}\hat{h}_{can,A}(B,h_{B}),
\end{equation}
where $\inf$ denotes the pointwise infimum and $A$ runs all the  
ample line bundles on $X$. 
Then $\hat{h}_{can}(B,h_{B})$ is  a well defined AZD on $K_{X} + B$ 
 with minimal singularities (cf. Definition \ref{minAZD}) depending only on $X$ and $(B,h_{B})$. \fbox{}
\end{theorem}
The proof of Theorem \ref{mainB} is parallel to that of Theorem \ref{main}. 
Hence we omit it.    
We also have the following  generalization of Theorem \ref{family}. 

\begin{theorem}\label{family2}
Let $f : X \longrightarrow S$ be a proper surjective projective morphism with connected fibers between complex manifolds such that 
for a general fiber $X_{s}$, $K_{X_{s}}$ is pseudoeffective. 
We set $S^{\circ}$ be the maximal nonempty Zariski open subset of $S$ such that $f$ is smooth over $S^{\circ}$ and $X^{\circ} = f^{-1}(S^{\circ})$.
Let $(B,h_{B})$ be a $\mathbb{Q}$-line bundle on $X$ with $C^{\infty}$-hermitian metric $h_{B}$ with semipositive curvature on $X$. 
Then there exists a unique singular hermitian metric $\hat{h}_{can}(B,h_{B})$ on 
$K_{X/S}+B$ depending only on $h_{B}$ such that 
\begin{enumerate}
\item[(1)] $\hat{h}_{can}(B,h_{B})$ has  semipositive curvature, 
\item[(2)] $\hat{h}_{can}(B,h_{B})\!\mid\!\!X_{s}$ is an AZD on $K_{X_{s}}+B|X_{s}$ with minimal singularities for 
every $s \in S^{\circ}$,
\item[(3)] There exists the union $F^{\prime}$ of at most countable union of proper subvarieties
of $S^{\circ}$ such that for every $s\in S^{\circ}\,\,\backslash\,\,F^{\prime}$, 
\begin{equation}
\hat{h}_{can}(B,h_{B})\!\mid\!X_{s} \leqq \hat{h}_{can}((B,h_{B})|X_{s})
\end{equation}
holds.
And $\hat{h}_{can}(B,h_{B})|X_{s} = \hat{h}_{can}((B,h_{B})|X_{s})$ holds outside of a set of measure $0$ on $X_{s}$ for almost every  $s\in S^{\circ}$. 
\fbox{} 
\end{enumerate} 
\end{theorem}
{\em Proof.}  The proof is almost parallel to that of Theorem \ref{family}. 
Let $q$ be the minimal positive integer such that $qB$ is a genuine line bundle. The only difference in the proof is that we extend 
\[
H^{0}(X_{s},\mathcal{O}_{X_{s}}(A|X_{s} + (mK_{X} + q\left\lfloor m/q\right\rfloor B_{s}))
\otimes \mathcal{I}(\hat{h}_{can}(B,h_{B})|X_{s})^{m-1}))
\]
by the induction on $m$ similarly as in Lemma \ref{ext}, 
where $A$ is a sufficiently ample line bundle on $X$ independent of $m$. 
 The rest of the proof is completly the same. 
Hence we omit it.  \fbox{} \vspace{3mm} \\

By Theorem \ref{family2} and the $L^{2}$-extension theorem (\cite[p.200, Theorem]{o-t}), we obtain the following corollary immediately.   

\begin{corollary}(\cite{s2})\label{pg2}
Let $f : X \longrightarrow S$ be a smooth projective family over 
a complex manifold $S$ and let $(B,h_{B})$ be a $\mathbb{Q}$-line bundle with 
$C^{\infty}$-hermitian metric $h_{B}$ with semipositive curvature on $X$. 
Then  for every $m\geqq 1$ such that $mB$ is Cartier, the twisted $m$-genus $h^{0}(X_{s},{\cal O}_{X_{s}}(m(K_{X_{s}} + B|X_{s})))$ is a locally constant function on $S$. \fbox{} 
\end{corollary}
{\em Proof}. 
We may and do assume that $S$ is the unit open disk $\Delta$ in $\mathbb{C}$. 
By the $L^{2}$-extension theorem (\cite[p.200, Theorem]{o-t}) 
and  the assertion (1) of Theorem \ref{family2}, for every $s\in S$, every element of 
\[
H^{0}(X_{s},\mathcal{O}_{X}(m(K_{X_{s}}+B|X_{s})
\otimes \mathcal{I}((\hat{h}_{can}(B,h_{B})^{m-1}|X_{s})\cdot h_{B}))
\]
extends to an element of 
\[
H^{0}(X,\mathcal{O}_{X}(K_{X}+ B + (m-1)(K_{X/S}+B))
\otimes \mathcal{I}(\hat{h}_{can}(B,h_{B})^{m-1}\cdot h_{B})). 
\]
By the assertion (2) of Theorem \ref{family2}, we see that for every $s\in S$ 
\begin{equation}
H^{0}(X_{s},\mathcal{O}_{X}(m(K_{X_{s}}+B|X_{s})
\otimes \mathcal{I}((\hat{h}_{can}(B,h_{B})^{m-1}|X_{s})\cdot h_{B}))
\simeq 
H^{0}(X_{s},\mathcal{O}_{X}(m(K_{X_{s}}+B|X_{s}))
\end{equation}
holds. 
Hence every element of $H^{0}(X_{s},\mathcal{O}_{X_{s}}(m(K_{X_{s}}+B|X_{s}))$ extends to an element of $H^{0}(X,\mathcal{O}_{X}(K_{X}+B+ (m-1)(K_{X/S}+B))
\otimes \mathcal{I}(\hat{h}_{can}(B,h_{B})^{m-1}))$.

 Since $s$ is arbitrary,  by the upper-semicontinuity of cohomologies, we see that 
 the twisted $m$-genus $h^{0}(X_{s},{\cal O}_{X_{s}}(m(K_{X_{s}} + B|X_{s})))$ is locally constant on $S$. \fbox{} \vspace{3mm} \\ 
The following corollary slightly improves Theorems \ref{family}  and \ref{family2}.  
\begin{corollary}\label{empty}
The sets $F$ and $F^{\prime}$ in Theorems \ref{family} and \ref{family2} respectively do not exist. \fbox{} 
\end{corollary}
{\em Proof}. 
We shall prove that $F$ is empty.  
We note that $E_{m} = f_{*}\mathcal{O}_{X}(A+mK_{X/S})$ (cf. (\ref{em})) used to define $\hat{h}_{can}$ is locally free over $S^{\circ}$, 
since $h^{0}(X_{s},\mathcal{O}_{X}(A|X_{s}+mK_{X_{s}})) (s\in S^{\circ})$ is locally constant 
over $S^{\circ}$ by Corollary \ref{pg2} and every element of 
    $H^{0}(X_{s},\mathcal{O}_{X}(A|X_{s}+mK_{X_{s}})) (s\in S^{\circ})$ 
    extends to a holomorphic section of $\mathcal{O}_{X}(A+mK_{X/S})$. 
By the same reason, for every $\ell, m \geqq 1$, we see that 
$E_{m}^{(\ell)} = f_{*}\mathcal{O}_{X}(\ell A+mK_{X/S})$ is locally free over $S^{\circ}$ and every element of 
    $H^{0}(X_{s},\mathcal{O}_{X}(\ell A|X_{s}+mK_{X_{s}})) (s\in S^{\circ})$ 
    extends to a holomorphic section of $\mathcal{O}_{X}(\ell A+mK_{X/S})$.  
Then by the construction of $\hat{h}_{can}$ in Section \ref{confam},
viewing the last part of  Section 3.4 (see the definition of $F$ just before (\ref{elm})),  we see that $F$ ought to be empty. 
The emptyness of $F^{\prime}$ follows from the parallel argument. 
 \fbox{} \vspace{3mm} \\ 
\noindent Now we complete the proof of Theorem \ref{family*}. \fbox{} 

\section{Generalization to KLT pairs}\label{KL}

In this section we shall generalize Theorems \ref{main} and \ref{family*}  to the case of KLT pairs.  This leads us to the proof of the invariance of logarithmic plurigenera (Theorem \ref{logplurigenera}).   
Here the essential new techniques are the perturbation of the log canonical bundle  by ample $\mathbb{Q}$-line bundles (cf. Section \ref{GC}) and the use of the dynamical systems of singular hermitian metrics (cf. Section \ref{Dynamical}).
  Here we make use the flexibility of $\mathbb{Q}$-line bundles and the nice  convergence properties of $L^{2/m}$-norms.

\subsection{Statement of the fundamental results}

 First we shall recall  the notion of  KLT pairs.
\begin{definition}\label{KLT}
Let $X$ be a normal variety and let $D= \sum_{i}d_{i}D_{i}$ be an effective  {\bf Q}-divisor such that $K_{X}+D$ is {\bf Q}-Cartier. 
If $\mu : Y \longrightarrow X$ is a log resolution  of the pair 
$(X,D)$, i.e., $\mu$ is a composition of successive blowing ups with smooth centers 
such that $Y$ is smooth and the support of $\mu^{*}D$ is a divisor with simple normal crossings, then we can write
\[
K_{Y} + \mu_{*}^{-1}D = \mu^{*}(K_{X}+D) + F
\]
with $F = \sum_{j}e_{j}E_{j}$ for the exceptional divisors $\{ E_{j}\}$, 
where $\mu_{*}^{-1}D$ denotes the strict transform of $D$. 
We call $F$ the discrepancy and $e_{j}\in  \mbox{\bf Q}$ the discrepancy
coefficient for $E_{j}$.  We define the {\bf log discrepancy}: 
$\mbox{\em ld}(E_{j};X,D)$ at $E_{j}$ by $\mbox{\em ld}(E_{j};X,D):= e_{j} +1$.

The pair $(X,D)$ is said to be {\bf  KLT}(Kawamata log terminal) 
(resp. {\bf LC}(log canonical)), if 
$d_{i}< 1$(resp. $\leqq 1$) for all $i$ and $e_{j} > -1$ (resp. $\geqq -1$)
for all $j$ for a log resolution $\mu : Y \longrightarrow X$.  For a pair $(X,D)$ with $D$ effective, we define the multiplier ideal sheaf   $\mathcal{I}(D)$ of $(X,D)$ by 
$\mathcal{I}(D):= \mu_{*}\mathcal{O}(\lceil F\rceil)$ and 
$CLC(X,D) = \mbox{\em Supp}\,\mathcal{O}_{X}/\mathcal{I}(D)$.  
We call $CLC(X,D)$ {\bf the center of log canonical singularities} of $(X,D)$. 
In this terminnology   $(X,D)$ (with $D$ effective) is KLT, if and only if $CLC(X,D) = \emptyset$.  
  
For an irreducible closed subset $W$ in $X$, we set 
\[
\mbox{\em mld}(\mu_{W};X,D) := \inf_{c_{X}(E) = \mu_{W}}\mbox{\em ld}(E;X,D) 
\] 
and call it the {\bf minimal log discrepancy} at the generic point of $W$ with respect to $(X,D)$, 
where $\mu_{W}$ denotes the generic point of $W$ and the infimum is taken 
over the all effective Cartier divisors $E$ on models of $X$  whose ceneter $c_{X}(E)$ is equal to $W$.    
$\square$ \end{definition}

The following is the counterpart of Theorem \ref{main} in the KLT case. 
\begin{theorem}\label{main2}
Let $(X,D)$ be a KLT pair such that  $X$ is a smooth projective variety. 
Suppose that $K_{X} + D$ is pseudoeffective.  

Then there exists  a singular hermtian metric  $\hat{h}_{can}$ 
on $K_{X} + D$ such that
\begin{enumerate}
\item[(1)] $\hat{h}_{can}$ is uniquely determined by the pair $(X,D)$ 
(see Remark \ref{invariance} below for the precise meaning of the uniqueness),
\item[(2)] $\hat{h}_{can}$ is an AZD  on $K_{X}+D$, i.e., 
\begin{enumerate}
\item[(a)]$\sqrt{-1}\,\Theta_{\hat{h}_{can}}$ is a closed semipositive current, 
\item[(b)] $H^{0}(X,{\cal O}_{X}(m(K_{X} +D))\otimes {\cal I}(\hat{h}_{can}^{m}))
\simeq H^{0}(X,{\cal O}_{X}(m(K_{X} + D)))$
holds for every $m \geqq 1$ such that $mD$ is an integral divisor \footnote{Without this condition ${\cal I}(\hat{h}_{can}^{m})$ is not well defined.}. 
\end{enumerate}
Moreover $\hat{h}_{can}$ is an AZD with minimal singularities (cf. Definition \ref{minAZD}). 
\end{enumerate}
We call $\,\hat{h}_{can}$ in Theorem \ref{main2} the  supercanonical AZD on $K_{X} + D$ on $(X,D)$. 
\fbox{}
\end{theorem}
In Theorem \ref{main2} we have used the same notation $\hat{h}_{can}$ as in 
Theorem \ref{super} for simplicity.   I think this will cause no confusion.  
We call $\hat{h}_{can}$ in Theorem \ref{main2} constructed 
 of $K_{X}+D$ as in Theorem \ref{main}.  The construction is essentially parallel to Theorem \ref{main} and will be given in the next subsection (cf. Theorem 
 \ref{logAZD}) 

As before, we study  the variation of the supercanonical AZD's 
for KLT pairs and prove the following semipositivity 
theorem similar to the  non logarithmic case in Theorem \ref{main}.  

\begin{theorem}\label{logfamily}
Let $f : X \longrightarrow S$ be a proper surjective projective morphism 
between  complex manifolds with connected fibers 
and let $D$ be an effective $\mathbb{Q}$-divisor on $X$  such that 
\begin{enumerate}
\item[(a)] $D$ is $\mathbb{Q}$-linearly equivalent to a $\mathbb{Q}$-line bundle $B$,
\item[(b)] The set: 
 $S^{\circ} := \{ s\in S|\,\,\mbox{$f$ is smooth over $s$ and $(X_{s},D_{s})$ is  KLT}\,\,\}$ 
is nonempty, 
\item[(c)] For every $X_{s} (s\in S^{\circ})$  , $K_{X_{s}} + D_{s}$ is pseudoeffective
\footnote{Here actually we only need to assume that for some  fiber $X_{s}(s\in S^{\circ})$. 
$(X_{s},D_{s})$ is KLT and $K_{X_{s}}+D_{s}$ is pseudoeffective. 
See Theorem \ref{PE} below.}. 
\end{enumerate}
Then there exists a singular hermitian metric $\hat{h}_{can}$ on 
$K_{X/S}+ B$ such that  
\begin{enumerate}
\item[(1)] $\hat{h}_{can}$ has semipositive curvature  current, 
\item[(2)] $\hat{h}_{can}\!\mid\!\!X_{s}$ is an AZD on $K_{X_{s}}+B_{s}$ for 
every $s \in S^{\circ}$,
\item[(3)] For every $s\in S^{\circ}$, $\hat{h}_{can}\!\!\mid\!X_{s} \leqq \hat{h}_{can,s}$ 
holds, where $\hat{h}_{can,s}$ denotes the supercanonical AZD on 
$K_{X_{s}}+B_{s}$. And $\hat{h}_{can}|\!X_{s} = \hat{h}_{can,s}$ holds outside of a set of measure $0$  on $X_{s}$ for almost every $s\in S^{\circ}$. 
\end{enumerate} 
 $\square$
\end{theorem}

\noindent In the proof of Theorem \ref{logfamily}, we do need to use the fact that $D$ is effective, since we need to use the semipositivity result similar 
to Lemma \ref{NSB-P}. 

\subsection{Construction of the supercanonical AZD's for  KLT pairs}
In this subsection, we shall construct the supercanonical AZD's for 
 KLT pairs similarly to Theorem \ref{main}. 
Let $(X,D)$ be a sub KLT pair such that $X$ is smooth and 
$K_{X} + D$ is pseudoeffective.  In this subsection we shall consider 
$D$ as a $\mathbb{Q}$-line bundle.  Because we are considering singular hermitian metrics on $K_{X}+D$ or its multiples, this is not a problem.  
 
Let $A$ be a sufficiently ample line bundle such that 
for any pseudoeffective singular hermitian line bundle $(L,h_{L})$ (cf. Definition \ref{pe}), 
${\cal O}_{X}(A+L)\otimes {\cal I}(h_{L})$ and 
${\cal O}_{X}(A+K_{X} + L)\otimes {\cal I}(h_{L})$ are globally generated. 
Such an ample line bundle $A$ exists by Proposition \ref{suf} below.   
Let $D = \sum_{i}d_{i}D_{i}$ be the irreducible decomposition of $D$ and 
for every $i$ we choose a nonzero global holomorphic section $\sigma_{D_{i}}$ of $\mathcal{O}_{X}(D_{i})$ with divisor $D_{i}$. 
For every positive integer $m$ such that $mD \in \mbox{Div}(X)$, we set 
\begin{equation}
\hat{K}_{m}^{A}:= \sup \left\{\,|\sigma|^{\frac{2}{m}} ;
\sigma \in \Gamma (X,{\cal O}_{X}(A +m(K_{X}+D))), 
\parallel\sigma\parallel_{\frac{1}{m}} = 1\right\}, 
\end{equation}
where $\sup$ denotes the pointwise supremum  and 
\[
\parallel\sigma\parallel_{\frac{1}{m}}:=\left|\int_{X} h_{A}^{\frac{1}{m}}\cdot h_{D}\cdot(\sigma\wedge\bar{\sigma})^{\frac{1}{m}}\right|,
\]
where 
\begin{equation}\label{HD}
h_{D}:= \frac{1}{\prod_{i}|\sigma_{D_{i}}|^{2d_{i}}}. 
\end{equation}
Here $\mid\sigma\mid^{\frac{2}{m}}$ is not a function on $X$, but the supremum
is takan as a section of the real line bundle $\mid\!A\!\mid^{\frac{2}{m}}\otimes \mid\!K_{X}+D\!\mid^{2}$ in the obvious manner.  
 Then 
\begin{equation}
\hat{h}_{m,A} : = (\hat{K}_{m}^{A})^{-1}
\end{equation}
is a singular hermitian metric on $m^{-1}A+(K_{X} + D)$
with semipositive curvature current. 
Then  $h_{A}^{-1/m}\cdot\hat{h}_{m,A}$ 
is a singular hermitian metric on $K_{X} + D$. 
Then 
\begin{equation}
\hat{h}_{can,A}: = \mbox{the lower envelope of}\,\, 
\liminf_{m\rightarrow\infty} h_{A}^{-\frac{1}{m}}\cdot\hat{h}_{m,A}
\end{equation}
is a singular hermitian metric on $K_{X}+D$ with semipositive curvature current,where $m$ runs all the  integers such that $mD\in \mbox{Div}(X)$.  Now we have the following theorem similar to Theorem \ref{main}. 
\begin{lemma}\label{logAZD}
$\hat{h}_{can,A}$  
and 
\begin{equation}
\hat{h}_{can}:= \mbox{the lower envelope of}\,\,\,\inf_{A}h_{can,A}
\end{equation}
are well defined and are AZD's of $K_{X} + D$ with minimal singularities,  where $A$ runs all the sufficiently ample line bundles on $X$.
$\square$
\end{lemma}
{\em Proof.} The proof of Lemma \ref{logAZD} is  almost parallel to 
the proof of Theorem \ref{super}.  Let $h_{0}$ be an AZD with minimal singularities on $K_{X} +D$ constructed in (\cite[Theorem 1.5]{d-p-s}) as in Section \ref{appendix}.  Then by H\"{o}lder's inequality,  
\[
\left|\int_{X}h_{A}^{\frac{1}{m}}\cdot h_{D}\cdot(\sigma\wedge\bar{\sigma})^{\frac{1}{m}}\right|
\leqq \left(\int_{X}h_{A}\cdot h_{0}^{m}\cdot \mid\sigma\mid^{2}\cdot (h_{0}^{-1}\cdot h_{D})\right)^{\frac{1}{m}}
\cdot \left(\int_{X}h_{0}^{-1}\cdot h_{D}\right)^{\frac{m-1}{m}}, 
\]
holds, where $m$ is a positive integer $m$ such that $mD \in \mbox{Div}(X)$, and $\sigma \in \Gamma (X,\mathcal{O}_{X}(A + m(K_{X}+D))$.
This inequality makes sense, because 
\[
\int_{X}h_{0}^{-1}\cdot h_{D} < + \infty
\]  
holds,  since $(X,D)$ is KLT. 
Hence we have the inequality:   
\[
\hat{K}_{m}^{A} \geqq K(A+m(K_{X}+D),h_{A}\cdot h_{0}^{m-1}\cdot h_{D})^{\frac{1}{m}}\cdot\left(\int_{X}h_{0}^{-1}\cdot h_{D}\right)^{-\frac{m-1}{m}}. 
\]
And  by Lemma \ref{dem} and Remark \ref{DDDD}, if $A$ is sufficiently ample, letting $m$ tend to infinity, we have that  
\begin{equation}\label{minca}
\hat{h}_{can,A} \leqq h_{0}\cdot \left(\int_{X}h_{0}^{-1}\cdot h_{D}\right). 
\end{equation}
holds.  On the other hand the upper estimate of $\hat{K}^{A}_{m}$ is obtained 
by the submeanvalue inequality for plurisubharmonic functions. 

Let $\sigma \in \Gamma (X,{\cal O}_{X}(A + m(K_{X}+D)))$ for some $m$ such that  $mD$ is integral. 
Let $(U,(z_{1},\cdots ,z_{n}))$ be a local coordinate on $X$ which is 
biholomorphic to the unit open polydisk $\Delta^{n}$ by the coordinate 
$(z_{1},\cdots ,z_{n})$. 
Taking $U$ sufficiently small, we may assume that $(z_{1},\cdots ,z_{n})$
is a holomorphic local coordinate on a neighborhood of the closure of $U$ and 
there exist local holomorphic frames $\mbox{\bf e}_{A}$ of $A$ 
and $\mbox{\bf e}_{D_{i}}$ of $D_{i}$ for every $i$  respectively on a neighborhood of the closure of $U$. 
Then there exists a bounded holomorphic function $f_{U}$ on $U$ such that 
\begin{equation}
\sigma = f_{U}\cdot \mbox{\bf e}_{A}\cdot (dz_{1}\wedge\cdots \wedge dz_{n})^{m}\cdot \left(\prod_{i}\mbox{\bf e}_{D_{i}}^{d_{i}}\right)^{m}
\end{equation}
holds.  
Suppose that 
\begin{equation}
\left|\int_{X}h_{A}^{\frac{1}{m}}\cdot h_{D}\cdot(\sigma\wedge\bar{\sigma})^{\frac{1}{m}}\right|= 1
\end{equation}
holds. 
Then since $h_{D}\cdot \prod_{i}|\mbox{\bf e}_{D_{i}}|^{2d_{i}}$ and 
$h_{A}(\mbox{\bf e}_{A},\mbox{\bf e}_{A})$ has positive lower bound 
on $U$, as (\ref{submean}) by the submeanvalue inequality for plurisubharmonic functions, 
we see that $|f_{U}|^{\frac{2}{m}}$ is bounded compact uniformly on $U$.  
Hence just as Lemma \ref{upper},  we have that  there exists a positive 
constant $C$ such that 
\begin{equation}\label{upper2}
\limsup_{m\to\infty}h_{A}^{\frac{1}{m}}\cdot h_{D}\cdot \hat{K}^{A}_{m} 
\leqq C \cdot \frac{dV}{{\prod_{i}\parallel\sigma_{D_{i}}\parallel^{2d_{i}}}}
\end{equation}
holds, where for every $i$, $\parallel\sigma_{D_{i}}\parallel$ denotes the hermitian norm 
of $\sigma_{D_{i}}$ with respect to a fixed $C^{\infty}$-hermitian metric 
on $D_{i}$ and $dV$ is a fixed $C^{\infty}$-volume form on $X$.  

Combining (\ref{minca}) and (\ref{upper2}), $\hat{h}_{can,A}$ is not identically $0$ and is  an AZD on $K_{X}+D$ with minimal singularities. 
This completes the proof of Lemma \ref{logAZD}. \vspace{3mm}\fbox{}\\
 \noindent By Lemma \ref{logAZD} and the definiton of $\hat{h}_{can}$, $\hat{h}_{can}$ is an AZD on $K_{X}+D$ with minimal singularities indeed.  And the rest of the proof is similar to that of Theorem \ref{main}.  This completes the proof of Theorem \ref{main2}. \fbox{} 

\begin{remark}\label{invariance}
In the above proof of Theorem \ref{main2}, $\hat{K}_{m}^{A}$ depends on 
the choice of $\{\sigma_{D_{i}}\}$.   
Nevertheless the singular volume form: 
\[
\frac{h_{A}^{\frac{1}{m}}\cdot\hat{K}_{m}^{A}}{\prod_{i}|\sigma_{D_{i}}|^{2d_{i}}}
\]
does not depend on the choice. Hence we see that the singular volume form: 
\[
\frac{\hat{h}_{can}^{-1}}{\prod_{i}|\sigma_{D_{i}}|^{2d_{i}}}
\]
is uniquely determined by $(X,D)$. In other words, $\hat{h}_{can}^{-1}$ is 
uniquely determined as a singular volume form on $X$ and it does not 
depend on the choice of $\{\sigma_{D_{i}}\}$. \fbox{}
\end{remark}
\subsection{Construction of supercanonical AZD's on  adjoint line bundles}\label{KLTL*}
Let $(L,h_{L})$ be a KLT singular hermitian $\mathbb{Q}$-line bundle (cf. Definition \ref{singKLT})  on a smooth projective variety $X$.  Suppose that $K_{X}+L$ is pseudoeffective.  
Let $A$ be a sufficiently ample line bundle on $X$ in the sense of Proposition \ref{suf} in Appendix and let $h_{A}$ be a $C^{\infty}$-hermitian metric on $A$. 
For a positive integer $m$ such that $mL$ is a genuine line bundle and $\sigma \in \Gamma (X,\mathcal{O}_{X}(A + m(K_{X}+L)))$, we set 
\begin{equation}\label{L1/m}
\parallel\sigma\parallel_{\frac{1}{m}}:= \left|\int_{X}h_{A}^{\frac{1}{m}}\cdot h_{L}\cdot(\sigma\wedge\bar{\sigma})^{\frac{1}{m}}\right|^{\frac{m}{2}}. 
\end{equation}
For $x\in X$, we set 
\begin{equation}
\hat{K}^{A}_{m}(x) := \sup \left\{ \mid\sigma\mid^{\frac{2}{m}}(x)\mid 
\sigma \in \Gamma (X,{\cal O}_{X}(A + m(K_{X}+L))), \parallel\sigma\parallel_{\frac{1}{m}} = 1
\right\}. 
\end{equation} 
We note that $\parallel\sigma\parallel_{\frac{1}{m}}$ is well defined by the assumption that $(L,h_{L})$ is KLT. 
We set 
\begin{equation}\label{LAAZD}
\hat{h}_{can,A}(L,h_{L}) := \mbox{the lower envelope of}\,\,
\liminf_{m\to\infty}h_{A}^{-\frac{1}{m}}\cdot (\hat{K}_{m}^{A})^{-1} 
\end{equation}
and 
\begin{equation}\label{LAZD}
\hat{h}_{can}(L,h_{L}):= \mbox{the lower envelope of}\,\,\inf_{A}\hat{h}_{can,A}(L,h_{L}), 
\end{equation}
where $A$ runs all the ample line bundles on $X$.  

\begin{theorem}\label{KLTAZD}
$\hat{h}_{can,A}(L,h_{L})$ and $\hat{h}_{can}(L,h_{L})$ defined respectively as (\ref{LAAZD}) and 
(\ref{LAZD}) are AZD's of $K_{X}+L$ with minimal singularities. 
We call $\hat{h}_{can}(L,h_{L})$ the supercanonical AZD on $K_{X} + L$ 
with respect to $h_{L}$. \fbox{}
\end{theorem}
The proof of Theorem \ref{KLTAZD} is completely parallel to 
the one of Theorem \ref{main2} above.  In fact, for example the lower estimate 
of $\hat{K}_{m}^{A}$ as follows. 
Let $h_{0}$ is a AZD with minimal singularities 
on $K_{X} + L$.  Then similarly as \ref{holder})  the inequality:
\[
\left|\int_{X}h_{A}^{\frac{1}{m}}\cdot h_{L}\cdot(\sigma\wedge\bar{\sigma})^{\frac{1}{m}}\right|
\leqq \left(\int_{X}h_{A}\cdot h_{0}^{m}\cdot \mid\sigma\mid^{2}\cdot (h_{0}^{-1}\cdot h_{L})\right)^{\frac{1}{m}}
\cdot \left(\int_{X}h_{0}^{-1}\cdot h_{L}\right)^{\frac{m-1}{m}}, 
\]
holds. where $m$ is a positive integer $m$ such that $mL$ is a genuine line bundle and  $\sigma \in \Gamma (X,\mathcal{O}_{X}(A + m(K_{X}+L)))$. 
We note that this inequality makes sense, since $(L,h_{L})$ is KLT.  And the rest is similar to that of Theorem \ref{main2}.  
Hence we omit it.
   
\subsection{Proof of Theorem \ref{logfamily}; Log general type case}\label{LGC}

In this subsection, we shall prove Theorem \ref{logfamily} under the following additional assumptions:
\begin{enumerate}
\item[(1)] Every fiber $(X_{s},D_{s})$  over $s \in S^{\circ}$ is of log general type, i.e., $K_{X_{s}}+D_{s}$ is big,
\item[(2)] $D$ is $\mathbb{Q}$-linearly equivalent to a genuine line bundle $B$. 
\end{enumerate}
In the following proof we shall consider $D$ as a line bundle 
$\mathcal{O}_{X}(B)$ and we shall abuse the notation 
$\mathcal{O}_{X}(D)$ instead of $\mathcal{O}_{X}(B)$, i.e., we shall 
fix a line bundle structure associated with the $\mathbb{Q}$-divisor $D$.   

Let $f : X \to S$ and $D$ be as in Theorem \ref{logfamily}. 
The construction of $\hat{h}_{can}$ in Theorem \ref{logfamily} is 
similar to that in Theorem \ref{family}.   
Here we shall assume that $S$ is of dimension $1$ for simplicity.  
The case of $\dim S > 1$ is treated parallel to Section \ref{dimgeq1}. 
The construction of $\hat{h}_{can}$ on the family is similar to 
Theorem \ref{family}.  More precisely we  replace $E_{m}$ in (\ref{em}) by 
\begin{equation}\label{efem}
E_{m} := f_{*}\mathcal{O}_{X}(A+ m(K_{X/S}+D)),  
\end{equation}
where $A$ is a sufficiently ample line bundle on $X$. 
We note that $E_{m}$ is locally free by the assumption: $\dim S = 1$. 
Let $E_{m,s}$ denotes the fiber of the vector bundle $E_{m}$ at $s$. 
Let $h_{A}$ be a $C^{\infty}$-hermitian metric on $A$ with strictly 
positive curvature on $X$. 
Let  $\sigma_{D}$ is a nonzero multivalued holomorphic section of 
$\mathcal{O}_{X}(D)$ with divisor $D$ and we set  
\begin{equation}\label{hD}
h_{D} = \frac{1}{|\sigma_{D}|^{2}}. 
\end{equation}
By fixing $\sigma_{D}$ 
we may identify  a holomorphic section $\tau$ of $\mathcal{O}_{X}(m(K_{X}+D))$ with a (multivalued) meromorphic $m$-ple canonical form $\tau /(\sigma_{D})^{m}$.
For $s\in S^{\circ}$ we set 
\begin{equation}\label{LKA}
\hat{K}^{A}_{m,s} : = \sup \{\,|\sigma|^{\frac{2}{m}}; \sigma \in E_{m,s}, 
\parallel\sigma\parallel_{\frac{1}{m}} = 1\},
\end{equation} 
where 
\begin{equation}
\parallel\sigma\parallel_{\frac{1}{m}} = \left|\int_{X_{s}}h_{A}^{\frac{1}{m}}\cdot h_{D,s}\cdot(\sigma\wedge\bar{\sigma})^{\frac{1}{m}}\right|^{\frac{m}{2}},  \end{equation} 
where $h_{D,s}:= h_{D}|X_{s}$. 
We set  
\begin{equation}
\hat{K}^{A}_{\infty,s}:= \mbox{the upper envelope of}\,\,\,\,\limsup_{m\to\infty}h_{A}^{\frac{1}{m}}\cdot\hat{K}^{A}_{m,s}. 
\end{equation}
and  define $\hat{K}^{A}_{\infty}$ by $\hat{K}^{A}_{\infty}|X_{s} = \hat{K}^{A}_{\infty,s}$. 
Then we set  
\begin{equation}
\hat{h}_{can,A} : = \mbox{the lower envelope of}\,\,\,(\hat{K}^{A}_{\infty})^{-1}. 
\end{equation}
and 
\begin{equation}
\hat{h}_{can}:= \mbox{the lower envelope of}\,\,\,\,\inf_{A}\hat{h}_{can,A},  
\end{equation} 
where $A$ runs all the ample line bundles on $X$. 
Since $h_{D}$ defined as (\ref{hD}) has semipositive curvature current, using 
Theorem \ref{NSB-P}, we see that $\hat{h}_{can}$ has semipositive curvature current on $X^{\circ}:= 
f^{-1}(S^{\circ})$ by the same argument as in the proof of Theorem \ref{family}
 in Section \ref{semipos}.  And we can  extend $\hat{h}_{can}$ to a singular hermitian metric on $K_{X/S}+D$ over the whole $X$ just as in Section \ref{semipos}. 
Hence we only need to prove the assertions (2) and (3) in Theorem \ref{logfamily}. 

For every $s\in S^{\circ}$, we also define the canonical singular hermitian metric  $\hat{h}_{can,s}$ 
on $K_{X_{s}}+D_{s}$ as in Lemma \ref{logAZD}.   We note that $\hat{h}_{can}|X_{s}$ may be different from  $\hat{h}_{can,s}$ for some $s\in S^{\circ}$. 
To prove Theorem \ref{logfamily}, we need to compare $\hat{h}_{can}|X_{s}$ 
 with  $\hat{h}_{can,s}$.  \vspace{3mm} \\ 
Let us fix $s \in S^{\circ}$.   
Let $U$ be a neighborhood of $s\in S^{\circ}$ (in $S^{\circ}$) which is biholomorphic to the unit open disk $\Delta$ in $\mathbb{C}$ by a local coordinate $s$.  Hereafter over $f^{-1}(U)$, 
we identify $K_{X/S}|U$ with $K_{X}|U$ by the isomorphism:
\begin{equation}\label{ID}
\otimes\,f^{*}ds : K_{X/S}|\,U\to 
K_{X}|\,U.
\end{equation}
   
First we shall assume that $K_{X_{s}}+D_{s}$ is big, i.e., $(X_{s},D_{s})$ is 
of log general type. The general case follows from this special case 
by considering $K_{X/S}+D + \epsilon A$ instead of $K_{X/S}+D$ and letting 
$\epsilon$ tend to $0$.  We shall discuss  in detail later.      

By \cite{b-c-h-m} there exists a modification 
\[
\mu_{s} : Y_{s}\to X_{s}
\]
such that $\mu_{s}^{*}(K_{X_{s}}+D_{s})$ has a Zariski decomposition:  
\begin{equation}\label{ZD}
\mu_{s}^{*}(K_{X_{s}}+D_{s}) = P_{s} + N_{s}, 
\end{equation}
i.e., $P_{s},N_{s}\in \mbox{Div}(Y_{s})\otimes \mathbb{Q}$ such that 
\begin{enumerate}
\item[(1)] $P_{s}$ is nef, 
\item[(2)] $N_{s}$ is effective, 
\item[(3)] $H^{0}(X_{s},\mathcal{O}_{X_{s}}(\lfloor m(K_{X_{s}}+D_{s})\rfloor))\simeq H^{0}(Y_{s},\mathcal{O}_{Y_{s}}(\lfloor mP_{s}\rfloor))$
 for every $m\geqq 0$. 
\end{enumerate}
We note that in this case $P_{s}$ is semiample (see \cite{b-c-h-m}).
\begin{lemma}\label{semiample} 
Let $h_{P_{s}}$ be a $C^{\infty}$-hermtian metric on $P_{s}$ with semipositive curvature and let 
$\tau_{N}$ be a multivalued holomorphic section of $N_{s}$ with 
divisor $N_{s}$.  Then 
\[
h_{P_{s}}\cdot \frac{1}{|\tau_{N}|^{2}}
\]
is an AZD on $\mu_{s}^{*}(K_{X_{s}}+D_{s})$ with minimal singularities.
\fbox{}
\end{lemma}
{\em Proof.} By Kodaira's lemma 
$H^{0}(Y_{s},\mathcal{O}_{Y_{s}}(m_{0}P_{s} - \mu_{s}^{*}A)) \neq  0$
holds for a sufficiently large positive integer $m_{0}$ with $m_{0}P_{s}$ is Cartier.  We take a nonzero element 
$\sigma_{0}\in H^{0}(X_{s},\mu_{s,*}\mathcal{O}_{Y_{s}}(m_{0}P_{s} - \mu_{s}^{*}A))$ and identify $\sigma_{0}$ as an element of 
$H^{0}(X_{s},\mathcal{O}_{X_{s}}(m(K_{X_{s}}+D_{s})-A))$ in the natural way.
Hence we have the inclusion:
\[
\otimes \,\mu_{s}^{*}\sigma_{0}\, : \mathcal{O}_{Y}(\mu_{s}^{*}(A +m(K_{s}+D_{s})))
\hookrightarrow  
\mathcal{O}_{Y}(\mu_{s}^{*}((m+m_{0})(K_{X_{s}}+D_{s}))). 
\]    
Then for every element $\sigma\in H^{0}(X_{s},\mathcal{O}_{X_{s}}(A + m(K_{X_{s}}+ D_{s})))$, by H\"{o}lder's inequality we see that 
\begin{equation}\label{SHO}
\int_{X_{s}}|\sigma_{0}\cdot \sigma|^{\frac{2}{m+m_{0}}}
\cdot h_{D}
\leqq \left(\int_{X_{s}}h_{A}^{\frac{1}{m}}\cdot |\sigma|^{\frac{2}{m}}
\right)^{\frac{m}{m+m_{0}}}\!\!\cdot\,
\left(\int_{X_{s}}h_{A}^{-\frac{1}{m_{0}}}\cdot |\sigma_{0}|^{\frac{2}{m_{0}}}
\right)^{\frac{m_{0}}{m+m_{0}}} 
\end{equation}
holds.  Now for every positive integer $\ell$, we set 
\begin{equation}
K_{\ell} :=  \sup \left\{\,\mid \sigma\mid^{\frac{2}{\ell}} ;
\sigma \in \Gamma (X_{s},{\cal O}_{X}(\ell(K_{X_{s}}+D_{s}))), \left|\int_{X_{s}}h_{D,s}\cdot(\sigma\wedge\bar{\sigma})^{\frac{1}{\ell}}\right| = 1\right\}.
\end{equation}
Then as Theorem \ref{canazd},  
\begin{equation}
h_{can} := \mbox{the lower envelope of}\,\,(\limsup_{\ell\to\infty}K_{\ell})^{-1}
\end{equation} 
is an AZD on $K_{X_{s}}+D_{s}$.
Then by (\ref{SHO}) we see that 
\[
K_{m+m_{0}} \geqq |\sigma_{0}|^{\frac{2}{m+m_{0}}}\cdot\left(\hat{K}_{m,A}\right)^{\frac{m}{m+m_{0}}}\cdot\left(\int_{X_{s}}h_{A}^{-\frac{1}{m_{0}}}\cdot |\sigma_{0}|^{\frac{2}{m_{0}}}
\right)^{-\frac{m_{0}}{m+m_{0}}} 
\]
holds.  Letting $m$ tend to infinity,  we see that 
\begin{equation}\label{QI}
h_{can} \leqq \hat{h}_{can,A}
\end{equation}  
holds.  By the definition of the Zariski decomposition and the semiampleness 
of $P_{s}$, we see that 
$h_{P_{s}}\cdot \frac{1}{|\tau_{N}|^{2}}$ is quasi-isometric to $h_{can}$,i.e.,
the ratio of these metrics is pinched by positive constants.    
Hence by (\ref{QI}) $h_{P_{s}}\cdot \frac{1}{|\tau_{N}|^{2}}$ is an AZD with 
minimal singularities on $\mu_{s}^{*}(K_{X_{s}}+D_{s})$.  \fbox{}\vspace{3mm}\\

Let $a$ be a positive integer such that $aP_{s}\in \mbox{Div}(Y_{s})$. 
Let us take the ample line bundle $A$ so that 
for every pseudoeffective singular hermitian line bundle $(L,h_{L})$ 
on $X_{s}$ 
\[
\mathcal{O}_{X_{s}}(A|X_{s}+j(K_{X_{s}}+D_{s})+L)\otimes \mathcal{I}(h_{L})
\]
is globally generated over $X_{s}$ for every $0\leqq j\leqq a$. 
This is certainly possible (see Proposition \ref{suff} in  Section 5.1).

Let us fix a $C^{\infty}$-hermitian metric $h_{ref,s}$ on $K_{X_{s}}+D_{s}$. 
The following lemma is similar to Lemma \ref{ext}. 
Hereafter we shall denote $A|X_{s}$ (resp. $h_{A}|X_{s}$)   by 
$A_{s}$ (resp. $h_{A,s}$) for simplicity.   

\begin{lemma}\label{ext2}
Let $h_{D}$ and the positive integer $a$ be as above. 
If we take $A$ sufficiently ample as above, we have the followings. 
\begin{enumerate}
 
\item[(1)] For every  positive integer $m$, every element of \\ 
$\Gamma (X_{s},{\cal O}_{X_{s}}(A_{s}+m(K_{X_{s}}+D_{s}))
\otimes {\cal I}(h_{D,s}\cdot\hat{h}_{can,s}^{a\lfloor m-1/a\rfloor}))$ 
extends to an element of 
  $\Gamma \left(f^{-1}(U),{\cal O}_{X}(A+m(K_{X}+D))\right)$
, where $h_{D,s}$ denotes the restriction $h_{D}|X_{s}$.
\item[(2)] There exists a positive constant $C$ such that 
\[
\hat{h}_{can}|X_{s} \leqq C\cdot\hat{h}_{can,s}
\]
holds on $X_{s}$. In particular $\hat{h}_{can}|X_{s}$ is an AZD on 
$K_{X_{s}}+D_{s}$ with minimal singularities (cf. Definition \ref{minAZD}).  \fbox{}
\end{enumerate} 
\end{lemma}

\noindent{\em Proof of Lemma \ref{ext2}}. 
We prove the lemma by induction on $m$. 
If $m=1$, then the $L^{2}$-extension theorem (\cite{o-t,o}) 
implies that every element of 
$\Gamma (X_{s},{\cal O}_{X_{s}}(A_{s} + (K_{X_{s}}+D_{s}))\otimes\mathcal{I}(h_{D,s})) = \Gamma (X_{s},{\cal O}_{X_{s}}(A + (K_{X_{s}}+D_{s})))$ (because $(X_{s},D_{s})$ is KLT)  
extends to an element of $\Gamma (f^{-1}(U),{\cal O}_{X}(A + (K_{X}+D))$.

Suppose that the extension is settled for $m-1 (m\geqq 2)$. 
 Let $\{\sigma_{1,s}^{(m-1)},\cdots ,\sigma_{N(m-1)}^{(m-1)}\}$
be a basis of  $\Gamma (X_{s},{\cal O}_{X_{s}}(A_{s}+ (m-1)(K_{X_{s}}+D_{s})
\otimes {\cal I}(h_{D}\cdot\hat{h}_{can}^{a\lfloor(m-2)/a\rfloor}))$. 
By the inductive assumption, we have already constructed holomorphic extensions:
\[
\{\tilde{\sigma}_{1,s}^{(m-1)},\cdots ,\tilde{\sigma}_{N(m-1),s}^{(m-1)}\}
\subset \Gamma (f^{-1}(U),{\cal O}_{X}(A+(m-1)(K_{X}+D)))
\]
of   $\{\sigma_{1,s}^{(m-1)},\cdots ,\sigma_{N(m-1),s}^{(m-1)}\}$
to $f^{-1}(U)$. 
We define the singular hermitian metric $\tilde{h}_{m-1}$ on $A+ (m-1)(K_{X} +D))|f^{-1}(U)$ by 
\begin{equation}\label{hm-1}
\tilde{h}_{m-1} :=
 \frac{1}{\sum_{j=1}^{N(m-1)}|\tilde{\sigma}^{(m-1)}_{j,s}|^{2}}.
\end{equation}
By the choice of $a$ and the fact that $(X_{s},D_{s})$ is KLT, $aP_{s}$ is integral on $Y_{s}$  and 
we have the inclusion: 
\begin{equation}\label{m-1infty}
{\cal O}_{X_{s}}\left(A_{s} + \left(m-1 - a\left\lfloor\frac{m-2}{a}\right\rfloor\right)(K_{X_{s}}+D_{s})\right) \otimes (\mu_{s})_{*}\mathcal{O}_{Y_{s}}\left(a\left\lfloor\frac{m-2}{a}\right\rfloor\cdot P_{s}\right)
\end{equation}
\[
\hspace{30mm}\hookrightarrow {\cal O}_{X_{s}}(A+ (m-1)(K_{X_{s}}+D_{s}))\otimes  {\cal I}(h_{D,s}\cdot\hat{h}_{can,s}^{a\lfloor m-2/a\rfloor}).
\] 
In fact by Lemma \ref{semiample}, $(\mu_{s})_{*}\mathcal{O}_{Y_{s}}\left(a\left\lfloor\frac{m-2}{a}\right\rfloor\cdot P_{s}\right)$ is nothing but the sheaf of germs of locally bounded holomorphic sections of $A + a\lfloor (m-2)/a\rfloor\cdot(K_{X_{s}}+D_{s})$ 
with respect to the metric $\hat{h}_{can,s}^{a\lfloor (m-2)/a\rfloor}$. 
Then since $(X_{s},D_{s})$ is KLT and ${\cal O}_{X_{s}}(A_{s}+ (m-1)(K_{X_{s}}+D_{s}))\otimes  {\cal I}(h_{D,s}\cdot\hat{h}_{can,s}^{a\lfloor (m-2)/a\rfloor})$ is globally generated over $X_{s}$ 
(by the choice of $A$), (\ref{m-1infty}) implies that
\begin{equation}\label{emel}
\tilde{h}_{m-1}|X_{s} = O (h_{A,s}\cdot\hat{h}_{can,s}^{a\lfloor (m-2)/a\rfloor}\cdot h_{ref,s}^{a\{ (m-2)/a\}}). 
\end{equation}
holds on $X_{s}$.  Here we have neglected the effect of the singularities $h_{D,s}$ by the choice of $a$ and $A$.   Apparently $\tilde{h}_{m-1}$ has a semipositive curvature current on $U$. 
Hence $\tilde{h}_{m-1}\cdot h_{D}$ 
is a singular hermitian metric on $(A + m(K_{X} + D)-K_{X})|U$ with 
semipositive curvature current.  
Then by the $L^{2}$-extension theorem (\cite{o-t,o}), we may extend every element of 
\[
\Gamma(X_{s},{\cal O}_{X_{s}}(A_{s}+m(K_{X_{s}}+D_{s}))\otimes\mathcal{I}(\tilde{h}_{m-1}\cdot  h_{D}|X_{s}))
\]
to an element of 
\[
\Gamma(f^{-1}(U),{\cal O}_{X}(A+m(K_{X}+D))\otimes\mathcal{I}(\tilde{h}_{m-1}\cdot  h_{D})). 
\]
And by (\ref{emel}), we have that
\begin{equation}\label{m-2}
\tilde{h}_{m-1}\cdot h_{D}|X_{s} = O(h_{A,s}\cdot\hat{h}_{can,s}^{a\lfloor(m-2)/a\rfloor}\cdot h_{ref,s}^{a\{(m-2)/a\}}\cdot h_{D,s}).
\end{equation} 
We note that  by  the choice of  $A$,
$\mathcal{O}_{X_{s}}(A_{s} + m(K_{X_{s}}+D_{s}))\otimes \mathcal{I}(\hat{h}_{can,s}^{a\lfloor(m-1)/a\rfloor}\cdot h_{D,s})$ is globally generated over $X_{s}$.  
Hence by (\ref{m-2}) we may extends every element of 
$H^{0}(X_{s},\mathcal{O}_{X_{s}}(A_{s}+ m(K_{X_{s}}+D_{s}))\otimes \mathcal{I}(\hat{h}_{can,s}^{a\lfloor(m-1)/a\rfloor}\cdot h_{D,s}))$ to an element ot 
$\Gamma (f^{-1}(U),\mathcal{O}_{X}(A + m(K_{X}+D)))$ and the estimate  
\begin{equation}\label{mest}
\tilde{h}_{m}|X_{s} = O (h_{A,s}\cdot\hat{h}_{can,s}^{a\lfloor (m-1)/a\rfloor}\cdot h_{ref,s}^{a\{(m-1)/a\}})
\end{equation}
holds by the same argument as above.    Hence by induction on $m$, we see that (\ref{mest}) holds for every 
$m\geqq 1$.   
(\ref{mest}) implies the inclusion:
\begin{equation}\label{mestinclusion}
H^{0}(X_{s},\mathcal{O}_{X_{s}}(A_{s} + m(K_{X_{s}}+D_{s}))\otimes\mathcal{I}(\hat{h}_{can,s}^{a\lfloor (m-1)/a\rfloor}))
\hookrightarrow E_{m,s} 
\end{equation}
and the assertion (1) of Lemma \ref{ext2} holds. 

By using H\"{o}lder's inequality as (\ref{b}) and the trivial inequality:
$a\lfloor (m-1)/a\rfloor\leqq m-1$, we shall transform the inclusion (\ref{mestinclusion}) to the lower estimates of $\hat{K}^{A}_{m}|X_{s}$ (cf. (\ref{LKA})):
\begin{equation}\label{KMS}
\hat{K}^{A}_{m}|X_{s} 
\geqq K\left(A_{s}+m(K_{X_{s}}+D_{s}),h_{A}\cdot\hat{h}_{can,s}^{m-1}\cdot h_{D,s}\right)^{\frac{1}{m}}\cdot  \left(\int_{X_{s}}\hat{h}_{can,s}^{-1}\cdot h_{D,s}\right)^{-\frac{m-1}{m}}. 
\end{equation}
We note that  by Lemma \ref{dem} and Remark \ref{DDDD},  we see that
\[
\limsup_{m\to\infty}h_{A,s}^{\frac{1}{m}}\cdot K(A_{s}+m(K_{X_{s}}+D_{s}),h_{A,s}\cdot
\hat{h}_{can,s}^{m-1}\cdot h_{D,s})^{\frac{1}{m}} =  
\hat{h}_{can,s}^{-1}
\]
holds.  
Then by the H\"{o}lder inequality, similarly as (\ref{holder}) and (\ref{b}) we have the estimate: 
\begin{equation}\label{estk}
\limsup_{m\to\infty}h_{A,s}^{\frac{1}{m}}\cdot \hat{K}^{A}_{m}|X_{s}
\geqq \left(\int_{X_{s}}\hat{h}_{can,s}^{-1}\cdot h_{D,s}\right)^{-1}\cdot
\hat{h}_{can,s}^{-1}.
\end{equation}
Hence by setting 
\[
C:= \int_{X_{s}}\hat{h}_{can,s}^{-1}\cdot h_{D,s},
\]
we have the estimate: 
\begin{equation}\label{dom}
\hat{h}_{can}|X_{s} \leqq  C\cdot\hat{h}_{can,s}
\end{equation}
This completes the proof of Lemma \ref{ext2}. \fbox{}
\vspace{3mm} \\ 

By Lemma \ref{ext2} and Theorem \ref{main2}, we see that $\hat{h}_{can}|X_{s}$ is an AZD on 
$K_{X_{s}}+D_{s}$ with minimal singularities. 
Since $s \in S^{\circ}$ is arbitrary, we complete the proof of the assertion (2) of Theorem \ref{logfamily}.  

Now we shall prove the assertion (3) of Theorem \ref{logfamily}. 
For $\ell,m \geqq 1$, we set $E^{(\ell)}_{m}:= f_{*}\mathcal{O}_{X}(\ell A + m(K_{X/S}+D))$.  We note that $E^{(\ell)}_{m}$ is a vector bundle, because we have assumed 
that $\dim S = 1$. 
We set 
\begin{equation}\label{ef}
F:= \{ s\in S^{\circ}|
\mbox{$E^{(\ell)}_{m,s} \neq \Gamma (X_{s},\mathcal{O}_{X_{s}}(\ell A_{s}+m(K_{X_{s}}+D_{s}))$ 
 for some $\ell,m\geqq 1$}\},   
\end{equation}
where $E^{(\ell)}_{m,s}$ denotes the fiber at $s$. 
Then by the definitions of $\hat{h}_{can,s}$ and $\hat{h}_{can}$, we see 
that for every $s\in S^{\circ}\backslash F$, 
\[
\hat{h}_{can}|X_{s} \leqq \hat{h}_{can,s}
\]
holds on $X_{s}$.   To prove that $F$ is empty, we need to prove 
that $h^{0}(X_{s},\mathcal{O}_{X_{s}}(\ell A_{s} + m(K_{X_{s}}+D_{s})))$ is locally 
constant on $S^{\circ}$ for every $\ell,m \geqq 1$.  But this follows from Theorem \ref{logplurigenera}.  Here we  note that to prove Theorem \ref{logplurigenera}.  
we need to use only the assertions (1) and (2) of Theorem \ref{logfamily}.   
 See Section \ref{INV} below.
 
\subsection{Proof of Theorem \ref{logfamily}; Non log general type case}\label{GC}
Next we shall consider the case that  $K_{X_{s}}+D_{s}$ is 
not necessarily big and $D$ is $\mathbb{Q}$-linearly equivalent to a genuine line bundle $B$.   In this case we shall consider  the perturbation: 
$K_{X/S}+D + \ell^{-1}A (\ell=1,2,\cdots)$. Since $K_{X/S}+D+\ell^{-1}A|X_{s}$ is big for every $s\in S^{\circ}$ by the assumption, 
we may apply the argument in the last subsection.  But there is a 
minor difference that $K_{X/S}+D+\ell^{-1}A|X_{s}$  is not Cartier.  
This is by no means  an essential difficulty. But we need to modify 
the argument in an obvious way.    Then we let  $\ell$ tend to infinity.

More precisely we argue as follows. Let $G = b(\ell)A$ be a positive multiple of $A$.     
We set 
\begin{equation}
E_{\ell,m}:= f_{*}\mathcal{O}_{X}(G + m(K_{X/S}+D)+\lfloor m/\ell\rfloor A)
\end{equation}
and 
\begin{equation}
\hat{K}^{G}_{\ell,m}(s) := \sup \left\{\,|\sigma|^{\frac{2}{m}};\,\, \sigma\in E_{\ell,m,s},\, \parallel\sigma\parallel_{\ell,m,s} = 1\right\}, 
\end{equation}
where 
\begin{equation}
\parallel\sigma\parallel_{\ell,m,s}:= \left|\int_{X_{s}}h_{A,s}^{\frac{b(\ell)+\lfloor m/\ell\rfloor}{m}}\cdot h_{D,s}\cdot(\sigma\wedge\bar{\sigma})^{\frac{1}{m}}\right|^{\frac{m}{2}}.
\end{equation}
We set 
\begin{equation}
\hat{K}^{G}_{\ell,\infty}:= \mbox{the upper envelope of}\,\,\, \limsup_{m\to\infty}\hat{K}^{G}_{\ell,m} 
\end{equation}
and 
\begin{equation}
\hat{h}_{can,\ell,G}:= \frac{1}{\hat{K}^{G}_{\ell,\infty}}. 
\end{equation}
And  we set 
\begin{equation}
\hat{h}_{can,\ell}:= \mbox{the lower envelope of}\,\,\,\inf_{G}h_{can,\ell,G},
\end{equation}
where $G$ runs all the positive multiples of $A$. 
Replacing $f: X\to S$ by $X_{s}\to \{ s\}$ we obtain 
$\hat{h}_{can,\ell,G,s}$ and $\hat{h}_{can,\ell,s}$. 
By \cite{b-c-h-m} there exists a modification 
\[
\mu_{\ell,s} : Y_{\ell,s}\to X_{s}
\]
such that $\mu_{\ell,s}^{*}(K_{X_{s}}+D_{s}+\ell^{-1}A)$ has a Zariski decomposition:  
\begin{equation}\label{ZDk}
\mu_{k,s}^{*}(K_{X_{s}}+D_{s}+\frac{1}{\ell}A) = P_{\ell,s} + N_{\ell,s}.  
\end{equation}
Let $a= a(\ell)$ be a positive integer such that $a\cdot P_{\ell,s}$ is Cartier. Here we shall use the same notation as before for simplicity. 

\begin{lemma}\label{ext3}
If we take $G$ (depending on $\ell$) sufficiently ample, we have the followings.
\begin{enumerate}
\item[(1)] For every  positive integer $m$, every element of \\
$\Gamma (X_{s},{\cal O}_{X_{s}}(G_{s}+  m(K_{X_{s}}+D_{s})+\lfloor m/\ell\rfloor A))
\otimes {\cal I}(h_{D,s}\cdot\hat{h}_{can,\ell,s}^{a\lfloor (m-1)/a\rfloor}))$ 
extends to an element of 
$\Gamma \left(f^{-1}(U),{\cal O}_{X}(G+m(K_{X}+D)+\lfloor m/\ell\rfloor A)\right)$
, where $G_{s} := G|X_{s}$ and $h_{D,s}$ denotes the restriction $h_{D}|X_{s}$. 
\item[(2)] There exists a positive constant $C$ independent of $\ell$ such that \begin{equation}\label{unifell}
\hat{h}_{can,\ell}|X_{s} \leqq C\cdot h_{A,s}^{\frac{1}{\ell}}\cdot\hat{h}_{can,s}\end{equation}
holds on $X_{s}$.
 \fbox{}
\end{enumerate} 
\end{lemma}
{\em Proof}.  The proof of the assertion (1)  is parallel to that of 
Lemma \ref{ext2}, i.e., we use the successive extensions. 
The only difference here is that  we need to tensorize $(A,h_{A})$ every $\ell$-steps.  Then as in Lemma \ref{ext2}, we see that  
we have that 
\begin{equation}\label{doml}
\hat{h}_{can,\ell}|X_{s} = O(\hat{h}_{can,\ell,s})
\end{equation}
holds as (\ref{dom}). 
  
Let us  prove the assertion (2).  By the assertion (1) we see that 
since $\hat{h}_{can,\ell}|X_{s}$ is an AZD with minimal singularities $\ell^{-1}A + (K_{X_{s}}+D_{s})$  by (\ref{doml}) and $\hat{h}_{can,s}$ is 
an AZD with minimal singularities on $K_{X_{s}} + D_{s}$.  
Since $A$ is ample  it is clear that
\begin{equation}\label{triv}
\hat{h}_{can,\ell}|X_{s} \leqq O\left(h_{A,s}^{\frac{1}{\ell}}\cdot\hat{h}_{can,s}\right)
\end{equation}
holds on $X_{s}$.  Our task is to find a constant $C$ independent of 
$\ell$ such that (\ref{unifell}) holds.  By (\ref{triv}) and the assertion (1), we have the inclusion:
\begin{equation}\label{I}
H^{0}(X_{s},\mathcal{O}_{X_{s}}(G_{s}+m(K_{X_{s}}+D_{s}) + \lfloor m/\ell\rfloor A )
\otimes \mathcal{I}(\hat{h}_{can,s}^{a\lfloor (m-1)/a\rfloor}\cdot h_{D,s}))
\hookrightarrow  E_{\ell,m,s}.
\end{equation}
Then  we tranform the inclusion (\ref{I}) to the inequality: 
\begin{equation}
\hat{K}^{G}_{\ell,m}|X_{s} 
\geqq 
 K\left(G_{s}+m(K_{X_{s}}+D_{s}) + \lfloor m/\ell\rfloor A_{s} ,h_{A,s}^{b(\ell)+\lfloor m/\ell\rfloor}\cdot\hat{h}_{can,s}^{m-1}\cdot h_{D,s}\right)^{\frac{1}{m}}\cdot  \left(\int_{X_{s}}\hat{h}_{can,s}^{-1}\cdot h_{D,s}\right)^{-\frac{m-1}{m}}
\end{equation}
obtained just as (\ref{KMS}) above. 
Hence letting $m$ tend to infinity, by Lemma \ref{dem}, 
\[
\hat{h}_{can,\ell}|X_{s} \leqq C\cdot\left(h_{A,s}^{\frac{1}{\ell}}\cdot\hat{h}_{can,s}\right)
\]
holds for
\[
C:=  \int_{X_{s}}\hat{h}_{can,s}^{-1}\cdot h_{D,s},
\]
hence it is independent of $\ell$. Hence the assertion (2) holds.    
\fbox{} \vspace{3mm} \\ 

\noindent Lemma \ref{ext3} implies that 
\begin{equation}
\hat{h}_{can,\infty}:= \mbox{the lower envelope of}\,\,\,\,\liminf_{\ell\to\infty}\hat{h}_{can,\ell} 
\end{equation}
exists and $\hat{h}_{can,\infty}|X_{s}$ is an  AZD on $K_{X_{s}}+D_{s}$ with 
minimal singularities (cf. Definition \ref{minAZD}).  In fact the upper estimate  of $\hat{h}_{can,\infty}$ follows from the assertion (1) of Lemma \ref{ext3} 
and the lower estimate follows from the same argument as in Section \ref{up}. 

Hence if $A$ is sufficiently ample,  by the $L^{2}$-extension theorem (\cite{o-t,o}), we see that for every $m\geqq 0$
every element of $H^{0}(X_{s},\mathcal{O}_{X_{s}}(A_{s} + m(K_{X_{s}}+D_{s}))\otimes \mathcal{I}(\hat{h}_{can,\infty}^{m-1}\cdot h_{D}|X_{s}))$ extends to an element of 
$H^{0}(f^{-1}(U),\mathcal{O}_{X}(A + m(K_{X/S} + D)))$.    
Since $\hat{h}_{can,\infty}|X_{s}$ is an AZD on $K_{X_{s}}+D_{s}$ with minimal singularities as above, the inclusion:   
\begin{equation}\label{INC}
\hspace{-20mm} H^{0}(X_{s},\mathcal{O}_{X_{s}}(A_{s} + m(K_{X_{s}}+D_{s}))\otimes \mathcal{I}(\hat{h}_{can,s}^{m-1}\cdot h_{D,s}))\subseteq 
\end{equation}
\[
\hspace{40mm} H^{0}(X_{s},\mathcal{O}_{X_{s}}(A_{s} + m(K_{X_{s}}+D_{s}))\otimes \mathcal{I}(\hat{h}_{can,\infty}^{m-1}\cdot h_{D}|X_{s}))
\] 
holds.  Hence by the $L^{2}$-extension theorem (\cite{o-t,o}), we may extend every element of 
$H^{0}(X_{s},\mathcal{O}_{X_{s}}(A_{s} + m(K_{X_{s}}+D_{s}))\otimes \mathcal{I}(\hat{h}_{can,s}^{m-1}\cdot h_{D,s}))$ to an element of $H^{0}(f^{-1}(U),\mathcal{O}_{X}(A + m(K_{X/S} + D)))$.    
Then  as (\ref{KMS}) we may transform the inclusion (\ref{INC}) to the inequality: 
\[
\hat{K}^{A}_{m}|X_{s} 
\geqq 
 K\left(A_{s}+m(K_{X_{s}}+D_{s}),h_{A}\cdot\hat{h}_{can,s}^{m-1}\cdot h_{D,s}\right)^{\frac{1}{m}}\cdot  \left(\int_{X_{s}}\hat{h}_{can,s}^{-1}\cdot h_{D,s}\right)^{-\frac{m-1}{m}},
\]
holds and repeating the same estimate as above (see (\ref{estk})), 
letting $m$ tend to infinity,  by Lemma \ref{dem}, we see that 
\[
\hat{h}_{can}|X_{s} \leqq \left(\int_{X_{s}}\hat{h}_{can,s}^{-1}\cdot h_{D,s}\right)\cdot \hat{h}_{can,s}. 
\]
holds on $X_{s}$ and $\hat{h}_{can}|X_{s}$ is an AZD on $K_{X_{s}}+D_{s}$ 
with minimal singularities.   
Since $s\in S^{\circ}$ is arbitrary, we complete the proof of the assertion (2) of Theorem \ref{logfamily}.  
The rest of the proof (the proof of the assertion (3)) is completely parallel to that of the previous 
subsection.   We complete the proof of Theorem \ref{logfamily}, 
assuming the boundary $D$ is $\mathbb{Q}$-linearly equivalent to 
a Cartier divisor.

\subsection{Dynamical systems of singular hermitian metrics}\label{Dynamical}
In this subsection,  we complete the proof of Theorem \ref{logfamily}.
Here we do not assume that the boundary $B$ is $\mathbb{Q}$-linearly equivalent to a genuine line bundle.  In Section \ref{T2}, we also give an alternative proof 
by using the ideas in \cite{e-p}. 

First we shall prove the following theorem. The technique used here 
is essentially the same as the Ricci iteration in \cite{ricci}.  

\begin{theorem}\label{PE}
Let $f : X \longrightarrow S$ be a proper surjective projective morphism 
between  complex manifolds with connected fibers 
and let $D$ be an effective $\mathbb{Q}$-divisor on $X$  such that 
the set: 
\[
S^{\circ} := \{ s\in S|\,\,\mbox{$f$ is smooth over $s$ and $(X_{s},D_{s})$ is  KLT}\,\,\}
\]
is nonempty. 
Suppose that $S$ is connected and for some $s_{0} \in S^{\circ}$, $K_{X_{s_{0}}}+D_{s_{0}}$ is pseudoeffective. 

Then  the followings hold. 
\begin{enumerate}
\item[(1)] $K_{X_{s}}+D_{s}$ is pseudoeffective for every $s \in S^{\circ}$. 
\item[(2)] $K_{X/S} + D$ is pseudoeffective (cf. Definition \ref{pe}). \fbox{}
\end{enumerate}   
\end{theorem}
\begin{remark}
Here the pseudoeffectivity is defined as Definition \ref{pe}.  Hence we do not 
assume the compactness of the base space $S$. \fbox{}
\end{remark}
{\em Proof.} The proof is quite similar to that of Lemma \ref{ext2}. 
The only essential difference is that we need the double induction, because 
of $D_{s}$ is not $\mathbb{Q}$-linearly equivalent to a Cartier divisor. 
First we may and do assume that $\dim S = 1$ without loss of generality. 
Let $A$ be a sufficiently ample line bundle on $X$  such that 
\begin{equation}\label{L0}
L_{0} : = A + (q-1)(K_{X/S}+B)
\end{equation}
is ample.    Let us fix a $C^{\infty}$-hermitian metric 
on $h_{L_{0}}$ on  $L_{0}$ with strictly positive curvature.  
 Then we define a singular hermitian metric $K_{X/S}+B + L_{0}|X^{\circ} 
= A + q(K_{X/S}+B)|X^{\circ}$ by    
\begin{equation}\label{H0}
h_{0}:= \hat{h}_{can}(L_{0}+B,h_{L_{0}}\cdot h_{D}). 
\end{equation}
Here we have used the relative version of Theorem \ref{KLTAZD}, i.e., 
we take the direct image $f_{*}\mathcal{O}_{X}(m(K_{X/S}+B+L_{0}))$ 
for every sufficiently divisible $m > 0$ and construct the metric 
just as in Thorem \ref{family*} by using the similar construction as in Theorem \ref{KLTAZD}. 
Then $h_{0}$ is of semipositive curvature current over $X^{\circ}$ by using 
Theorem \ref{NSB-P} as in Section \ref{semipos} and it extends to  a singular hermitian metric with semipositive curvature on $q(K_{X/S} + B) + A$ by the same argument as in Section \ref{dimgeq1}.    
Now we set 
\begin{equation}
L_{1} : = (q-1)(K_{X/S} + B) + \frac{q-1}{q}A    
\end{equation}
and define the singular hermitian metric $h_{1}$ on 
$K_{X/S} + B + L_{1}$ over $X^{\circ}$ by 
\begin{equation}
h_{1}:= \hat{h}_{can}(L_{1}+B,h_{0}^{\frac{q-1}{q}}\cdot h_{D}).  
\end{equation}
Similarly as $h_{0}$, $h_{1}$ is of semipositive curvature current over $X^{\circ}$ and it extends to a singular hermitian metric 
on $K_{X/S}+B + L_{1}$ with semipositive curvature current. 
Inductively for every positive integer $m$, we set 
\begin{equation}\label{LM}
L_{m} : = (q-1)(K_{X/S} + B) + \left(\frac{q-1}{q}\right)^{m}\!\!A    
\end{equation}
and 
\begin{equation}
h_{m}:  = \hat{h}_{can}(L_{m}+B,h_{m-1}^{\frac{q-1}{q}}\cdot h_{D}).
\end{equation}
Then by induction on $m$, using Theorem \ref{NSB-P}, we see that $h_{m}$ has semipositive curvature for every $m$. 
The above inductive construction is not well defined apriori, since we do not 
assume the pseudoeffectivity of $(K_{X/S} + B)|X_{s}$ for every $s\in S^{\circ}$. But the well definedness of $h_{m}$ can be verified by successive extensions as follows.

Let $U$ be a neighborhood of $s_{0}$ in $S^{\circ}$ which is biholomorphic to the  unit disk $\Delta$ in $\mathbb{C}$.  We may assume that $s_{0} = 0$ 
on $U \simeq \Delta$. 
Let us assume the followings: 
\begin{enumerate}
\item[(1)] We have already defined the singular hermitian metric $h_{m-1}$ on  $K_{X/S} + B +L_{m-1}$ with semipositive curvature.  
\item[(2)] $h_{m-1}|X_{0}$ is an AZD of $K_{X_{0}}+B_{0} + L_{m-1}|X_{0}$ with minimal singularities. 
\end{enumerate}
These assumptions are certainly satisfied for $m-1 = 0$, if we take $A$ sufficiently ample (see (\ref{L0}) and (\ref{H0})).  
Under these assumptions, we shall prove the followings:  
\begin{enumerate} 
\item[$\mbox{(A0)}_{m}$] $K_{X_{s}}+ B_{s}+ L_{m}|X_{s}$ is pseudoeffective 
for every $s\in S^{\circ}$.  Hence the singular hermitian metric $h_{m}$ on $K_{X/S} + B + L_{m}$ is well defined and has semipositive curvature. 
\item[$\mbox{(B0)}_{m}$] $h_{m}|X_{0}$ is an AZD of $K_{X_{0}}+B_{0} + L_{m}|X_{0}$ with minimal singularities. 
\end{enumerate}
Let $H$ be a sufficiently ample line bundle  on $X$ in the sense of 
Proposition \ref{suff} and let 
$h_{H}$ be a $C^{\infty}$-hermitian metric on $H$ with strictly positive curvature.

  We shall construct the singular hermitian metric $\tilde{h}_{m,\ell}$  
on $H + \lfloor\ell(K_{X/S}+B + L_{m})\rfloor|f^{-1}(U)$ with semipositive curvature  for every $\ell \geqq 0$  by induction on $\ell$ as follows.
Let $h_{A}$ be a $C^{\infty}$-hermitian metric on $A$ with strictly positive curvature.  

For  $\ell = 0$, we set $\tilde{h}_{m,0} := h_{H}$.  
Suppose that we have already constructed $\tilde{h}_{m,\ell-1}$ for some $\ell\geqq 1$.    
We shall extend  every element of 
\begin{equation}\label{enough}
H^{0}(X_{0},\mathcal{O}_{X_{0}}(H +\lfloor\ell(K_{X/S}+B + L_{m})\rfloor)
\otimes \mathcal{I}(\tilde{h}_{m,\ell-1}\cdot h_{m-1}^{\frac{q-1}{q}}\cdot h_{D}|X_{0}))
\end{equation}
to an element of 
\[
H^{0}(f^{-1}(U),\mathcal{O}_{X}(H +\lfloor\ell(K_{X/S}+B + L_{m})\rfloor)\otimes \mathcal{I}(\tilde{h}_{m,\ell-1}\cdot h_{m-1}^{\frac{q-1}{q}}\cdot h_{D}))
\] 
by the $L^{2}$-extension theorem (\cite{o-t,o}). In fact we use the semipositively curved metric:  
\[
h_{A}^{\delta_{\ell}}\cdot \tilde{h}_{m,\ell-1}\cdot h_{m-1}^{\frac{q-1}{q}}\cdot h_{D}
\]
on $H +\lfloor\ell(K_{X/S}+B + L_{m})\rfloor - K_{X/S}$, 
where 
\[
\delta_{\ell}:= \left\lfloor\ell\left(\frac{q-1}{q}\right)^{m}\right\rfloor
-\left\lfloor(\ell -1)\left(\frac{q-1}{q}\right)^{m}\right\rfloor
\]
to apply the $L^{2}$-extension theorem. 
Extending a set of basis of $H^{0}(X_{0},\mathcal{O}_{X_{0}}(H+\lfloor\ell(K_{X/S}+B + L_{m})\rfloor)
\otimes \mathcal{I}(\tilde{h}_{m,\ell-1}\cdot h_{m-1}^{\frac{q-1}{q}}\cdot h_{D}|X_{0}))$ by the $L^{2}$-extension theorem, we define the singular hermitian metric $\tilde{h}_{m,\ell}$ on  $H + \lfloor\ell(K_{X/S}+B + L_{m})\rfloor|f^{-1}(U)$    with semipositive curvature  just as (\ref{hm-1}) above.    
Here we note that since $(X_{0},D_{0})$ is KLT, 
\[
H^{0}(X_{0},\mathcal{O}_{X_{0}}(H+\lfloor\ell(K_{X/S}+B + L_{m})\rfloor)
\otimes \mathcal{I}(\tilde{h}_{m,\ell-1}\cdot h_{m-1}^{\frac{q-1}{q}}\cdot h_{D}|X_{0}))
\]
contains the subspace:  
\[
H_{(\infty)}^{0}(X_{0},\mathcal{O}_{X_{0}}(H+\lfloor\ell(K_{X/S}+B + L_{m})\rfloor),(\tilde{h}_{m,\ell-1}\cdot h_{m-1}^{\frac{q-1}{q}}\cdot h_{A}^{\delta_{\ell}})|X_{0}\cdot h_{B,0})
\]
of $H^{0}(X_{0},\mathcal{O}_{X_{0}}(H+\lfloor\ell(K_{X/S}+B + L_{m})\rfloor))$ consisting the bounded holomorphic sections 
with respect to $(\tilde{h}_{m,\ell-1}\cdot h_{m-1}^{\frac{q-1}{q}}\cdot h_{A}^{\delta_{\ell}})|X_{0}\cdot h_{B,0}$, where $h_{B,0}$ is a $C^{\infty}$-hermitian metric on $B|X_{0}$. 
Let $h_{m,0,min}$ be an AZD of $K_{X_{0}}+B_{0} + L_{m}|X_{0}$ with minimal singularities (cf. Section \ref{appendix}).    
We shall use this fact for the estimate of $\tilde{h}_{m,\ell}|X_{0}$ as 
the use of (\ref{m-1infty}) in the proof of Lemma \ref{ext2}. 
We note that 
\[
h_{m-1}|X_{0} = O\left(h_{m,0,min}\cdot h_{A,0}^{\frac{1}{q}\left(\frac{q-1}{q}\right)^{m-1}}\right)
\] 
holds by the assumption that $h_{m-1}|X_{0}$ is an AZD of 
$K_{X_{0}}+B_{0} + L_{m-1}|X_{0}$ with minimal singularities. 
 Then  using \cite{b-c-h-m} again, as (\ref{mest}) in Lemma \ref{ext2}, by induction on $\ell$,   we see that 
\begin{equation}\label{EST0}
\tilde{h}_{m,\ell}|X_{0} = O\left(h_{H,0}\cdot h_{m,0,min}^{\ell}\cdot h_{A,0}^{-\{\ell\left(\frac{q-1}{q}\right)^{m}\}}\right) 
\end{equation}
holds, where $h_{H,0}:= h_{H}|X_{0}$ (Actually as in Lemma \ref{ext2}, we have a slightly better estimate).  Hence by the sufficiently ampleness of $H$, $\{ h_{m,\ell}\}_{\ell=0}^{\infty}$ is well defined on $U$. 
In particular $K_{X/S} + B + L_{m}$ is pseudoeffective on $f^{-1}(U)$ and $h_{m}$ is well defined on $f^{-1}(U)$ because the pseudoeffectivity on the fiber 
is closed under specialization over $S^{\circ}$.  
We transform  (\ref{EST0}) into the estimate:
\begin{equation}
h_{m}|X_{0} = O(h_{m,0,min})
\end{equation} 
as (\ref{dom}) by the same argument as (\ref{mest}) $\sim$ (\ref{dom}). 
And we see that $h_{m}|X_{0}$ is an AZD with minimal singularities on $K_{X_{0}}+B_{0} + L_{m}|X_{0}$.  Hence the induction works. 
 In this way, $\{\tilde{h}_{m,\ell}\}_{\ell=0}^{\infty}$ is well defined for every $m\geqq 0$.  And this implies that  $(K_{X/S}+ B + L_{m})|f^{-1}(U)$ is pseudoeffective for every $m\geqq 0$.   Letting $m$ tend to infinity, we see that $K_{X_{s}} + B_{s}$ is 
 pseudoeffective for every $s\in U$.   
The openness of the pseudoeffectivity of $K_{X_{s}}+B_{s}+ L_{m}|X_{s}$ is obtained just by repeating the  above argument.  
Hence for every  $s\in S^{\circ}$, $K_{X_{s}}+B_{s} + L_{m}|X_{s}$ 
is pseudoeffective because the pseudoeffectivity is closed under specializations.   This implies that  $h_{m}$ is well defined  and is a singular hermitian metric 
on $K_{X/S} + B + L_{m}$ with semipositive curvature  over $X^{\circ}:= f^{-1}(S^{\circ})$ by Theorem \ref{NSB-P} and  induction on $m$.  And it extends to a singular hermitian metric with semipositive curvature 
just as in Sections  3.2 and \ref{dimgeq1}.   
By the semipositivity of the curvature of $h_{m}$ we see that  
 $K_{X/S} + B + L_{m}$ is pseudoeffective.

Replacing   $0= s_{0}$ by an arbitrary $s\in S^{\circ}$ and repeating the above argument,  we prove the followings 
by induction on $m$.  
\begin{enumerate}
\item[$\mbox{(A)}_{m}$] 
$K_{X/S} + B+ L_{m}$ is pseudoeffective on $X$ and 
$K_{X_{s}}+ B_{s} + L_{m}|X_{s}$ is pseudoeffective for every $s\in S^{\circ}$.
\item[$\mbox{(B)}_{m}$] 
$h_{m}|X_{s}$ is an AZD with minimal singularities on $K_{X_{s}}+B_{s} + L_{m}|X_{s}$ for every $s \in S^{\circ}$. 
\end{enumerate}
Hence  $h_{m}$ is well defined and has semipositive curvature.   Letting $m$ tend to infinity,  we see that $K_{X/S} + B$ is pseudoeffective and  
 $K_{X_{s}} + B_{s}$ is pseudoeffective for every $s\in S^{\circ}$.   \fbox{}   \vspace{3mm} \\
Now we shall complete the proof of Theorem \ref{logfamily}.  We shall use 
the same notation as above. 
In the above proof of Theorem \ref{PE}, we have seen  that 
$h_{m}|X_{s}$ is an AZD with minimal singularities on $K_{X_{s}}+B_{s} + L_{m}|X_{s}$.  
Hence again similarly as  Lemma \ref{ext3}, (2), by the induction on $m$,  we see that there exists  a positive constant $C$ independent of $m$ such that for every $m\geqq 1$, 
\begin{equation}
h_{m}|X_{s} \leqq \exp \left(C\cdot \sum_{k=0}^{m-1} \left(\frac{q-1}{q}\right)^{k}\right)\cdot h_{A,s}^{\left(\frac{q-1}{q}\right)^{m}}\cdot \hat{h}_{can,s}^{q}
\end{equation}
holds.      Letting $m$ tend to infinity,  we see that 
\begin{equation}\label{AA}
\liminf_{m\to\infty}\,\,\left(h_{A}^{-\left(\frac{q-1}{q}\right)^{m}}\!\!\!\!\cdot h_{m}\right)|X_{s}\leqq \exp(C\cdot q)\cdot \hat{h}_{can,s}^{q}
\end{equation}
holds.  Since the lower estimate of the left-hand side is obtained as 
in Section \ref{upper}, the lower semicontinuous envelope of the left-hand side is a well defined singular hermitian metric with semipositive curvature.  
We note that  
\begin{equation}\label{BB}
\hat{h}^{q}_{can}|X_{s} = O\left(\liminf_{m\to\infty}\,\,\left(h_{A}^{-\left(\frac{q-1}{q}\right)^{m}}\!\!\!\!\cdot h_{m}\right)|X_{s}\right)
\end{equation}
holds , since $\hat{h}_{can}$ is an AZD with minimal singularities on $K_{X/S}+D$ and $\liminf_{m\to\infty}\,\,\left(h_{A}^{-\left(\frac{q-1}{q}\right)^{m}}\!\!\!\!\cdot h_{m}\right)$ is a singular hermitian metric on 
$K_{X/S} + D$ with semipositive curvature. 
Combining (\ref{AA}) and (\ref{BB}), we see that 
\[
\hat{h}_{can}|X_{s} = O(\hat{h}_{can,s})
\]
holds.  This completes the proof of the asseretion (2) in Theorem \ref{logfamily}.   The rest of the proof is identical to the one of Theorem \ref{family*}. \fbox{} 

\subsection{Variation of supercanonical AZD's for relative adjoint line bundles of KLT $\mathbb{Q}$-line bundles}

We can generalize Theorem \ref{logfamily} to the case of the family of 
relative adjoint bundles of KLT singular hermitian line bundles. 
 
\begin{theorem}\label{logfamilyL}
Let $f : X \longrightarrow S$ be a proper surjective projective morphism 
between  complex manifolds with connected fibers 
and let $(L,h_{L})$ be a pseudoeffective singular hermitian $\mathbb{Q}$-line bundle on $X$  such that for a general fiber $X_{s}$, $(L,h_{L})|X_{s}$ is KLT(cf. Definition \ref{singKLT}). 
We set 
\[
S^{\circ}:= \left\{s\in S|\,\,\mbox{$f$ is smooth over $s$ and 
$(L,h_{L})|X_{s}$ is well defined and KLT}\,\,\right\}.
\]
Then there exists a singular hermitian metric $\hat{h}_{can}(L,h_{L})$ on $K_{X/S} + L$ such that 
\begin{enumerate}
\item[(1)] $\hat{h}_{can}(L,h_{L})$ has semipositive curvature  current, 
\item[(2)] $\hat{h}_{can}(L,h_{L})\!\mid\!\!X_{s}$ is an AZD on $K_{X_{s}}+L_{s}$  
(with minimal singularities) for 
every $s \in S^{\circ}$,
\item[(3)] For every $s\in S^{\circ}$, $\hat{h}_{can}(L,h_{L})\!\!\mid\!X_{s} \leqq \hat{h}_{can}((L,h_{L})|X_{s})$ 
holds, where $\hat{h}_{can}((L,h_{L})|X_{s})$ denotes the supercanonical 
AZD on $K_{X_{s}}+L_{s}$ with respect to $h_{L}|X_{s}$ (cf. Theorem \ref{KLTAZD}). 
And $\hat{h}_{can}(L,h_{L})|\!X_{s} = \hat{h}_{can}((L,h_{L})|X_{s})$ holds outside of a set of measure $0$  on $X_{s}$ for almost every $s\in S^{\circ}$. $\square$  
\end{enumerate} 
\end{theorem}

\noindent  The proof of Theorem \ref{logfamilyL} is completely parallel 
to the one of Theorem \ref{logfamily} above.  

Assume that  for every positive rational number $\varepsilon$, $L_{s}+\varepsilon A|X_{s} (s\in S^{\circ})$ is $\mathbb{Q}$-linearly equivalent to an effective $\mathbb{Q}$-divisor $D_{\varepsilon,s}$ such that 
$(X_{s},D_{\varepsilon,s})$ is KLT.   Then we may  apply \cite{b-c-h-m} 
as in Section \ref{GC} and can prove Theorem \ref{logfamilyL}  similarly 
as Theorem \ref{logfamily}. To assure this assumption, we need the following lemma.  
\begin{lemma}\label{APP}
Let $(F,h_{F})$ be a pseudoeffective singular hermitian $\mathbb{Q}$-line bundle on a smooth projective variety $M$ such that the curvature 
$\sqrt{-1}\,\Theta_{h_{F}}$ dominates a $C^{\infty}$-K\"{a}hler form on $M$.
Suppose that $(F,h_{F})$ is KLT (cf. Definition \ref{singKLT}).  
Then there exists an effective $\mathbb{Q}$-divisor $V$ on $M$ such that 
$L$ is $\mathbb{Q}$-linearly equivalent to $V$ and $(M,V)$ is KLT.  
\fbox{}
\end{lemma}   
{\em Proof}.  By the assumption and Nadel's vanishing theorem (\cite[p.561]{n}), we see that for every sufficiently large $m$ such that $mF$ is Cartier, $\mathcal{O}_{M}(mF)\otimes \mathcal{I}(h_{F}^{m})$ is globally 
generated.  Take such a sufficiently large $m$ and let $\sigma$ be a general nonzero element of 
$H^{0}(M,\mathcal{O}_{M}(mF)\otimes \mathcal{I}(h_{F}^{m}))$ and set 
$V= m^{-1}(\sigma)$, where $(\sigma)$ denotes the divisor associated with 
$\sigma$.  Then $(M,V)$ is KLT by the global generation property. \fbox{} 
\vspace{3mm} \\  
The rest of the proof of Theorem \ref{logfamilyL} is parallel to the one of 
Theorem \ref{logfamily}. Hence we omit it. In fact we just need to replace $h_{D,s}$ by $h_{L,s}: = h_{L}|X_{s}$.  \fbox{} \vspace{3mm} \\ 

\noindent The following pseudoeffectivity theorem is similar to 
\cite[Theorem 0.1]{b-p}.  The advantage  is that we deal with $\mathbb{Q}$-line bundles and without assuming the existence of sections on the special 
fiber.  But our theorem has the additional KLT assumption.    

\begin{theorem}\label{PE2}
Let $f : X \longrightarrow S$ be a proper surjective projective morphism 
between  complex manifolds with connected fibers 
and let $L$ be a $\mathbb{Q}$-line bundle on $X$  with a singular hermitian 
metric $h_{L}$ with semipositive curvature.  We assume that $S$ is quasiprojective or Stein.    
Suppose that the set:  
\[
S^{\circ}:= \{s\in S|\,\mbox{$f$ is smooth over $s$, $(L,h_{L})|X_{s}$ is well defined and KLT}\,\,\}
\]
is nonempty and  there exists  some $s_{0} \in S^{\circ}$ such that $K_{X_{s_{0}}}+L|X_{s_{0}}$ is pseudoeffective.  

Then $K_{X/S} + L$ is pseudoeffective on $X$. \fbox{}
\end{theorem}
The proof of Theorem \ref{PE2} is parallel to that of Theorem \ref{PE},
if we use the perturbation as in Section \ref{GC} and Lemma \ref{APP}.   
Hence we omit it.  The assumption that $S$ is quasiprojective or Stein is used 
to globalize Lemma \ref{APP} on $X$.  



\subsection{Proof of Theorems \ref{logplurigenera} and \ref{logplurigeneraL}}
\label{INV}
In this subsection we shall prove Theorems \ref{logplurigenera} and \ref{logplurigeneraL}. \vspace{3mm} \\   
Let $f : X \longrightarrow S$ be a proper surjective projective morphism 
between complex manifolds with connected fibers.  
Let $D$ be an effective $\mathbb{Q}$-divisor on $X$  such that 
\begin{enumerate}
\item[(a)] $D$ is $\mathbb{Q}$-linearly equivalent to a $\mathbb{Q}$-line bundle $B$, 
\item[(b)] The set:  
$S^{\circ} := \{ s\in S|\,\,\mbox{$f$ is smooth over $s$ and $(X_{s},D_{s})$ is 
KLT}\,\,\}$ is nonempty. 
\end{enumerate}
If  $K_{X_{s_{0}}}+D_{s_{0}}$ is pseudoeffective for some $s_{0}\in S^{\circ}$,  then for every $s\in S^{\circ}$, $K_{X_{s}}+D_{s}$ is pseudoeffective by 
Theorem \ref{PE}.  Hence we may and do assume for every $s\in S^{\circ}$, 
$K_{X_{s}}+D_{s}$ is pseudoeffective.  In fact otherwise we see that 
$P_{m}(X_{s},B_{s})$ is identically $0$ on $S^{\circ}$ for every $m\geqq 1$.    

Since the problem is local, to prove Theorem \ref{logplurigenera}, we may and do assume that $S$ is the unit open disk  in $\mathbb{C}$ and $S^{\circ} = S$. 
Let $\hat{h}_{can}$ be the relative supercanonical AZD on $f: (X,D) \to S$@
as in Theorem \ref{logfamily}.  Let $m$ be an arbitrary positive integer such that $mB$ is a genuine line bundle. 
Let $s\in S^{\circ}$ and let $\sigma \in H^{0}(X_{s},\mathcal{O}_{X_{s}}(m(K_{X_{s}}+B_{s})))$ be an arbitrary nonzero element.
 Since $\hat{h}_{can}|X_{s}$ is an AZD with minimal singularities (see Definition \ref{minAZD}) by Theorems \ref{main2} and  \ref{family} (or Lemma \ref{ext2}), 
\[
\hat{h}_{can}|X_{s} = O\left(\frac{1}{|\sigma|^{\frac{2}{m}}}\right)
\] 
holds.  
Hence $\hat{h}_{can}^{m}(\sigma,\sigma)$  is bounded on $X_{s}$. 
We note that 
\[
|\sigma|^{2}\cdot (\hat{h}_{can}^{m-1}|X_{s})\cdot h_{D,s} 
= |\sigma|^{2}\cdot (\hat{h}_{can}^{m}|X_{s})\cdot ((\hat{h}_{can}^{-1}|X_{s})\cdot h_{D,s})
\]
holds and $(\hat{h}_{can}^{-1}|X_{s})\cdot h_{D,s}$ is a locally integrable singular volume form on $X_{s}$, since $(X_{s},D_{s})$ is KLT.  
Hence we see that 
\[
\int_{X_{s}}|\sigma|^{2}\cdot (\hat{h}_{can}^{m-1}|X_{s})\cdot h_{D,s}
\]
is bounded.  
By Theorem \ref{family} and the $L^{2}$-extension theorem (\cite{o-t,o}), 
 we may extend $\sigma$ to an element of 
$H^{0}\left(X,\mathcal{O}_{X}(m(K_{X}+B))\right)$.  
Since $s\in S^{\circ}$ is arbitrary, noting the upper-semicontinuity theorem for cohomologies, we see that $P_{m}(X_{s},B_{s})
 = \dim H^{0}(X_{s},\mathcal{O}_{X_{s}}(m(K_{X_{s}}+B_{s})))$ is 
locally constant over $S^{\circ}$. \fbox{} \vspace{3mm} \\

\noindent The proof of Theorem \ref{logplurigeneraL} is simlar to the one 
of Theorem \ref{logplurigenera}.  Hence we omit it. 

\subsection{Semipositivity of the direct image of pluri log canonical systems}
\label{SEM}
The semipositivity of the direct image of the relative pluricanonical system has been studied 
in many papers such as \cite{f,ka1,v1,v}. 
But in the case of the relative pluri log canonical systems,  not so much is known except \cite[p.175,Theorem 1.2]{ka2}.  

Let $f : X \to S$ be a proper projective morphism between complex manifolds with connected fibers.  Let $D$ be an effective $\mathbb{Q}$-divisor on $X$  such that 
$S^{\circ}:= \{ s\in S|\,\mbox{$f$ is smooth over $s$ and $(X_{s},D_{s})$ is KLT}\,\,\}$ is nonempty.  Suppose that $D$ is $\mathbb{Q}$-linearly equivalent to a 
$\mathbb{Q}$-line bundle $B$.  Let $m$ be a positive integer such that $mB$ 
is a genuine line bundle.   
Then by Theorem \ref{logplurigenera}, the direct image:  
\begin{equation}
F_{m}:= f_{*}\mathcal{O}_{X}(m(K_{X/S}+B)) 
\end{equation}
is locally free over $S^{\circ}$.   
By Theorem \ref{logfamily}, the relative supercanonical AZD  $\hat{h}_{can}$ exists on $K_{X/S}+B$ and has semipositive curvature current. We define the metric $h_{m}$ on $F_{m}|S^{\circ}$ by  
\begin{equation}
h_{m}(\sigma,\sigma^{\prime}):= (\sqrt{-1})^{n^{2}}\int_{X_{s}}\hat{h}_{can}^{m-1}\cdot h_{D}\cdot\sigma\wedge\overline{\sigma^{\prime}} \hspace{10mm}(\sigma,\sigma^{\prime}\in F_{m,s}), 
\end{equation}
where $h_{D}$ is the metric defined as (\ref{hD}) and 
$n:= \dim X -\dim S$. 
Then by Theorem \ref{logfamily} and \cite[Theorem 3.5]{b-p}, we have the following theorem.

\begin{theorem}\label{direct} 
The locally bounded metric $h_{m}$ on $F_{m}|S^{\circ}$ is semipositive in the sense $h_{m}$ gives a singular hermitian metric with semipositive curvature on the tautological line bundle 
$\mathcal{O}(1)$ on $\mathbb{P}(F^{*}_{m}|S^{\circ})$.  \fbox{} 
\end{theorem}
\begin{remark}\label{extdeg}
It is trivial to generalize Theorem \ref{direct} in the case of 
the direct image of the multi adjoint line bundle of a generically KLT line bundles. 

If Conjecture \ref{delta*} is true and $\dim S = 1$, it is not difficult to see that 
$h_{m}$ gives a singular hermitian metric with semipositive curvature on the tautological line bundle 
$\mathcal{O}(1)$ on the whole $\mathbb{P}(F^{*}_{m})$. \fbox{}
\end{remark}
For the different treatments such as  weak semistability of the direct images of 
pluri log canonical systems, see  \cite{canonical,ricci, tu10}.  
In these papers, the canonical metric comes from a log canonical bundle on 
the base space of an Iitaka fibrations and the construction of the metric is 
quite different from the one here. 
   
\section{Appendix}

Here we collect  miscellaneous facts.    
   
\subsection{Choice of the sufficiently ample line bundle $A$}\label{suff} 
 In this subsection, we shall prove the following proposition. 
\begin{proposition}\label{suf} 
Let $X$ be a smooth projective $n$-fold. 
Then there exists an ample line bundle $A$ on $X$ 
such that for every pseudoeffective singular hermitian line bundle $(L,h_{L})$
on $X$,
${\cal O}_{X}(A + L)\otimes {\cal I}(h_{L})$ and 
 ${\cal O}_{X}(K_{X}+ A + L)\otimes {\cal I}(h_{L})$ are globally generated. \fbox{}
\end{proposition}
{\em Proof}. 
We construct such an  $A$ by using  $L^{2}$-estimates. 
In fact let $g$ be a K\"{a}hler metric on $X$ and for every $x$, $d_{x}$ 
denotes the distance function from $x$ and let $R > 0$ denotes the infimum of 
the injective radius on $(X,g)$.  Let $\rho$ be a $C^{\infty}$-function on 
$[0,R)$ such that 
\begin{enumerate}
\item[(1)] $0\leqq \rho \leqq 1$,
\item[(2)] $\mbox{Supp}\,\rho \subset [0,\frac{2}{3}R]$, 
\item[(3)] $\rho \equiv 1$ on $[0,\frac{1}{3}R]$. 
\end{enumerate} 
Then we may take an ample line bundle $A$ and a $C^{\infty}$-hermitian metric 
$h_{A}$ such that \\
 $\sqrt{-1}\left(\Theta_{h_{A}} + 2n\partial\bar{\partial}\left(\rho (d_{x})\cdot \log d_{x}\right)\right)$ 
and 
 $\mbox{Ric}_{g} + \sqrt{-1}\left(\Theta_{h_{A}} + 2n\partial\bar{\partial}\left(\rho (d_{x})\cdot \log d_{x}\right)\right)$ 
are closed strictly positive $(1,1)$ current on $X$ for every $x\in X$. 
Then by Nadel's vanishing theorem \cite[p.561]{n}, 
for every pseudoeffective singular hermitian line bundle $(L,h_{L})$
on $X$,
${\cal O}_{X}(A + L)\otimes {\cal I}(h_{L})$ and 
 ${\cal O}_{X}(K_{X}+ A + L)\otimes {\cal I}(h_{L})$ are globally generated. 
\fbox{} \vspace{3mm} \\
 
Let $h_{A}$ be a a $C^{\infty}$-hermitian metric on $A$ 
with strictly positive curvature as above.   
Let us fix a $C^{\infty}$-volume form $dV$ on $X$. 
By the $L^{2}$-extension theorem (\cite{o}) we take  
a sufficiently ample line bundle $A$ so that for every $x\in X$ and 
for every pseudoeffective singular hermitian line bundle $(L,h_{L})$, there exists a bounded interpolation operator: 
\[
I_{x} : A^{2}(x,(A + L)_{x},h_{A}\cdot h_{L},\delta_{x})
\rightarrow A^{2}(X,A + L,h_{A}\cdot h_{L},dV)
\]
such that the operator norm of $I_{x}$ is bounded by 
a positive constant independent of $x$ and $(L,h_{L})$, 
where $A^{2}(X,A + L,h_{A}\cdot h_{L},dV)$ denotes the Hilbert space defined by
\[
 A^{2}(X,A + L,h_{A}\cdot h_{L},dV) := \left\{ \sigma\in \Gamma (X,{\cal O}_{X}(A + L)\otimes {\cal I}(h_{L}))| \int_{X}\mid\sigma\mid^{2}\cdot h_{A}\cdot h_{L}\cdot dV < + \infty\right\} 
\]
with the $L^{2}$-inner product:  
\[
(\sigma ,\sigma^{\prime}):= \int_{X}\sigma\cdot\overline{\sigma^{\prime}}\cdot h_{A}\cdot h_{L}\cdot dV 
\]
and  $A^{2}(x,(A + L)_{x},h_{A}\cdot h_{L},\delta_{x})$ is defined 
similarly, where $\delta_{x}$ is the Dirac measure supported at $x$.
We note that if $h_{L}(x) = + \infty$, then 
$A^{2}(x,(A + L)_{x},h_{A}\cdot h_{L},\delta_{x}) = 0$.  

\subsection{Analytic Zariski decompositions and singular hermitian metrics with minimal singularities}\label{appendix}

In this paper we have used the notion of AZD's (cf. Definition \ref{azd}). 
We note that there is a similar but {\bf different} notion : singular hermitian metrics with minimal singularities introduced in \cite{d-p-s} (see Definiton 
\ref{minAZD} below).   I would like to explain the difference of these two notions here. 

According to \cite{d-p-s}, an AZD is constructed 
for any pseudoeffective line bundle $L$ as follows. 
Let $h_{L}$ be any $C^{\infty}$-hermitian metric on $L$. 
Let $h_{0}$ be an AZD on $K_{X}$ defined by the lower envelope of : 
\[
 \inf \left\{ h \mid \mbox{$h$ is a singular hermitian metric on $L$ with  $\sqrt{-1}\,\Theta_{h}\geqq 0$,\,$h \geqq h_{L}$}\right\}, 
\]
where the $\inf$ denotes the pointwise infimum. 
This construction is exactly the same as (\ref{hzero}) above.  
Then by the classical theorem of Lelong (\cite[p.26, Theorem 5]{l}) it is easy to verify that $h_{0}$ is an AZD on $L$ (cf. \cite[Theorem 1.5]{d-p-s}).
By the definition, $h_{0}$ is of  minimal singularities in the following sense. 
\begin{definition}\label{minAZD}
Let $L$ be a pseudoeffective line bundle on a smooth projective variety $X$.
An AZD $h$ on $L$  is said to be  {\bf a singular hemitian metric with minimal singularities} or {\bf an AZD with minimal singularities}, if for any singular hermitian metric $h^{\prime}$ on $L$ with semipositive curvature current, there exists a positive constant $C$ such that
\[
h \leqq   C \cdot h^{\prime}
\]
holds on $X$.  In particular for any AZD $h^{\prime}$ on $L$ the above inequality holds for some positive constant $C$. \fbox{} 
\end{definition}
We note that any AZD's  with minimal singularities are quasi-isometric, i.e., 
any two AZD's with minimal singularities $h_{1},h_{2}$ on a common line bundle 
$L$, there exists a positive constant $C > 1$ such that 
\[
C^{-1}\cdot h_{2} \leqq h_{1} \leqq C\cdot h_{2}
\]
holds.  In particular for any AZD with minimal singularities $h$ on a line bundle $L$, 
the multiplier ideal $\mathcal{I}(h^{m})$ is uniquely determined for every $m$. And the above construction of an AZD is very easy.  In the above sense, 
the AZD with minimal singularities is very canonical.   \vspace{3mm} \\
\noindent But in general, an AZD is not with minimal singularities as follows.  
\begin{example}
Let $X$ be a smooth projective variety and let $D$ be a divisor with simple normal crossings on $X$.   Suppose that $K_{X} + D$ is ample. 
Then there exists a complete K\"{a}hler-Einstein form $\omega_{E}$ on $X\backslash D$ with $-\mbox{\em Ric}_{\omega_{E}} = \omega_{E}$ and $\omega_{E}$ extends to a closed positive current on $X$ 
with vanishing Lelong numbers and $[\omega_{E}] = 2\pi c_{1}(K_{X}+D)$ (\cite{ko}).  
The metric $h:= (\omega_{E}^{n})^{-1} (n = \dim X)$ 
is a singular hermitian metric on $K_{X}+ D$ with strictly positive 
curvature on $X$.  Let $D = \sum_{i}D_{i}$ be the irreducible decomposition 
of $D$ and let $\sigma_{i}$ be a nontrivial global section of $\mathcal{O}_{X}(D_{i})$ with divisor $D_{i}$. 
$h$ is an AZD on $K_{X} + D$, but $h$ has logarithmic singularities along $D$, i.e., there exists a $C^{\infty}$-hermitian metric on $h_{0}$ 
on $K_{X}+ D$ such that 
\[
h = h_{0}\cdot\prod_{i}|\log \parallel\sigma_{i}\parallel\!|^{2}, 
\]
where $\parallel\sigma_{i}\parallel$ denotes the hermitian norm of $\sigma_{i}$
with respect to a $C^{\infty}$-hermitian metric on $\mathcal{O}_{X}(D_{i})$
respectively.
Hence $h$ blows up along $D$.  In particular $h$ is not of minimal singularities.\fbox{}
\end{example}
As above,  even in the case of ample line bundles, some natural AZD's are not 
of minimal singularities.   Indeed the notion of AZD's is much broader than the notion of singular hermitian metrics with minimal singularities.  Much more general singular K\"{a}hler-Einstein metrics on  LC pairs of log general type was considered in \cite{ricci}.  More precisely in the paper, we have considered  singular K\"{a}hler-Einstein metrics on LC pairs $(X,D)$ of log general type such that the inverse of the K\"{a}hler-Einstein volume form is an AZD on  $K_{X}+D$.  In that case the AZD is not 
necessarily of minimal singularities as is seen in the above examples.

And also it is not clear whether the canonical AZD $h_{can}$ defined in 
Section \ref{hcan} has  minimal 
singularities.      

The above examples indicate us that we had better not to restrict ourselves to consider  AZD's with 
minimal singularities to consider broader canonical singular hermitian 
metrics. 

\subsection{An alternative proof of Theorem \ref{logfamily}}\label{T2}

In this subsection, we shall give an alternative proof of Theorem \ref{logfamily} by using 
the argument in \cite{e-p}.  Here we assume the results in Sections \ref{LGC} and \ref{GC}.  The reason why we present an alternative proof is that although 
the proof itself is far more complicated than the one in Section \ref{Dynamical}, it may indicate the way how to handle the extension without assuming the 
bigness.  Hence it may have an independent interest. 

The strategy of the proof is as follows.  
We follows the argument in \cite{e-p} when the fibers over $S^{\circ}$ is 
of log general type.  The key point of the proof is 
we subdivide the extension into several steps by  using  the 
logarithmic vanishing theorem as \cite[Theorems 2.9, 3.2]{e-p} for the 
extension of holomorphic sections similar to \cite[Propositions 4.1 and 4,2]{e-p}. But the theorem requires the bigness of the line bundle on every component 
of the log canonical centers.  In our case this condition need not be satisfied. Hence we perturb the log canonical bundle by adding ample $\mathbb{Q}$-line bundles as in Section \ref{GC}.  
We note that this condition is stated in \cite{e-p} by the language of augumented base locus $B_{+}$ (cf. \cite{e-p}).  Then as in Section \ref{GC}, we take the 
limit.   Since we have already known how to transform the inclusion of 
the multiplier ideals into the estimate of the canonical singular hermitian metric as Lemmas \ref{ext2} and \ref{ext3}, this part is essentially nothing new.   Hence the  proof here is just a combination of the perturbation and the argument in  
\cite{e-p} using the estimate of canonical singular hermitian metrics in 
Sections \ref{LGC} and \ref{GC}. Since we follow the argument in \cite{e-p}, 
we do not repeat the proof, when we just borrow the argument 
in \cite{e-p}. 

 Let us assume that $D$ is $\mathbb{Q}$-linearly equivalent to a $\mathbb{Q}$-line bundle $B$. As for the proof of the assertion (1) of Theorem \ref{logfamily}, nothing changes and it follows from Theorem \ref{NSB-P}. 
Hence we only need to verify the assertions (2) and (3) of Theorem \ref{logfamily}.  For simplicity we shall consider the case: $\dim S = 1$. 
To prove the assertion (3), we may and do  assume that $S =  S^{\circ} = \Delta$ hold.  
 
\begin{lemma}\label{COMP} For every $s\in S^{\circ}$, there exists a positive constant $C_{+}$ depending on $s$ such that 
\[
\hat{h}_{can}|X_{s} \leqq C_{+}\cdot \hat{h}_{can,s}
\]
holds on $X_{s}$. \fbox{}
\end{lemma}
{\em Proof}. 
Let $U$ be an open neighborhood of $s$ in $S^{\circ}$ which is biholomorphic 
to the unit open polydisk $\Delta^{k}(k = \dim S)$ as above. 
We shall identify $K_{X/S}|f^{-1}(U)$ with $K_{X}$ by 
\begin{equation}
\otimes f^{*}dt : 
K_{X/S} \longrightarrow  K_{X}, 
\end{equation}
where $t$ is the standard coordinate 
on $\Delta$.  Let $A$ be an ample $\mathbb{Q}$-line bundle on $X$.  By \cite{b-c-h-m} the relative log canonical 
ring $\oplus_{m=0}^{\infty}f_{*}\mathcal{O}_{S}(\lfloor m(K_{X/S}+ A + D)\rfloor)$ 
is locally finitely generated.     
Let $q$ be a positive integer such that $q(K_{X/S}+A+ D)$ is integral. 
In this case the relative log canonical ring is big at $s$, i.e., 
the image of  
\begin{equation}\label{GR}
\oplus_{m=0}^{\infty}f_{*}\mathcal{O}_{X}(mq(K_{X/S}+A+ D))_{s} 
\to \oplus_{m=0}^{\infty}H^{0}(X_{s},\mathcal{O}_{X_{s}}(mq(K_{X_{s}}+A_{s}+D_{s}))) 
\end{equation} 
is a subring of maximal growth. 
Let $\mu_{s} : Y_{s}\to X_{s}$ be a modification such that Zariski decomposition of the image of (\ref{GR}): 
\begin{equation} 
\mu_{s}^{*}(K_{X_{s}}+A_{s} + B_{s}) = P_{s} + N_{s}\,\,\,\,(P_{s},N_{s}\in \mbox{Div}(Y_{s})\otimes\mathbb{Q})
\end{equation} 
exists as (\ref{ZD}), i.e., $P_{s}$ is nef, $N_{s}$ is effective and $H^{0}(Y_{s},\mathcal{O}_{Y_{s}}(\lfloor mqP_{s}\rfloor))$ is isomorphic to the image of 
$f_{*}\mathcal{O}_{X}(mq(K_{X/S}+ D))_{s} 
\to H^{0}(X_{s},\mathcal{O}_{X_{s}}(mq(K_{X_{s}}+B_{s})))$.  In this case 
$P_{s}$ is semiample by the finite generation of the relative log canonical ring. Since $P_{s}$ is semiample, and $(X_{s},D_{s})$ is KLT, 
again by \cite{b-c-h-m}, there exists a Zariski decomposition 
\begin{equation}\label{ZDD}
\mu_{s}^{*}(K_{X_{s}}+ A + B_{s})= Q_{s} + E_{s}\,\,\,\, (Q_{s},E_{s}\in \mbox{Div}(Y_{s})
\otimes \mathbb{Q})
\end{equation}
of $\mu_{s}^{*}(K_{X_{s}}+B_{s})$ such that $Q_{s}$ is semiample.
Hereafter for simplicity, we shall assume that $D_{red} + X_{s}$ is a simple normal crossing divisor.  The general case is handled by taking a log resolution of $(X,D + X_{s})$ as in \cite[Proposition 5.4]{e-p}.  We shall review the notion of adjoint ideals in \cite{e-p}.    

\begin{definition}(Adjoint ideals (\cite{e-p}))\label{ADJ} 
Let $\Gamma$ be a reduced simple normal crossing divisor on $X$ and 
$\mathfrak{a}\subset \mathcal{O}_{X}$ an ideal sheaf such that no log-canonical center of $\Gamma$ is contained in $Z(\mathfrak{a})$.  Let $f : Y \to X$ 
be a common log-resolution for the pair $(X,\Gamma)$ and the ideal $\mathfrak{a}$, ad write $\mathfrak{a}\cdot\mathcal{O}_{Y} = \mathcal{O}_{Y}(-E)$.  
We set 
\begin{equation}\label{adj}
\mbox{\em Adj}_{\Gamma}(X,\mathfrak{a}^{\lambda}):= 
f_{*}\mathcal{O}_{Y}(K_{Y/X}- f^{*}\Gamma + 
\sum_{\mbox{\em ld}(\Gamma,D_{i}) = 0}D_{i} -\lfloor \lambda E\rfloor),
\end{equation}
where the sum appearing in the expression is taken over all divisors on $Y$ having log-discrepancy $0$ with respect to $\Gamma$, i.e., among those appearing in  $\Gamma^{\prime}$ in the expression $K_{Y}+\Gamma^{\prime} = f^{*}(K_{X}+\Gamma)$. 
We note that 
\[
\mbox{\em Adj}_{\Gamma}(X,\mathfrak{a}^{\lambda}) \subset \mathcal{I}(X,\mathfrak{a}^{\lambda}).
\]
For a graded system $\mathfrak{a}_{*} = \{ \mathfrak{a}_{m}\}$ of ideal sheaves,we set 
\[
\mbox{\em Adj}_{\Gamma}(X,\mathfrak{a}_{*}^{\lambda})
:= \mbox{\em Adj}_{\Gamma}(X,\mathfrak{a}_{m}^{\lambda/m})
\]
for $m$ sufficiently divisible.  
For a $\mathbb{Q}$-effective line bundle $L$, we set 
\[
\mbox{\em Adj}_{\Gamma}(X,\parallel L\parallel)
:= \mbox{\em Adj}_{\Gamma}(X,\mathfrak{b}_{*}),
\]
where $\mathfrak{b}_{*}$ denotes the graded system of the base ideal 
of $L$.  We call $\mbox{\em Adj}_{\Gamma}(X,\parallel L\parallel)$ the asymptotic adjoint ideal of $L$ with respect to $\Gamma$.  These definitions can be generalized to the case of 
a pair $(X,\Lambda)$ of $X$ and an effective $\mathbb{Q}$-divisor $\Lambda$ provided $\mbox{\bf B}(L)\cup \mbox{\em Supp}\,\Lambda$ does not contain any LC center 
of $\Gamma$, where $\mbox{\bf B}(L)$ denotes the stable base locus of $L$.  \fbox{}
\end{definition}
Let us interpret Definition \ref{ADJ} in terms of singular hermitian metrics. 
To do this we introduce the following notion. 
\begin{definition}\label{algsing}
Let $(L,h_{L})$ be a singular hermitian line bundle on a complex manifold 
$W$.   We say that $h_{L}$ is of algebraic singularities, 
if $h_{L}$ is written locally as   
\[
h_{L} = h_{0}\cdot\left(\sum_{i=1}^{N} |f_{i}|^{2}\right)^{-\alpha},
\]
where $f_{1},\cdots ,f_{N}$ are local holomorphic functions, 
$h_{0}$ is a local $C^{\infty}$-hermitian metric on $L$ and  $\alpha$ is a 
positive number. \fbox{}
\end{definition}
Let $(L,h_{L})$ be a singular hermitian line bundle of algebraic singularities.
Then there exists a modification $\mu : V \to W$ such that
\begin{enumerate}
\item[(1)] The exceptional divisor of $\mu$ is a simple normal crossing divisor, 
\item[(2)]  There exists an effective $\mathbb{R}$-divisor $E$ on $V$ such that 
$\mbox{Supp}\, E$ is a simple normal crossing divisor on $V$ and 
$\mathcal{I}(h_{L}^{m}) = \mu_{*}\mathcal{O}_{V}(K_{V/W} -\lfloor mE\rfloor)$ holds for every positive integer $m$. 
\end{enumerate}
Let $\Gamma$ be a reduced simple normal crossing divisor on $W$.  
By using this divisor $E$, we can define the adjoint ideal 
$\mbox{Adj}_{\Gamma}(W;h_{L})$ as (\ref{adj}), if the singular set $\mbox{Sing}\, h_{L}$ of $h_{L}$ does not contain 
any components of $\Gamma$.  This notion is used to rewrite the argument in
 \cite{e-p} in terms of singular hermitian metrics and to avoid the 
 use of asymptotic multiplier ideals in \cite{e-p}.  For a pseudoeffective 
 singular hermitian line bundle $(L,h_{L})$ on a smooth projective variety $W$,  we say that $(L,h_{L})$ is big, if $\,\,\limsup_{m\to\infty} m^{-\dim W}\cdot h^{0}(W,\mathcal{O}_{W}(mL)\otimes \mathcal{I}(h_{L}^{m}))$ is positive.         
The following proposition follows from the same argument as in 
the proof of \cite[Theorems 2.9 and 3.2]{e-p}.  This is nothing but the singular hermitian version of Norimatsu's vanishing  theorem (\cite{no}).   
\begin{proposition}\label{logvanishing}
Let $X$ be a smooth projective variety and let $\Gamma$ be a reduced 
simple normal crossing divisor on $X$. 
Let $(L,h_{L})$ be a singular hermitian line bundle on $X$ such that 
\begin{enumerate}
\item[(1)] $h_{L}$ is of algebraic singularities, 
\item[(2)] $(L,h_{L})$ is big and the restriction $(L,h_{L})|\Gamma_{j}$ to every irrducible component  $\Gamma_{j}$ of $\Gamma$ is well defined and big.
\end{enumerate}  
Then 
\[
H^{q}(X,\mathcal{O}_{X}(K_{X}+ L + \Gamma)\otimes \mbox{\em Adj}_{\Gamma}(X,h_{L})) = 0
\]
holds for every $q \geqq 1$. \fbox{}
\end{proposition}
The following extension theorem is similar to \cite[Proposition 4.2]{e-p}.
The proof follows from Proposition \ref{logvanishing}.  
\begin{proposition}\label{logextension}
Let  $X$ be a smooth projective variety and $S \subset X$ be a smooth divisor. 
Let $\Gamma$ be an effective integral divisor on $X$ such that $S + \Gamma$ is a reduced simple normal crossing divisor.  Let $\Gamma = \sum \Gamma_{j}$ be 
the irreducible decomposition. Let $(L,h_{L})$ be a big pseudoeffetive line bundle  on $X$ with algebraic singularities, 
such that no log-canonical center of $(X,S+\Gamma)$ is contained in 
$\mbox{\em Sing}\,h_{L}\cup \mbox{\em Supp}(\Lambda)$ and $(L,h_{L})|\Gamma_{j}$ is big for every $j$.  If $A$ is an integral nef divisor on $X$, then the sections in 
\[
H^{0}(S,\mathcal{O}_{S}(K_{S}+A_{S}+\Gamma_{S}+L_{S})\otimes\mbox{\em Adj}_{\Gamma_{S}}(S,h_{L}|S)))
\]
are in the image of the restriction:  
\[
H^{0}(X,\mathcal{O}_{X}(K_{X}+S+A+\Gamma+L))
\to 
H^{0}(S,\mathcal{O}_{S}(K_{S}+A_{S}+\Gamma_{S}+L_{S})).
\]
\fbox{}
\end{proposition}

First we shall assume that $Q_{s}$ is big.
We shall replace $D_{s}$ by 
\begin{equation}\label{adjustment}
(D_{s} - \mu_{s,*}E_{s})_{\geqq 0}
\end{equation}
where for a Weil divisor $F = \sum a_{i}F_{i}$,  $F_{\geqq 0}:= \sum \max (a
_{i},0)F_{i}$, $\mu_{s} : Y_{s} \to X_{s}$ be the morphism as in (\ref{ZDD}) and $E_{s}$ is the effective $\mathbb{Q}$-divisor in 
(\ref{ZDD}). And we shall change the coefficients  of $D$ so that  $D|X_{s} = D_{s}$ holds.

Let $A$ be an ample $\mathbb{Q}$-divisor on $X$. 
Let $k$ be a sufficiently divisible positive integer such that $kD$  and 
$kA$  are integral divisors.   
Now we proceed as \cite{e-p}.  First we shall decompose $kD$ as  
\begin{equation}
kD = \Delta_{1} + \cdots + \Delta_{k-1}
\end{equation}
such that each $\Delta_{i}$ is a reduced simple normal crossing divisor. 
This is possible, since $(X,D)$ is KLT and $\mbox{Supp}\, D$ is a simple normal crossing divisor.  
We set
\begin{equation}
M:=  kK_{X}+ kA + \Delta_{1} + \cdots +\Delta_{k-1}. 
\end{equation} 
Let $H$ be a sufficiently ample line bundle on $X$ which such that 
for any pseudoeffective singular hermitian line bunlde $(L,h_{L})$ on $X_{s}$, $\mathcal{O}_{X_{s}}(H+L+ (\ell+1) K_{X}+ kA+ \Delta_{1}+ \cdots + \Delta_{\ell})\otimes \mathcal{I}(h_{L})$ is globally generated for every $0\leqq \ell \leqq k-1$.   For $0\leqq \ell \leqq k-1$, we extend every element of 
\begin{equation}
H^{0}(X_{s},\mathcal{O}_{X_{s}}(H+ mM  +(\ell+1)(K_{X}+X_{s})+ kA + \Delta_{1}+ \cdots + \Delta_{\ell})
\otimes \mathcal{I}(h_{Q}^{km}))
\end{equation}
to an element of 
\begin{equation}
H^{0}(X,\mathcal{O}_{X}(H+mM +(\ell+1)K_{X}+ kA+ \Delta_{1}+ \cdots + \Delta_{\ell}))
\end{equation}
by induction on $m$ and $\ell$.  The assertion is trivial when $m= \ell = 0$. 
Suppose that we have already costructed the extension for some $m$ and $\ell -1$. 
Then we take a basis of $\{\sigma^{(m)}_{\ell,1},\cdots ,\sigma^{(m)}_{\ell,N(m,\ell)}\}$ of  
\begin{equation}
H^{0}(X_{s},\mathcal{O}_{X_{s}}(H+ mM  +\ell(K_{X}+X_{s})+ kA + \Delta_{1}+ \cdots + \Delta_{\ell})
\otimes \mathcal{I}(h_{Q}^{km}))
\end{equation}
and extend the basis to a set of sections 
$\{\tilde{\sigma}^{(m)}_{\ell,1},\cdots ,\tilde{\sigma}^{(m)}_{\ell,N(m,\ell)}\}$
in
\begin{equation}
H^{0}(X,\mathcal{O}_{X}(H+mM +\ell K_{X}+ kA+ \Delta_{1}+ \cdots + \Delta_{\ell}))
\end{equation}
by Proposition \ref{logextension}.    
We define the singular hermitian metric $h_{m,\ell}$ on $mM + H + \ell K_{X}+\Delta_{1}+\cdots +\Delta_{\ell}$ by 
\begin{equation}
h_{m,\ell} := \frac{1}{\sum_{j=1}^{N(m,\ell)} |\tilde{\sigma}^{(m)}_{\ell,j}|^{2}}
\end{equation}
which is apparently of algebraic singularities and with semipositive curvature. 
As in the proof of \cite[Proposition 5.4]{e-p}, by the induction on $\ell$, if we take $H$ to be sufficiently ample, by the induction on $\ell$, we see that 
\begin{equation}
\mathcal{I}(h_{Q}^{m}) \subseteq 
\mbox{Adj}_{\Delta_{\ell}}(X_{s};h_{m,\ell}|X_{s})
\end{equation}
holds for every  $0\leqq \ell \leqq k-1$.  We note that $\mbox{Adj}_{\Delta_{\ell}}(X_{s};h_{m,\ell}|X_{s})$ 
is well defined, since  by the choice of $H$,  
$\mathcal{O}_{X_{s}}((H+mM +(\ell+1) K_{X}+ kA+ \Delta_{1}+ \cdots + \Delta_{\ell})|X_{s})\otimes \mathcal{I}(h_{Q}^{km})$ is globally generated and the adjustment (\ref{adjustment}) implies   $\mbox{Bs}|((H+mM +(\ell+1) K_{X}+ kA+ \Delta_{1}+ \cdots + \Delta_{\ell})|X_{s})\otimes \mathcal{I}(h_{Q}^{km})|$ does not contain any irreducible component of $\Gamma$.  Moreover $|((H+mM +(\ell+1) K_{X}+ kA+ \Delta_{1}+ \cdots + \Delta_{\ell})|X_{s})\otimes \mathcal{I}(h_{Q}^{km})|$ defines a birational map   
from each component of $\Gamma$ into a projective space by the effect of
the ample line bundles $H$ and $A$\footnote{We note that in \cite{e-p}, to carry out the similar argument, they have assumed that the restricted base locus $B_{-}(M)$ does not contain any closed subset $W$ with minimal log discrepancy (cf. Definition \ref{KLT}) at its generic point 
 $\mbox{mld}(\mu_{W};X,X_{s}+ D) < 1$, which intersects $X_{s}$ but  is different from $X_{s}$ itself, because they have used the asymptotic adjointideals.}.      
Hence we may apply Proposition \ref{logextension}.
This completes the inductive construction of the metrics $\{ h_{m,\ell}\}_{\ell=0}^{k-1}$. Then we set $h_{m+1,0} = h_{m,k-1}$ and continue the induction.  
This completes the construction of the metrics $\{ h_{m,\ell}\}_{m,\ell\geqq 0}$.  

Let $h_{A}$ be a $C^{\infty}$-hermitian metric on $A$ with strictly positive curvature.    Let $\hat{h}_{can,D}(A,h_{A})$ (resp. $\hat{h}_{can,D,s}((A,h_{A})|X_{s})$)  be the supercanonical AZD on $K_{X/S}+A+D$ (resp. $K_{X_{s}}+A_{s}+D_{s}$)
with respect to $h_{A}$  similar to Theorem \ref{mainB} in  Section \ref{T1}.

Hence by the similar argument as in the proof of Lemma \ref{ext2} in Section \ref{LGC}, using H\"{o}lder's inequality, the existence of the sequece of metrics  $\{ h_{m,\ell}\}_{m,\ell\geqq 0}$ 
implies that 
\begin{equation}
\hat{h}_{can,D}(A,h_{A})|X_{s} = O(\hat{h}_{can,D,s}(A,h_{A})|X_{s})
\end{equation}
holds.    Replacing $(A,h_{A})$ by $(\epsilon A,h_{A}^{\epsilon})(\epsilon \in \mathbb{Q}^{+})$ and letting $\epsilon$ tend to $0$, by the argument in Section \ref{GC}, 
we completes the proof of Lemma \ref{COMP}. 
\fbox{} \vspace{3mm} \\

\noindent The rest of the proof is similar as the one of Theorme \ref{family*}.
This completes the proof of Theorem \ref{logfamily}. 
 \fbox{}
\small{

}

\noindent Author's address\\
Hajime Tsuji\\
Department of Mathematics\\
Sophia University\\
7-1 Kioicho, Chiyoda-ku 102-8554\\
Japan \\
h-tsuji@h03.itscom.net 
\end{document}